\documentclass[12pt,psfig]{article}
\usepackage{amsmath, amsthm, amssymb,amstext}
\usepackage{graphicx,epsfig,subfigmat,epstopdf}
\usepackage{amsfonts}
\usepackage{amscd}
\usepackage{mathrsfs}
\usepackage{enumerate,hyperref}
\usepackage{enumerate,mdwlist}
\usepackage{latexsym}
\usepackage[all]{xy}
\usepackage{cite,xcolor}
\usepackage[ruled]{algorithm2e}

\newcommand{\E}{\mathscr{E}}
\newcommand{\B}{\mathscr{B}}
\newcommand{\D}{\mathscr{D}}
\newcommand{\R}{\mathcal{R}}
\newcommand{\RM}{\mathscr{R}}
\newcommand{\M}{\mathcal{M}}
\newcommand{\Itr}{{\rm Itr}}
\newcommand{\p}{\mathscr{P}}
\newcommand{\s}{\mathscr{S}}
\newcommand{\eg}{{\rm e.g., }}
\newcommand{\ie}{{\rm i.e., }}
\newcommand{\Rem}{{\rm Rem}}

\newtheorem{thm}{Theorem}[section]

\newtheorem{lem}[thm]{Lemma}

\theoremstyle{definition}
\newtheorem{defn}[thm]{Definition}
\newtheorem{exm}[thm]{Example}
\newtheorem{rem}[thm]{Remark}

\numberwithin{equation}{section}
\newtheorem{proc}[thm]{Procedure}
\numberwithin{definition}{subsection}
\def\be {\begin{equation}}
\def\ee {\end{equation} }
\def\ba {\begin{eqnarray} }
\def\ea {\end{eqnarray} }
\def\bes {\begin{equation*} }
\def\ees {\end{equation*} }
\def\bas {\begin{eqnarray*} }
\def\eas {\end{eqnarray*} }
\def\bpr {\begin{proof} }
\def\epr {\end{proof} }
\def \Gr {Gr\"obner }
\def \lex {{\rm lex}}
\def \alex {{\rm alex}}

\def \Itr {{\rm Itr }}
\def \la {\langle }
\def \ra {\rangle }
\def \Ri {\right }
\def \Le {\left }
\def \Terms {{\rm Terms}}

\def\LM{{\mathrm{LM}}}
\def\LC{{\mathrm{LC}}}
\def\LT{{\mathrm{LT}}}
\def\LCM{{\mathrm{LCM}}}

\makeatletter \oddsidemargin.2375in \evensidemargin \oddsidemargin
\def\be {\begin{equation}}
\def\ee {\end{equation} }
\def\ba {\begin{eqnarray} }
\def\ea {\end{eqnarray} }
\def\bes {\begin{equation*} }
\def\ees {\end{equation*} }
\def\bas {\begin{eqnarray*} }
\def\eas {\end{eqnarray*} }
\def\bpr {\begin{proof} }
\def\epr {\end{proof} }
\def \Gr {Gr\"obner }
\def \lex {{\rm lex}}
\def \alex {{\rm alex}}

\def \Itr {{\rm Itr }}
\def \la {\langle }
\def \ra {\rangle }
\def \Ri {\right }
\def \Le {\left }
\def \Terms {{\rm Terms}}

\begin{document}
\baselineskip=18pt
\renewcommand {\thefootnote}{\dag}
\renewcommand {\thefootnote}{\ddag}
\renewcommand {\thefootnote}{ }

\pagestyle{empty}

\begin{center}
                \leftline{}
                \vspace{-0.0 in}
{\Large \bf Symbolic local bifurcation analysis of scalar smooth maps
} \\ [0.1in]

{\large Majid Gazor$^{*}$\footnote{$^*\,$Corresponding author. Phone: (98-31) 33913634; Fax: (98-31) 33912602; Email: mgazor@cc.iut.ac.ir; mahsa.kazemi@math.iut.ac.ir. } and Mahsa Kazemi}

\vspace{0.1in}

{\small {\em Department of Mathematical Sciences, Isfahan University of Technology,
\\[-0.5ex] Isfahan 84156-83111, Iran }
}

\end{center}


\noindent

\vspace{0.1in}

\noindent

\begin{abstract}
The local zero structure of a smooth map may {\it qualitatively} change, when the map is subjected to small perturbations. The changes may include births and/or deaths of zeros. The {\it qualitative properties} are defined as the invariances of an appropriate {\it equivalence relation}. The occurrence of
a {\it qualitative change} in the zero structures is called a {\it bifurcation} and the map is named a {\it singularity}. The local bifurcation
analysis of singularities has been extensively studied in {\it singularity theory} and many powerful algebraic tools have been developed for their study.
However, there does not exist any symbolic computer-library for this purpose.  We suitably generalize some powerful tools from algebraic geometry for correct implementation of the results from singularity theory. We provide some required criteria along with rigorous proofs for efficient and cognitive computer-implementation.
We have accordingly developed a {\sc Maple} end-user friendly library, named ``\texttt{Singularity}'', for an efficient and complete local bifurcation analysis of real zeros of scalar smooth maps. We have further written a comprehensive user-guide for \texttt{Singularity}. The main features of \texttt{Singularity} are briefly illustrated along with a few examples. 
\end{abstract}

\vspace{0.10in}
\noindent
{\it Keywords}: \ Singularity and bifurcation theory; Persistent bifurcation diagram classification; Transition sets; Ideal membership problem; Standard and Gr\"obner bases.

\vspace{0.1in}
\noindent {\it 2010 Mathematics Subject Classification}: 37G10; 13P10; 58K50; 58K60.



\section{Introduction}

Many real life problems may result in the analysis of local zeros (around a zero solution, named a {\it base point}) of a smooth map
\be\label{m0}
f:\mathbb{R}^{n}\times \mathbb{R}^{m}\longrightarrow \mathbb{R}^{n}, \quad f(x,\alpha)=0, \quad \hbox{for }\; x=(x_1,\ldots, x_n)\in \mathbb{R}^{n}, \alpha=(\alpha_1, \ldots, \alpha_m)\in \mathbb{R}^{m}.
\ee We refer to each \(x_i\) by a state variable, \(n\) by state dimension and each \(\alpha_i\) by a parameter. Note that locating singular base points of a smooth map is
related to finding roots of nonlinear systems and is not the purpose of our work here; see \cite{Govaerts,MatcontACM,GovaertsBook}. Hence, we may assume that the base point is the origin and \(f(0, 0)=0\). Equation \eqref{m0} may demonstrate a surprising change on the solution set when the parameters vary. This occurs when the Jacobian matrix of \(f\) does not have a full rank. In this case we say that \(f\) is {\it singular}.

Equation \eqref{m0} may appear by direct mathematical modeling of a singular engineering problem or indirect through reduction methods such as
{\it Liapunov-Schmidt reduction}; \eg see \cite[Chapter VII]{GolubitskySchaefferBook}, \cite[Pages 156--162]{GovaertsBook} and \cite{Govaerts}.
For example, Equation \eqref{m0} appears in the study of equilibria and limit cycles of dynamical systems or steady-state solutions of partial differential equations. In fact, the theory described here is known as a ``natural framework'' for {\it equilibrium bifurcation theory}; see \cite{Govaerts}. Using Liapunov-Schmidt reduction, we can reduce the state dimension so that the Jacobian matrix at the origin is the zero matrix. Thus in this paper, we assume that
\be\label{m00} n=1 \quad \hbox{ and } \quad \frac{\partial f}{\partial x}(0,0)=0.\ee
We will deal with the case of multi-state dimensional problems in a future project. When \(f\) is a singular map and the parameters \(\alpha\) vary,
the number of solutions for Equation \eqref{m0} may change and any of such changes is called a {\it bifurcation}. The equation \(f(x,\alpha)=0\) is
called a {\it bifurcation problem} and the set
\be\lbrace (x,\alpha):f(x,\alpha)=0\rbrace\ee is called a {\it bifurcation diagram}; see \cite{Melbourne87,Keyfitz,GovaertsBook,GolubitskySchaefferBook,GatermannHosten,GatermannLauterbach,Gaffney,ArmbrusterKredel,Armbruster,MurdBook} for our main references of this subject.

\pagestyle{myheadings} \markright{{\footnotesize {\it M. Gazor and M. Kazemi \hspace{1.50in} {\it Bifurcation analysis for local zeros of smooth maps }}}}


Many powerful algebraic tools have been developed for local bifurcation analysis of zeros in Equation \eqref{m0}. Armbruster \cite{Armbruster} proposed a cognitive use of \Gr basis and encouraged a systematic implementations of the existing results in a computer. Yet to the best of our knowledge, there does not exist any (end-user) symbolic library for the {\it local bifurcation analysis of zeros of smooth maps}. This is a long overdue contribution despite its importance and wide applications. In the last two decades, there has been a considerable progress in development of computer algebra systems so that an efficient {\it symbolic} implementation of the results in singularity theory is now feasible; see \cite{GovaertsBook,Govaerts, JepsonSpence85, JepsonSpenceCliffe85} for numeric approaches.


A contribution here is an exposition of tools from algebraic geometry to the context of (locally) smooth maps (germs) for correct symbolic implementation of the results in bifurcation theory, where the involved subtleties are explicitly illustrated by counterexamples from bifurcation problems. Due to the infinite nature of Taylor series of smooth maps, the computations are performed modulo a given degree. We provide a sufficient condition for a given degree whose truncation does not lead to error. The default work of \texttt{Singularity} tests the condition and does the computations modulo an optimal degree. However, this approach adds a computational cost. Further, smooth maps involve flat functions (functions with zero Taylor series) and this may cause unnecessarily complicated formulas. Thereby, it is fundamentally helpful to use a ring smaller than the ring of smooth maps when it is feasible. Unlike the ring of formal power series, the associated computations in the ring of polynomials or fractional maps are exact and no truncation is required. We provide conditions along with rigorous proofs for the possible efficient implementations using ideals generated in either the rings of polynomials, fractional maps, formal power series or smooth germs; see
\cite[Chapters 1--4]{GolubitskySchaefferBook}, \cite[Chapters 6--7]{GovaertsBook}, \cite[Sections 6.2 and 6.3]{MurdBook} and \cite{Melbourne87,Armbruster,GatermannLauterbach,ArmbrusterKredel}. \texttt{Singularity} is adapted accordingly. We use identical notations and terminologies from \cite{GolubitskySchaefferBook,GovaertsBook,Melbourne87,Keyfitz,BeckerBook,CoxUsing,CoxLittleIdeals} as far as it is feasible. We have further written a user guide \cite{GazorKazemiUser} and constructed a comprehensive built-in {\sc Maple} help for \texttt{Singularity}.



\texttt{Singularity} computes a variety of algebraic structures associated with singular scalar maps including {\it tangent} and {\it restricted tangent spaces,} {\it high order term ideals,} and the {\it intrinsic ideals} associated with ideals of both {\it finite} and {\it infinite codimension.} Our {\sc Maple} library derives {\it low and intermediate order terms,} {\it normal forms, universal unfolding,} and {\it transition sets.} \texttt{Singularity} efficiently simplifies {\it intermediate} and {\it high order terms}. It, further, generates {\it persistent bifurcation diagrams (plot or animation),} and estimates the {\it transformations} transforming {\it contact-equivalent scalar maps to each other}. Finally, \texttt{Singularity} solves {\it the recognition problem for normal forms and universal unfoldings}. An interesting capability of \texttt{Singularity} is the classification of persistent bifurcation diagrams by generating an automatic list. This latter capability is in fact enabled by using a powerful built-in {\sc Maple} package called {\tt RegularChains} \cite{Marc}.

The rest of this paper is organized as follows. Singularity theory and bifurcation analysis of Equation \eqref{m0}-\eqref{m00} is discussed in
Section \ref{Sec1Int}. We further explain how singularity theory is related to {\it ideal membership problem} in algebraic geometry.
Section \ref{Sec2} describes how to treat the ideal membership problem. In this direction, computational algebraic tools such as standard and \Gr bases for ideals in three different rings, and the concept of {\it finite codimension} ideals are introduced. Truncation degree and alternative ring computations are discussed in Section \ref{TrunDegRing}. {\it Intrinsic ideals} and their associated ideal representations are discussed in Section \ref{ItrRepresent}. We further explain a procedure for computing the intrinsic part of an ideal or a vector space. Section \ref{Sec4} gives our suggestions on how to implement some objects and results from singularity theory. These implementations include high order term ideals, tangent spaces, transition sets, persistent bifurcation diagram classifications, normal forms and the universal unfolding. The capabilities of the main features of \texttt{Singularity} along with a few examples are sketched in Section \ref{SecFeatures}. Finally, Section \ref{Sec8} presents how {\tt Singularity} can be used in analysis and controller designs for a bucking problem.

\section{Bifurcation theory}\label{Sec1Int}

In this section we briefly describe the bifurcation theory and how our library can help the analysis. Due to the local nature of the problem \eqref{m0}-\eqref{m00}, we recall the notion of {\it smooth germs} around a base point. Two maps are considered as
{\it germ-equivalent} when both maps are locally identical; more precisely, when there exists a neighborhood of the base point so that both maps are equal
on the neighborhood. A {\it germ} is a {\it germ-equivalence} class of a smooth map. We denote \(\E\) for the set of all scalar smooth germs whose base point is
the origin. From now on we merely work with elements of \(\E\) rather than a scalar smooth map; see \cite[156]{GovaertsBook}.

Following \cite{GolubitskySchaefferBook} and \cite[Chapter 7]{GovaertsBook}, we study the local zeros of maps when there is only one distinguished parameter denoted by \(\lambda\), \ie
\be\label{eq1}
g:\mathbb{R}\times \mathbb{R}\longrightarrow \mathbb{R}, \quad g(x,\lambda)=0,\quad g(0,0)=g_{x}(0,0)=0, \; \hbox{ and } \; m:=1.
\ee The effect of additional parameters may be treated by study of their small perturbations. The main goal of this theory is to classify {\it qualitative} types of Equation \eqref{eq1} and its arbitrary small perturbations. In order to achieve this goal, we first define a {\it qualitative property} as a property that is
invariant under an appropriate equivalence relation. Here, we use {\it contact-equivalence} relation:
\begin{itemize}
  \item We say that the germs \(g\) and \(h\) are contact-equivalent when there exist a smooth germ \(S(x,\lambda)\) and diffeomorphic germs \(X\)
  and \(\Lambda\) such that
\be\label{ContEqvl} g(x,\lambda)= S(x,\lambda)h(X(x,\lambda),\Lambda(\lambda)),  \ee where  \(S(x,\lambda)>0,\) \(X_{x}(x,\lambda)>0\) and
\(\Lambda^{'}(\lambda)>0.\)
\end{itemize} 

\begin{figure}[h]
\begin{center}
\subfigure[\label{a}]{\includegraphics[width=.26\columnwidth,height=.2\columnwidth]{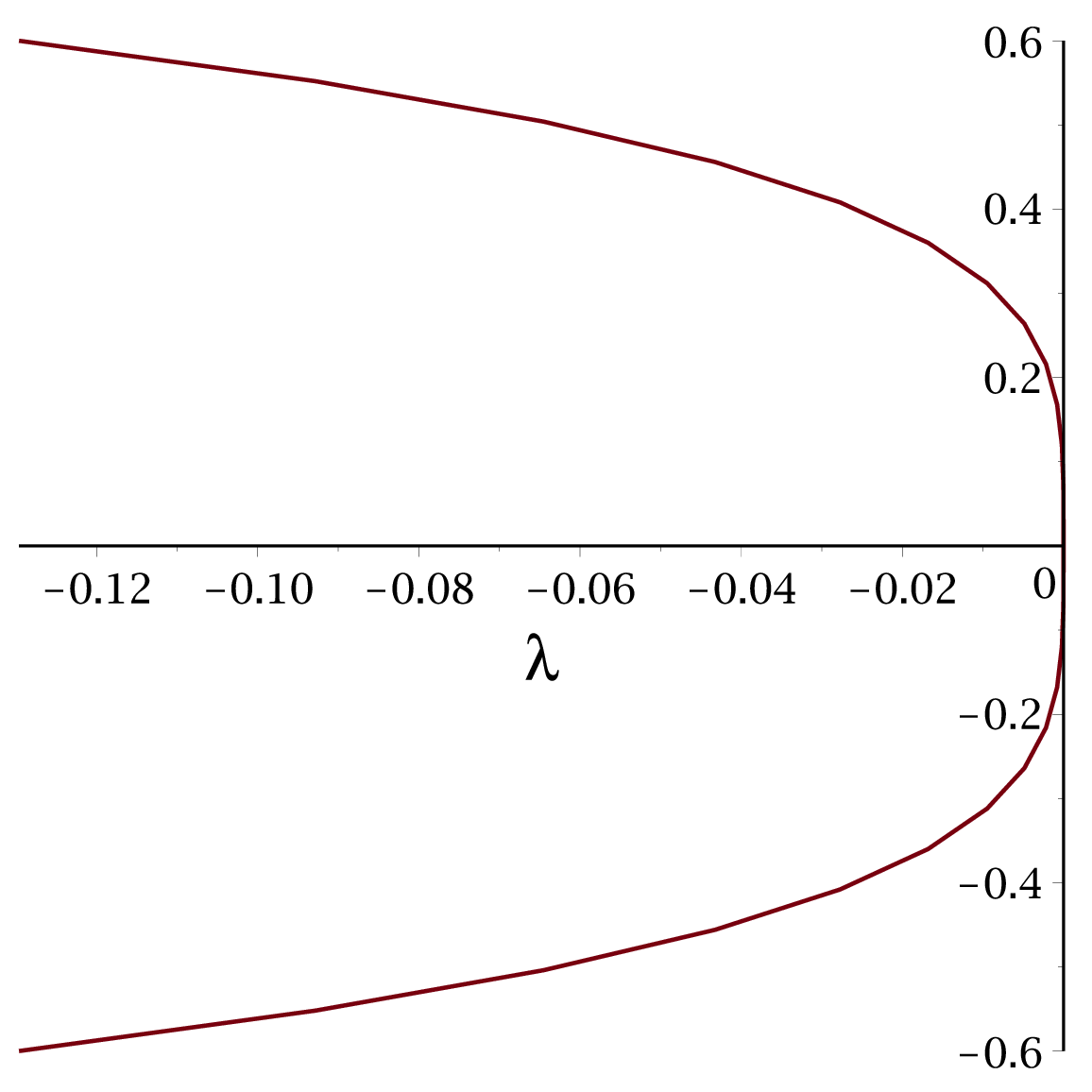}}
\subfigure[\label{b}]{\includegraphics[width=.26\columnwidth,height=.2\columnwidth]{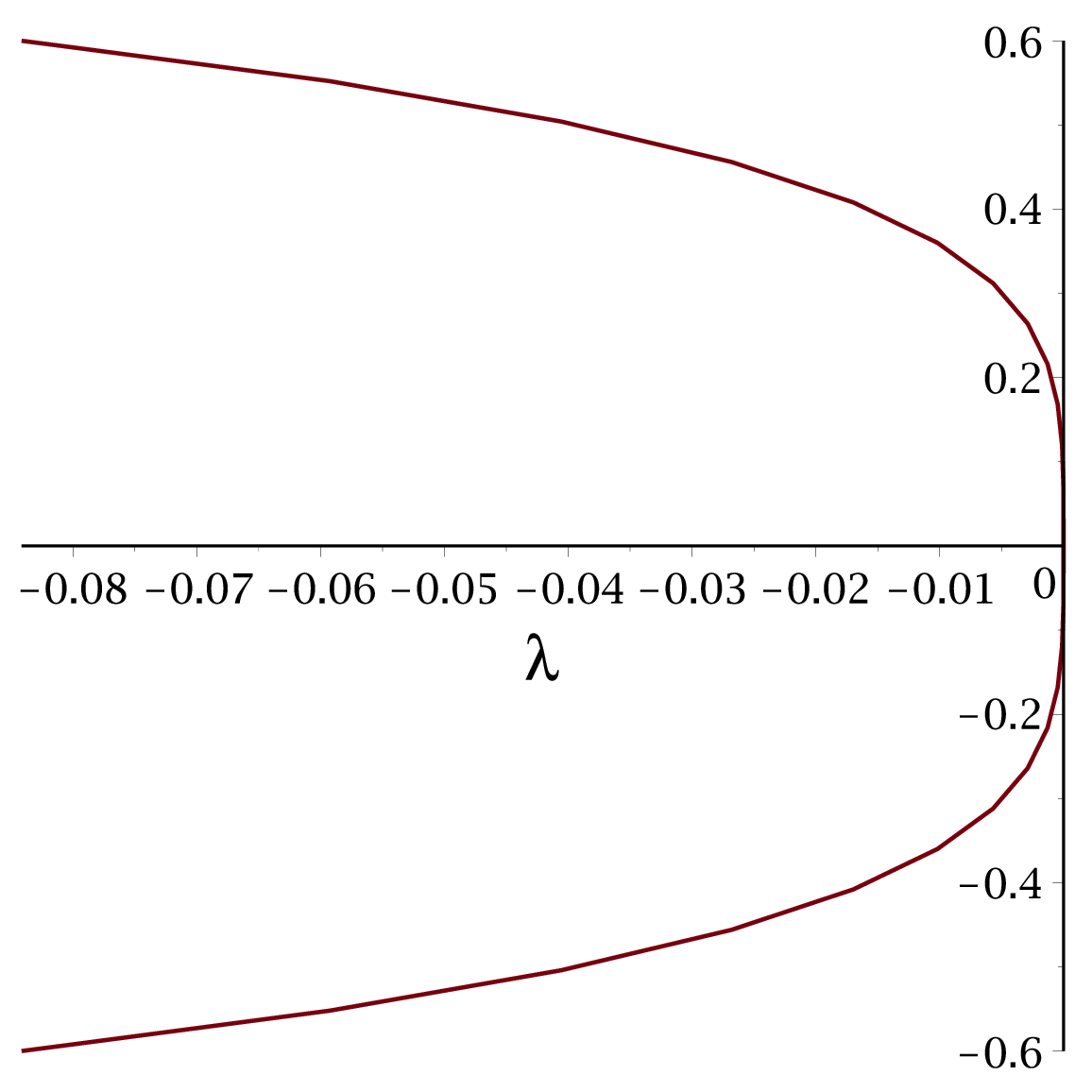}}
\subfigure[\label{c}]{\includegraphics[width=.26\columnwidth,height=.2\columnwidth]{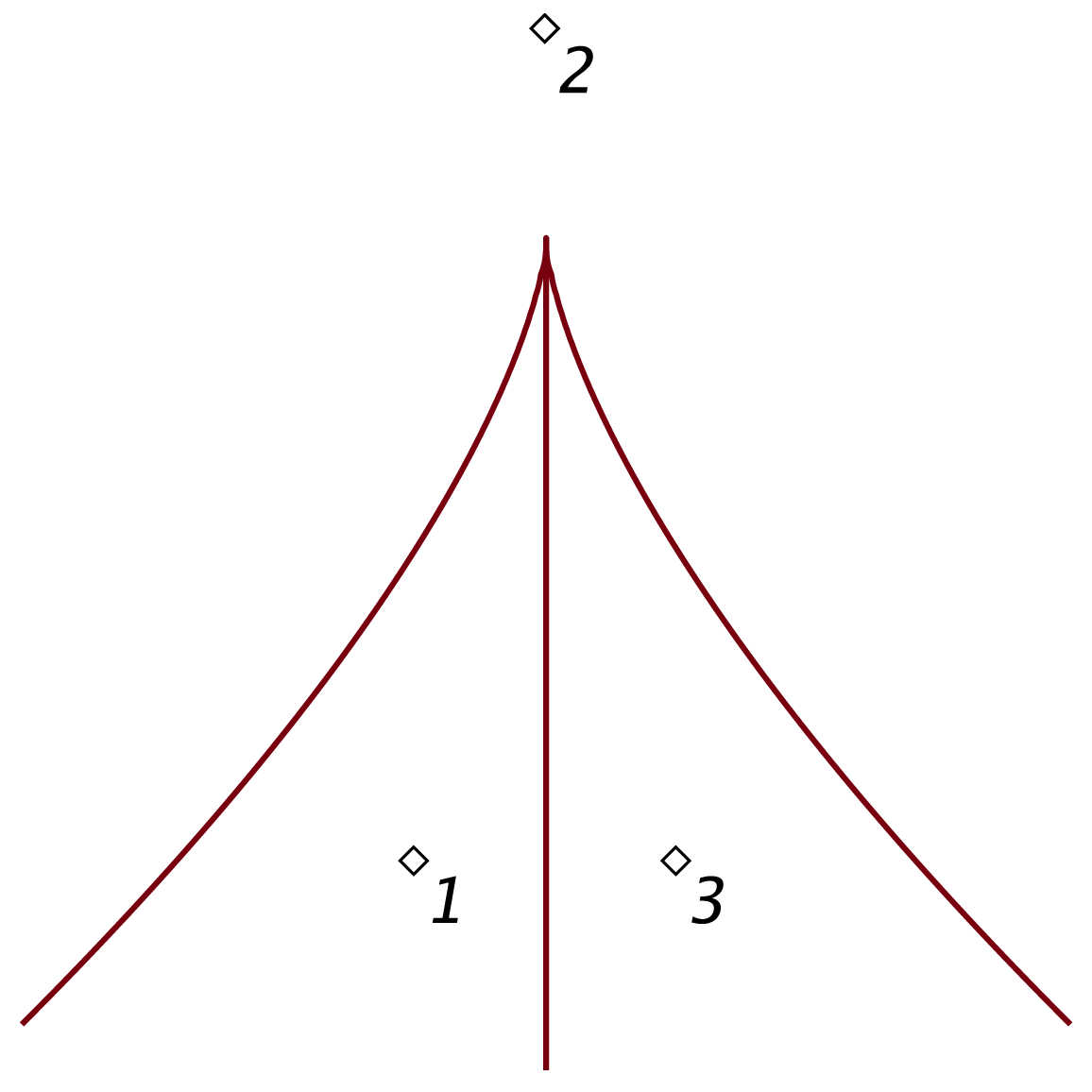}}
\subfigure[]{\includegraphics[width=.26\columnwidth,height=.2\columnwidth]{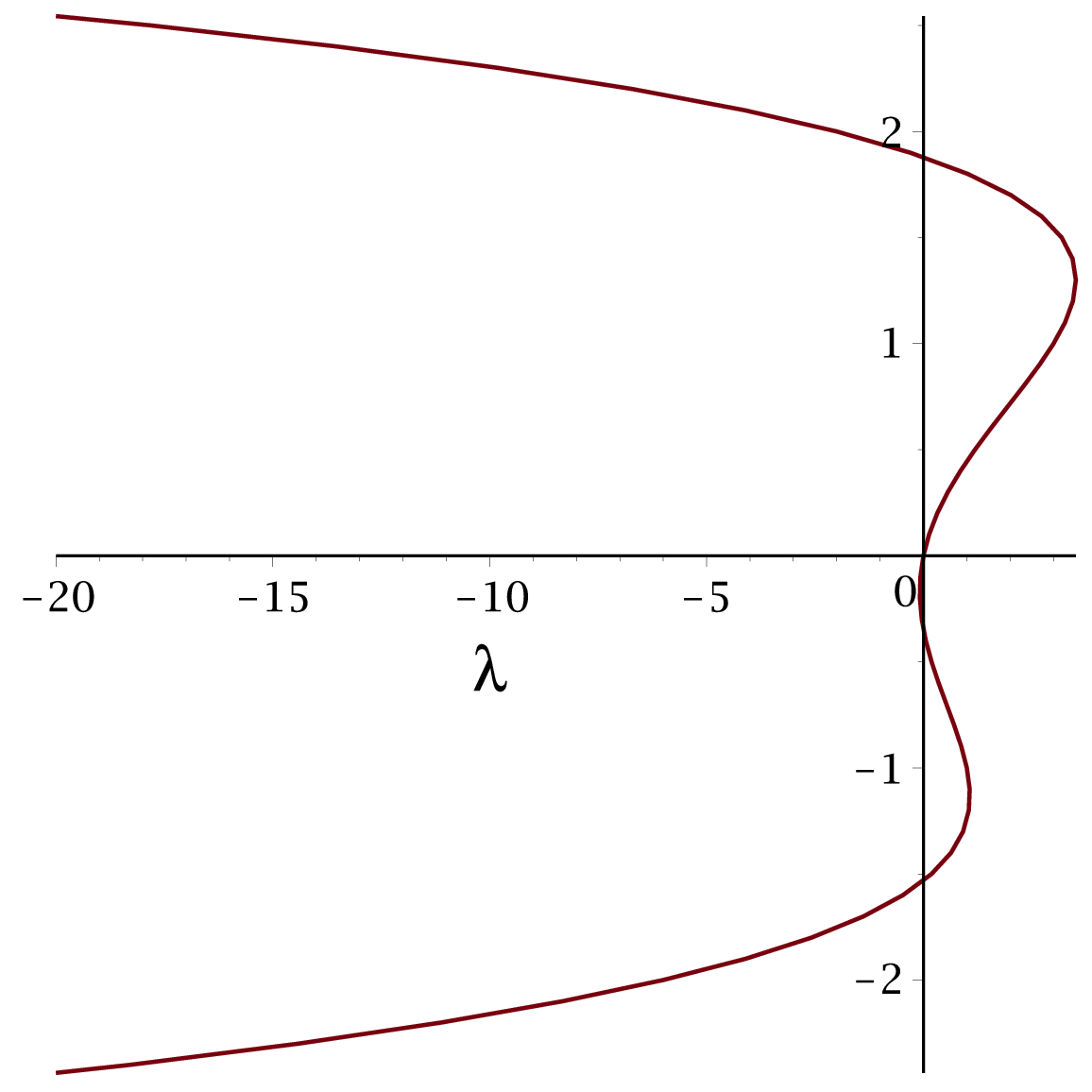}}
\subfigure[]{\includegraphics[width=.26\columnwidth,height=.2\columnwidth]{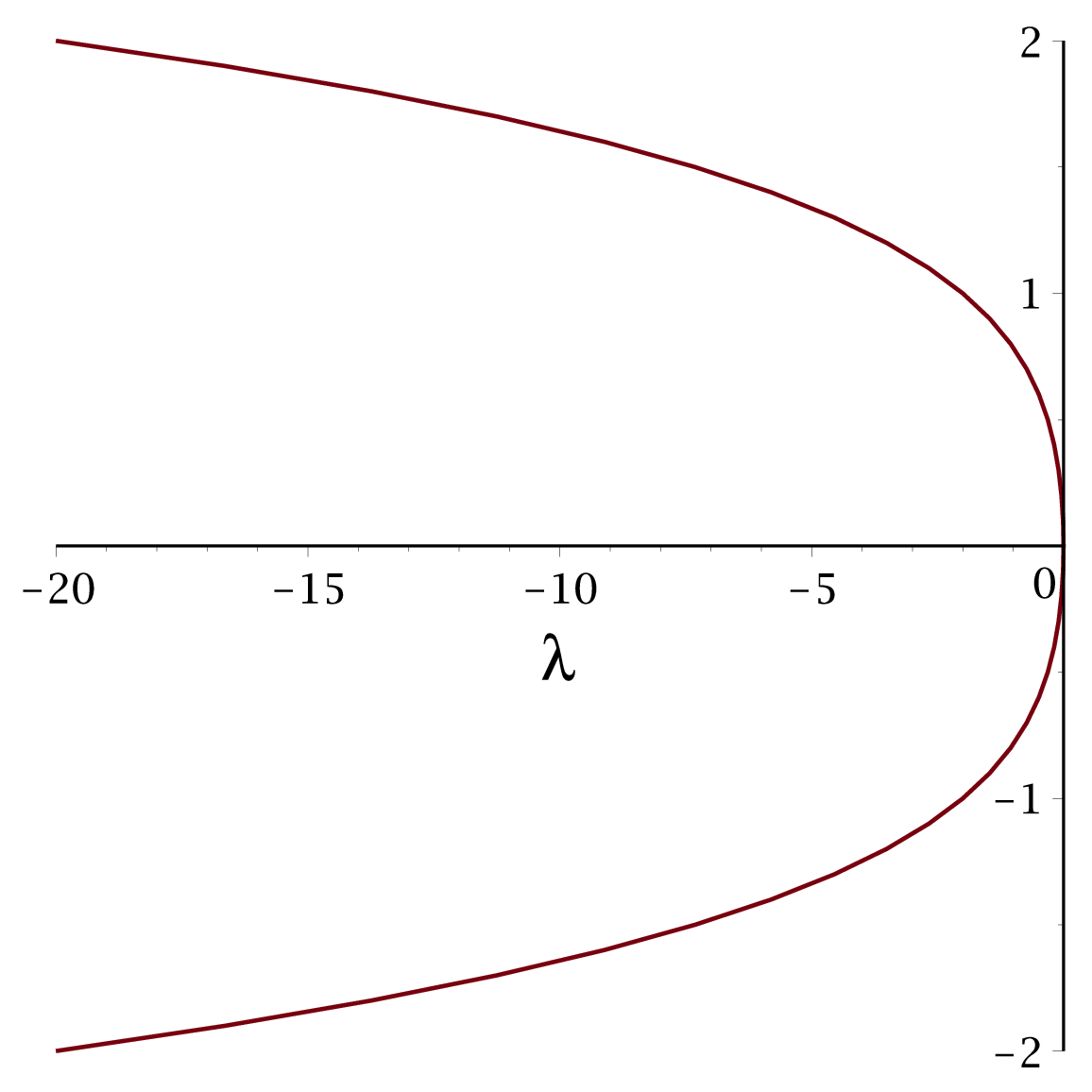}}
\subfigure[]{\includegraphics[width=.26\columnwidth,height=.2\columnwidth]{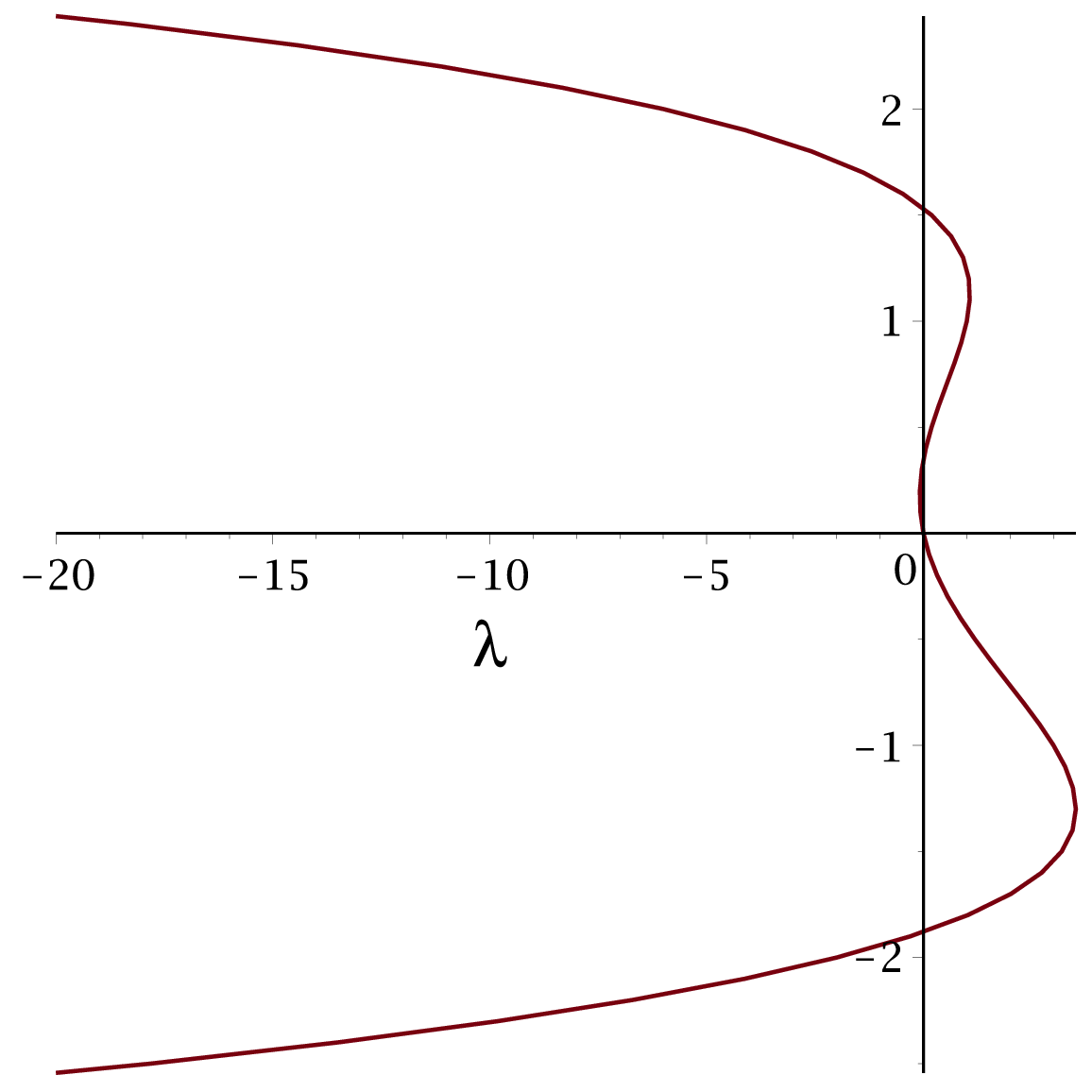}}
\caption{\small The vertical axes stand for the state variable \(x.\) Diagrams \ref{a} and \ref{b} depict bifurcation diagrams of \(g\) in Equation \eqref{Exmp1} and its normal form \(f,\) respectively. The diagram \ref{c} shows the transition set and the regions 1-3 associated with Equation \eqref{UnivExmp}. The second row illustrates the persistent bifurcation diagrams associated with these regions. }\label{Fig1}
\end{center}
\end{figure}

The principal idea in bifurcation theory lies in the notion of {\it normal forms}. To study the local zeros of \eqref{eq1}, we choose a representative (say \(f\)) from contact-equivalent class of \(g\) that is considered to be the {\it simplest} for the analysis and call it a {\it normal from}. In order to compute the normal form of a singular germ, we need to compute certain ideals in the {\it local ring of smooth germs}.
This signifies the importance of the well-known ideal membership problem in algebraic geometry, that is, deciding what kinds of germs belong to an ideal
generated by a given set. One may study the zero structures of the normal forms and then, conclude about the solution behavior of Equation \eqref{eq1}.
For instance, let
\be\label{Exmp1} g(x, \lambda):= \exp(x^{2})+2\cos(x)-3+\sin(\lambda).\ee Using the command \({\tt Normalform(g, [x, \lambda], 5)}\) in {\tt Singularity}, we obtain its normal form by \(f(x, \lambda):= \frac{7}{12} x^4+\lambda.\) The bifurcation diagrams of \(f\) and \(g\) are depicted in Figure \ref{a} and \ref{b}. Here, \({\tt Transformation(g, f, [x, \lambda], 5)}\) provides an approximation
\bas
X(x,\lambda)&:=& \lambda+ x+\lambda^2 +\lambda x,\\
S(x,\lambda)&:=&1+\frac{1}{6}\lambda^{2}-\frac{7}{12}\lambda^{3}-\frac{7}{3}\lambda^{2}x-\frac{7}{2}\lambda x^{2}-\frac{7}{3}x^{3}-\frac{833}{360}\lambda^{4}-\frac{28}{3}\lambda^{3}x-14\lambda^{2}x^{2}-\frac{28}{3}\lambda x^{3}-\frac{7}{3}x^{4}
\eas
modulo degree 5 for \((X(x, \lambda), S(x, \lambda))\) which transforms \(g\) into \(f\) via Equation \eqref{ContEqvl}. The locally invertible transformation \(X(x, \lambda)\) sends the bifurcation diagram of \(f\) into that of \(g\).

Real life problems can not be perfectly modeled by a system of equations and {\it imperfections} are always inevitable. Furthermore, the singularity of a germ \(g\) implies that the zeros of a small perturbation of \(g,\) say
\be\label{Unfold}
h(x, \lambda, \beta)=0, \hbox{  where } h(x, \lambda, 0)=g (x, \lambda),
\ee
may behave substantially different than what the zero structure associated with \(g=0\) does. The parameterized germ \(h\) in \eqref{Unfold} is called an {\it unfolding} of \(g.\) Hence, {\it modeling imperfections} and the possible existence of additional parameters in a model are the main obstacles of simply using normal forms of a singular germ for the qualitative understanding of a real life problem. Thus, the approach needs to be refined through the notion of {\it universal unfolding}. In fact, we are interested to find a parameterized family like
\be\label{eq4}
G(x, \lambda, \alpha)=0, \quad \alpha\in \mathbb{R}^p,
\ee such that for any small perturbation \(\epsilon p(x, \lambda, \epsilon)\), the germ \(g+\epsilon p(x, \lambda, \epsilon)\) would be contact-equivalent to
\(G(x, \lambda, \alpha(\epsilon))\) for some germ \(\alpha(\epsilon).\) We call such germ \(G\) a {\it versal unfolding} of \(g.\)
A versal unfolding with the minimum possible number of parameters is called a {\it universal unfolding} for \(g\); see \cite[Definitions 1.1 and 1.3, Pages 120--121]{GolubitskySchaefferBook} and
\cite[Definitions 6.4.2 and 6.4.3, Page 176]{GovaertsBook}.
The number of parameters in a universal unfolding is named the {\it codimension}
of \(g.\) The universal unfolding accommodates any possible modeling imperfections, any arbitrary small perturbation and also the existence of any possible
number of parameters (in addition to the distinguished parameter \(\lambda).\) In order to derive the universal unfolding of a given singularity, the
computation of a vector space called {\it tangent space} or instead a basis for its complement is required. Using {\tt Singularity}'s command \({\tt UniversalUnfolding(g, [x, \lambda], 5, normalform)}\) for \(g\) in \eqref{Exmp1} gives rise to
\be\label{UnivExmp} G(x, \lambda, \alpha_1, \alpha_2):= x^4+\lambda+\alpha_{1}x+\alpha_{2}x^2.\ee

{\it Bifurcation diagram classification} of the universal unfolding gives an insight to the zero structure of a germ and any of its perturbations. This is
studied by the notion of {\it persistence} in the bifurcation diagrams. In fact a bifurcation diagram is called {\it persistent} when all its small perturbations remain self contact-equivalent, and otherwise it is called {\it nonpersistent}. Finding nonpersistent systems and their associated subset of the parameter space
(the so-called {\it transition set} \(\Sigma\)) play a central role in this classification. More precisely, all parametric germs associated with
parameters in a connected component of \(\Sigma^c,\) complement of \(\Sigma,\) are contact-equivalent. Therefore by choosing a parameter from each connected
component of \(\Sigma^c,\) a complete list of persistent bifurcation diagrams modulo contact-equivalent is obtained. The non-persistence may either originate
from a local nature or be caused by the singular boundary conditions. Local nonpersistent bifurcation diagrams are determined with families of germs associated
with three semi-algebraic parameter spaces of codimension one; i.e., {\it bifurcation} \(\mathscr{B}\), {\it hysteresis} \(\mathscr{H}\) and
{\it double limit point} \(\mathscr{D}\); see \cite[Page 140]{GolubitskySchaefferBook} for details. Nonpersistent germs associated with boundary conditions add
extra complications into the solution dynamics, when bifurcation diagrams cross the boundary;
see \cite[Pages 154--165]{GolubitskySchaefferBook}. Given our description for any finite codimension singularity, the connected components of the complement set
of the transition  set, \ie \(\Sigma^{c},\) provide the qualitative classification of persistent bifurcation diagrams. The command \(\verb|TransitionSet|(G, [\alpha_1, \alpha_2], [x, \lambda])\) for \(G\) in \eqref{UnivExmp} gives
\bes\mathscr{B}:=\emptyset, \quad \mathscr{H}:=\{(\alpha_1, \alpha_2)\,|\, 8\alpha_{2}^{3}=-27\alpha_{1}^{2}\}, \quad \mathscr{D}:=\{(\alpha_1, \alpha_2)\,|\,\alpha_{1}=0, \alpha_{2}\leq 0 \rbrace,\ees and \(\Sigma:= \mathscr{H}\cup \mathscr{D}.\) The transition set \(\Sigma\) is plotted in Figure \ref{c}. \(\verb|PersistentDiagram|(G, [x, \lambda],\break \verb"plot", \verb"ShortList")\) plots a complete list including all contact-inequivalent types of persistent bifurcation diagrams for \(G\); see the second row in Figure \ref{Fig1}.

\section{Ideal membership problem and tools from algebraic geometry}\label{Sec2}

Given our description in the previous section, we mainly need to address the {\it ideal membership problem} in the ring of smooth germs. We call an {\it ideal basis} for the (finite) {\it generators} of a finitely generated ideal. The {\it ideal membership problem} refers to the question on whether or not an element belongs to a predefined ideal. Although a {\it finitely generated ideal} is defined by an ideal basis, yet it is not an easy task to understand whether a given element belongs to the ideal or not. There are two different ways to facilitate the ideal membership problem. One is to find a convenient ideal presentation from which the ideal membership can be easily understood. This is feasible for certain ideals and will be addressed in Section \ref{ItrRepresent}. The second is a computational approach based on a division algorithm. Dividing a given element by the ideal basis (generators), the element certainly belongs to the ideal when the remainder is zero. However for non-zero remainders with respect to an arbitrary ideal basis, the division may not help to conclude about the membership of the element in the ideal. B. Buchberger in 1965 introduced an {ideal basis} (called {\it \Gr basis}) for ideals in the polynomial ring on which zero/nonzero remainders would indeed conclude the ideal membership; see \cite{BuchAlg,BuchCrit,BuchSystem,BuchThesis}.


Armbruster and Kredel \cite{ArmbrusterKredel} suggested a cognitive use of \Gr basis for computing universal unfolding; also see \cite{Armbruster}.
Then, Gatermann and Lauterbach \cite{GatermannLauterbach} extended the tools for equivariant bifurcation problems. Wright and Cowel in
\cite{CowelWright} noticed that a local version of \Gr basis is indeed the appropriate tool to work with it in this theory. In fact, naive uses of \Gr basis yield wrong results for many singularities (\eg see Example \ref{m15}). The main reason is that \Gr basis may work with polynomial germs, while the ring of smooth germs is a considerably larger ring than the polynomial germs. Recall that a set may generate a larger ideal in a ring than what it generates in a subring. Since Wright and Cowel's remark, the only result in this direction is due to Gatermann and Hosten \cite{GatermannHosten}. This is a fundamental contribution, but it is yet incomplete. They used (mixed) {\it standard basis} for ``\emph{mixed modules}'' over fractional maps and multi-dimensional state variables. They discussed neither of Taylor series truncations, smooth maps and germs, nor computations with smooth maps of finite or infinite codimensions. Their algorithms are indeed useful for possible implementations in a computer algebra system to simplify certain fractional maps with finite codimension. Although they did not discuss their approach limitations, their suggested algorithms are limited to the cases 
when their associated (restricted) tangent space is a zero dimensional ideal in the ring of fractional (germs) maps. However, they are not yet sufficient to compute {\it normal forms} and the associated normalizing transformations. 
We further remark that any useful implementation of the results in a computer algebra system for their real life applications in bifurcation theory needs to treat arbitrary {\it singular smooth maps}. As far as our information is concerned, there does not yet exist any contribution to discuss on how to implement the results in bifurcation theory for actual bifurcation analysis of singularities; see Section \ref{Sec4} in this direction. Our suggested algorithms also provide a systematic approach for detection and treatment of infinite and high codimensional singularities and their possible persistent bifurcation diagram classifications.

The localization of the polynomial-germ ring gives rise to the fractional germs, \ie germs with a fractional representative in their germ-equivalent class. Although the set of all fractional germs is a local ring, yet it is still a much smaller ring than the ring of smooth germs. The other alternative is the local germ ring of all formal power series. This is a larger ring than the fractional germs and it is, perhaps, suitable due to Borel lemma. Borel lemma indicates the existence of a one-to-one correspondence modulo flat functions between smooth germs and formal power series through their Taylor series expansion. 

\subsection{Standard and \Gr bases for ideals}

Let \(K\) be a field of characteristic zero; in particular we are interested in the field of real numbers. For our convenience, we simply identify a given germ with a convenient representative of that germ. For instance, we talk about the polynomial germ ring over the field \(K\) and denote it by \(K[x, \alpha]\) while we mean the ring of all smooth germs whose germ-equivalent class have a polynomial representative. The quotient ring of \(\E\) over the ideal of all flat germs is denoted by \(K[[x, \alpha]]\) due to the fact that it is ring-isomorphic to the ring of formal power series. Thus, we call \(K[[x, \alpha]]\) the germ ring of {\it formal power series}. Further, we identify members of \(K[[x, \alpha]]\) with the infinite Taylor expansion of their representative.
Since \(K[x, \lambda]\) and \(K[[x, \lambda]]\) are Noetherian rings, we can guarantee termination of most algorithms during the computations.

Denote \(\deg(X^\alpha):= \alpha_{1}+ \ldots+ \alpha_{n}\) when \(X^{\alpha}:= x_1^{\alpha_{1}}x_2^{\alpha_{2}}\cdots x_n^{\alpha_n},\)
\(X:=(x_1, x_2, \ldots, x_n)\) and \(\alpha:=(\alpha_1, \alpha_2, \ldots, \alpha_n)\in Z_{\geq 0}^n.\) Any expression like \(X^{\alpha}\) is called a
{\it monomial germ} while a {\it term} in \(\E\) means a monomial germ along with its coefficient. We define the lexicographic ordering \(\prec_{\lex}\)
on monomial germs \(X^{\alpha}\) as follows:
\bes X^{\alpha}\prec_{\lex}X^{\beta}\hbox{ when } \alpha_i-\beta_i\hbox{ is negative for } i:= \inf \{j\,|\, \alpha_j-\beta_j\neq0\}.\ees

\begin{defn}
A {\em local order} \(\prec\) is a total ordering (every two terms are comparable) and furthermore,
\begin{itemize}
  \item  For any \(\alpha, \beta, \gamma\in Z_{\geq 0}^n\), the condition \(X^\alpha \prec X^\beta\) implies \(X^{\alpha+\gamma}\prec X^{\beta+\gamma}\).
  \item \(x_i \prec 1\) for all \(i=1, \ldots, n\). We further assume that \(0\prec X^\alpha\) for any \(\alpha.\)
\end{itemize}
\end{defn} Due to Dickson's lemma (see \cite[Page 251]{Becker} and \cite[Page 71]{CoxLittleIdeals}), every (infinite) set of monomials have a maximum with
respect to any arbitrary local order; here, the condition \(X^\beta| X^\alpha\) implies \(X^\alpha\preceq X^\beta.\)
An important example of a local order is {\it anti-graded lexicographic} ordering  \(\prec_{\alex}\) defined by
\begin{equation*}
X^\alpha\prec_{\alex} X^\beta \Leftrightarrow \left\{
\begin{array}{lcl}
\deg(X^\alpha) > \deg(X^\beta),&&\\
&&\\
\deg(X^\alpha) = \deg(X^\beta)&\text{and}& X^\alpha\prec_{\lex} X^\beta.
\end{array}
\right.
\end{equation*}

The localization of the polynomial germ ring is defined as
\bes \mathscr{R}:= \left\{\frac{f}{g} \mid f, g\in K[x, \lambda], g(0,0) \neq 0 \right\}\ees
whose unique maximal ideal is generated by \(x\) and \(\lambda.\) It is common to denote \(\RM\) with \(K[x, \lambda]_{\langle x, \lambda\rangle}.\)
We call \(\RM\) the {\it ring of fractional germs}. Throughout this paper, we denote
\bes \mathcal{R} \hbox{ for either of the local rings } K[[x, \lambda]], \RM, \hbox{ or } \E,\ees
unless it is explicitly stated.

\begin{defn}
Let \(f\in\mathcal{R}\) and \(\prec\) be a local order. The infinite jet and \(k\)-jet of \(f\) are defined by its Taylor series expansion around the origin and denoted by
\bes J^{\infty}(f):= \sum_{(i,j)} a_{ij}x^{i}\lambda^{j}, \; \hbox{and }\; J^{k}(f):= \sum_{i+j\leq k} a_{ij}x^{i}\lambda^{j} \quad \hbox{ for }  a_{ij}\in K.\ees
The set of terms of \(f\) are defined by
\bes \Terms(f):= \left\{ a_{ij} x^{i} \lambda^{j}\,|\, a_{ij}\neq 0\right\},\ees
\ie all terms appearing in \(J^{\infty}(f).\) When \(\Terms(f)\neq\emptyset,\) the leading term of \(f\) is
\bes \LT(f)=\max\, \Terms(f),\ees
\ie \(\LT(f)\in \Terms(f)\) and \(p\prec \LT(f)\) for any \(p\in \Terms(f)\setminus \{\LT(f)\}.\) The germ \(g\) is flat iff \(\Terms(f)=\emptyset.\) We define
\(\LT(f):=0,\) when \(f\) is a flat germ, \ie \(J^{\infty}(f)=0.\) The coefficient and monomial of the leading term are respectively called the {\it leading coefficient} (\(\LC\)) and {\it leading monomial} (\(\LM\)). For the case of \(f\in K[x, \lambda]\) and \(\prec_{\lex},\) we may similarly define \(\LT(f), \LC(f)\) and \(\LM(f).\)
\end{defn}

Now we present some definitions, terminologies and theorems from \cite{CoxUsing,CoxLittleIdeals,GreuelPfister,Becker}. These are suitably
modified and generalized to fit in our purpose.

\begin{defn}[Remainder]\label{m11}
A {\it remainder} of a germ \(f\) with respect to the set of germs \(G=\{g_{1},\ldots, g_{m}\}\subset \E\) and the local order  \(\prec\) is defined as a
germ \(\Rem(f,G,\prec)\in\E\) so that
\begin{itemize}
\item[(1)] \(\Rem(f,G,\prec)=f-\sum_{i=1}^{m}q_{i}g_{i}\), for some \(q_{1},\ldots,q_{m}\in\E\) so that \(\LT(q_{i}g_{i})\preceq \LT(f)\).
\item[(2)] No term of \(\Rem(f,G,\prec)\) is divisible by any of \(\LT(g_{i})\) for \(i=1,\ldots, m\).
\end{itemize}
\end{defn}
By replacing \(\E\) with \(K[[x, \lambda]]\) in Definition \ref{m11}, the remainders in the ring of \(K[[x, \lambda]]\) is readily defined; also see
\cite[Pages 251-252]{Becker}. The same is true for the case of polynomial germ ring \(K[x, \lambda]\) provided that the local order \(\prec\) would be
instead a monomial ordering like \({\prec_{\lex}}.\) The remainders in the ring of \(\RM\) is defined in \cite[Page 170]{CoxUsing} using {\it Mora normal form algorithm}. Usually, the terminology of {\it normal form} is used rather than {\it remainder}. However, we choose {\it remainder} as the other may cause confusion with the {\it normal forms} of germs in singularity theory. The division here is related to Malgrange preparation theorem and Mather division theorem; see \cite[Corollaries A.6.2 and A.7.2, Theorem A.7.1]{MurdBook}.

When \(f\) is flat, its remainder with respect to any set of germs and local order is flat. The remainder is not necessarily unique even modulo flat germs; also see \cite[Pages 251-252]{Becker}. In fact, the remainder is unique modulo flat germs when \(G\) is a {\it standard basis} (or \Gr basis for the case of \(K[x, \lambda]\)), defined as follows.

For an ideal \(I\) in \(\mathcal{R},\) we define the leading term ideal \(\LT(I)_{\mathcal{R}}\) by
\begin{equation*}
\LT(I)_{\mathcal{R}}:=\langle \LT(f):f\in I\rangle_{\mathcal{R}}.
\end{equation*}

\begin{defn}[Standard basis]\label{SBasis}
Let \(I\) be an ideal in \(\mathcal{R}\) with a finite generating set \(\lbrace g_1, \ldots, g_m\rbrace\subset I\). The set \(\lbrace g_1, \ldots, g_m\rbrace\)
is called a {\em standard basis} of \(I\) when
\begin{equation}\label{SBCondition}
 \LT(I)_{\mathcal{R}}=\langle \LT(g_1), \ldots, \LT(g_m)\rangle_\mathcal{R}.
\end{equation}
The set \(\lbrace g_1, \ldots, g_m\rbrace\subseteq \mathcal{R}\) is called a standard basis in \(\mathcal{R}\) when it is a standard basis for the ideal
\(\langle g_1, \ldots, g_m\rangle_\mathcal{R}\).
\end{defn}
\begin{rem}\label{New}
\begin{itemize}
\item[(a)] The set \(\lbrace g_1, \ldots, g_m\rbrace\subset K[x, \lambda]\) is called \Gr basis with respect to \(\prec_{\lex}\) when the condition
 \eqref{SBCondition} holds.
\item[(b)] Any finite set of germs \(\{g_i\;|\; 1\leq i\leq m\}\) in \(\E\) is a standard basis in \(\E\) iff the set of formal power series
\(\{J^\infty(g_i)\}\) is a standard basis in \(K[[x, \lambda]].\)
\item[(c)] Let \(\mathcal{R}:= K[x, \lambda], K[[x, \lambda]]\) or \(\RM\) and
\(\lbrace g_1, \ldots, g_m\rbrace\subset I\) be a finite set. Then,
\bes \LT(I)_{\mathcal{R}}=\langle \LT(g_1), \ldots, \LT(g_m)\rangle_\mathcal{R}\ees implies that
\(\lbrace g_{1},\ldots, g_{m}\rbrace\subset I\) is a generating set for \(I\). This has a simple argument as follows; also see \cite[Page 206]{BeckerBook}. For any
  \(f\in I,\)
  \bes r:=\Rem \big(f, \{g_i|\, i=1, \ldots, m\}, \prec\big)\in I.\ees
  Thus, \(\LT(r)\in \LT(I)_\mathcal{R}\) and \(\LT(r)\) must be factored by \(\LT(g_i)\) for some \(i.\) The latter implies that \(r=0,\) otherwise this contradicts with \(r\) being a remainder.
\end{itemize}
\end{rem}

Computation of standard basis uses the notion of \(S\)-germs.
Let \(f, g\in \R\) and \(\preceq\) be a local order. Then, {\it S-germ} of \(f\) and \(g\) is defined by
\begin{equation*}
S(f,g)=\left\{
\begin{array}{lcl}
\frac{\LCM\big(\LM(f),\LM(g)\big)}{\LT(f)}f-\frac{\LCM\big(\LM(f),\LM(g)\big)}{\LT(g)}g& \hbox{ if }& J^{\infty}(f)J^{\infty}(g)\neq0,\\
0 & \hbox{ for }& J^{\infty}(f)J^{\infty}(g)=0.
\end{array}
\right.
\end{equation*}
Here \(\LCM\) stands for the least common multiple for a pair of monomials.

\begin{thm}(See \cite{Hironaka}, \cite[Theorem 2.6]{GazorKazemiArxiv}, and Hironaka theorem on \cite[Page 252]{Becker}) \label{m10}
Let \(\R=\E\) or \(K[[x, \lambda]]\), \(G=\lbrace g_{1},\ldots, g_{m}\rbrace\subsetneq \R,\) \(0\neq f\in\R\) and \(\prec\) be a local order. Then,
\begin{itemize}
\item[(a)] Always \(f\) has a remainder with respect to \(\prec\) and \(G.\)
\item[(b)] The set \(G\) is a standard basis iff the remainder of \(f\) with respect to \(G\) and
\(\prec\) is unique modulo flat germs.
\item[(c)]{(Buchberger's Criterion)} The set \(G=\lbrace g_{1},\ldots, g_{m}\rbrace\) is a standard basis iff for all \(i, j,\) the expression \(\Rem(S(g_{i},g_{j}),G,\prec)\) is a flat germ.
\item[(d)] The set \(G\) is a standard basis iff \(\Rem(f,G,\prec)\)
is flat for all \(f\in\langle G\rangle_{\R}\).
\end{itemize}
\end{thm}

In order to compute the standard basis for an ideal in \(K[[x, \lambda]]\) (and also in \(\E\)), one needs to sequentially enlarge and update
its set of generators by only adding the {\it non-flat} remainders (with respect to the updated generators and a given local order) of the \(S\)-germs of the
generator pairs. This process is usually known as {\it Buchberger algorithm} and it terminates when the ascending ideals generated by the leading terms of the
updated generators stops any further enlargement. The Buchberger algorithm is finitely terminated due to the fact that \(K[[x, \lambda]]\) is a
Noetherian ring.

\begin{exm}
\begin{itemize}
\item[(a)] The division of a polynomial by a set of polynomials in our division algorithm may involve formal power series; for example we have
\(\Rem(1, 1-x, \prec_{\alex})=0=1-(\sum_{n=0}^{\infty}x^{n})(1-x).\)
\item[(b)] The remainder of a polynomial divided by a polynomial may give rise to an {\em infinite} formal power series. For instance let
\(G:=\{x\lambda-x^{2}\lambda^{2}-x^{4}\}.\) Then,
\bas h&:=&\Rem(x^{2}\lambda, G, \prec_{\alex})=x^{5}+ x^{9}+ 2x^{13}+ 5x^{17}+ 14x^{21}+ \ldots
\\&=& x^{2}\lambda-\big(x+x^{2}\lambda+x^{5}+ x^{3}\lambda^{2}+ 2x^{6}\lambda+ x^{4}\lambda^{3}+ \ldots\big)(x\lambda- x^{2}\lambda^{2}- x^{4}).
\eas
Since the generator of an ideal with only one generator is always a standard basis, the remainder here is unique.
\item[(c)] This example is to show that a finite set of polynomial ideal basis in \(\E\) may lead to a standard basis that includes non-polynomial germs.
Let \(G=\{f:=x\lambda-x^{2}\lambda^{2}-x^{4}, g:=\lambda-x\lambda-x\lambda^2-x^3\}\).  Then, \(S(f, g)=x^2\lambda.\) It is easy to verify that \(S=\{ f, g, h\}\)
is a standard basis, where \(h\) is defined in part (b).
\end{itemize}
\end{exm}

The first and second Buchberger criteria are applied for efficient computation of standard basis in
\texttt{Singularity}. However, we will not discuss them in this paper.
\begin{defn}[Reduced standard germ basis]\label{rsgb}
Let \(G=\lbrace g_{1},\ldots, g_{n}\rbrace\) be a standard basis and \(\LC(g_{i})=1\) for \(i=1, \ldots, n.\) When
\begin{itemize}
   \item \(\LT(g_{i})\nmid p\) for all \(p\in \Terms(g_{j})\), except for when \(p=\LT(g_{i})\) and \(i=j\),
\end{itemize} the set \(G\) is called a {\it reduced standard basis}.
\end{defn}

Given a local order \(\prec\), \cite[Theorem 2.1, Page 255]{Becker} states that any ideal in \(K[[x, \lambda]]\) has a unique reduced standard basis.
\begin{thm}(See \cite[Theorem 2.9]{GazorKazemiArxiv})\label{m1}
With respect to any local order \(\prec\), every finitely generated ideal \(I\subseteq \E\) has a reduced standard basis. The standard basis is
unique modulo flat germs.
\end{thm}
The following Theorem along with \cite[Equation (2.2)]{GazorKazemiArxiv} and Buchberger algorithm provide computational guidelines on how to compute a reduced standard basis
in \(\E.\)
\begin{thm}(See \cite[Page 255]{Becker} and \cite[Theorem 2.10]{GazorKazemiArxiv})\label{m14}
Let \(I=\langle f_{1},\ldots, f_{n}\rangle_{\E}\).
\begin{itemize}
  \item There always exist germs \(g_{1},\ldots, g_{m}\in\E\) for \(m\leq n\) so that
\begin{equation}
I=\langle g_{1},\ldots, g_{m}\rangle_{\E}
\end{equation}
and \(\LT(g_{i})\nmid p\) for all \(p\in \Terms(g_{j})\), except for when \(i=j\) and \(p=\LT(g_{i})\).
  \item  Furthermore, assume that \(\{f_{i}\}^n_{i=1}\) is a standard basis in \(\E.\) Then, the set \(\{g_{i}\}^m_{i=1}\)  (for \(m\leq n\)) can be chosen
  so that it is a reduced standard basis.
\end{itemize}
\end{thm}
\begin{exm}\label{m12}
Let \(f\) be a non-zero flat germ, and \(I=\langle G\rangle_{\E}\) where
\bes G:=\big\lbrace g_{1}:=\lambda-\lambda \exp(x), g_{2}:=x-\sin(x), g_{3}:=\lambda x+\lambda^3+\lambda^{2} \ln(1+x)+f\big\rbrace. \ees
Since \(\LT(G)=\lbrace \lambda x, x^{3}\rbrace\) and \(\LT(I)=\langle \lambda x, x^{3}, \lambda^{3}\rangle_{\E}\), \(G\)
is not a standard basis with respect to \(\prec_{\alex}\). Due to
\bes \Rem\big(S(g_{1},g_{2}), G, \prec_{\alex}\big)=0, \quad \Rem\big(S(g_{2}, g_{3}), G, \prec_{\alex}\big)= -x^2 f\ees
and \(\Rem(S(g_{1},g_{3}), G, \prec_{\alex})=\lambda^{3}+f\), we add \(\lambda^{3}+f\) to the basis, i.e., \(S_{1}:=\lbrace g_{1}, g_{2}, g_{3},
\lambda^{3}+f\rbrace\). Now \(S_{1}\) is a standard basis, yet it is not a reduced standard basis. Given Theorem \ref{m14} and its proof,
\begin{eqnarray*}
&\LT(g_{1})\mid \LT(g_{3}), \quad \Rem\big(g_{3}, S_{1}\setminus\{g_3\}, \prec_{\alex}\big)= 0, \quad S_2:=\{g_{1}, g_{2}, \lambda^{3}+f\} \hbox{ and } & \\
&\Rem\big(\lambda x+g_{1}, S_{2}, \prec_{\alex}\big)=\Rem\big(g_{2}-\frac{x^3}{6}, S_{2}, \prec_{\alex}\big)=0,&
\end{eqnarray*} we define
\begin{equation*}
S_{3}:=\lbrace \lambda x, x^{3}, \lambda^{3}+f \rbrace.
\end{equation*}
Now \(S_{3}\) is a reduced standard basis.
\end{exm}

\subsection{Finite codimension}

Finite codimensional ideals demonstrate an important role in ideal presentation in Subsection \ref{ItrRepresent}. Let \(\mathcal{M}_\E:=\langle x, \lambda \rangle_{\E}\), \(\mathcal{M}_\RM=\langle x, \lambda\rangle_{\RM},\) \(\mathcal{M}_{K[[x, \lambda]]}=\langle x, \lambda \rangle_{K[[x, \lambda]]},\) and
\(\mathcal{M}_{K[x, \lambda]}=\langle x, \lambda \rangle_{K[x, \lambda]}\).

\begin{defn}\label{m4} An ideal \(I\) in \(\E\) (or in \(\RM,\) \(K[[x, \lambda]]\)) is said to have a
{\it finite codimension} when \({\M_\E}^{k}\subseteq I\) (\({\M_\RM}^k\subseteq I, {\M_{K[[x, \lambda]]}}^k\subseteq I\)) for \(k\in \mathbb{N}.\)
Equivalently, \(I\) has a complement (vector) subspace in \(\E\) with finite dimension.
\end{defn}

Now we compare the ideals in \(K[x, \lambda]\) with those in \(\RM,\) \(K[[x, \lambda]]\) and \(\E.\)
\begin{thm}(See \cite[Theorem 2.13]{GazorKazemiArxiv})\label{FcodRemainder}
Let \(I\) be a finite codimensional ideal in \(\E.\)
\begin{itemize}
\item[(1)] The ideal \(I\) has a unique reduced standard {\em polynomial} germ basis.
\item[(2)] Assume that \(G:=\{g_i\}^n_{i=1}\) is a standard basis for \(I\) and \(f\in \E.\) Then, the remainder is always a {\em polynomial germ}, \ie
\(\Rem(f, G, \prec_{\alex})\in K[x, \lambda].\) In particular, we have
\(\deg(\Rem(f, G, \prec_{\alex}))\leq k\) when \({\M_\E}^{k+1}\subseteq I.\)
\end{itemize}
\end{thm}
\begin{exm}
Let \(S= \{f_i\}^n_{i=1}\) be a reduced standard basis whose generated ideal is of
finite codimension. Then, \(\Rem(f, S, \prec_{\alex})= f-\sum q_i f_i\in K[x, \lambda].\) However, \(q_i\) may not always be a polynomial germ.
Consider the ordering \(\lambda\prec_{\alex}x\) and define
\be\label{specialexm} I:=\Le\langle g_{1}:=2\lambda^{3}-3\lambda^{2}x+x^{5}, g_{2}:=-3x\lambda^2+5x^{5}, g_{3}:=-3\lambda^{3}+5x^{4}\lambda\Ri\rangle.\ee
The ideal \(I\) is finite codimensional since  \({\M_\E}^{6}\subset I.\) Here,
\(G=\left\lbrace g_1, g_2, g_3, \frac{4}{3}x^{5}-\frac{10}{9}x^{4}\lambda \right\rbrace\) is a standard basis and
\bes S=\left\lbrace f_1:= x\lambda^{2}-\frac{25}{18}x^{4}\lambda, f_2:=-\frac{1}{3}g_{3}, f_3:=x^{5}-\frac{5}{6}x^{4}\lambda\right\rbrace\ees is a reduced
standard basis. Now we have
\bes x^3\lambda^3=x^{3} f_{2}+ \frac{25}{18} x^5q f_{1} +\frac{5}{3}\lambda x^{2}q f_{3},\ees where
\(q:= \sum^\infty_{i=0} (\frac{125}{108})^{i}x^{2i}\), \(\Rem(x^3\lambda^3, S, \prec_{\alex})=0\) while \(q_1:=\frac{25}{18}x^5q\) and
\(q_2:= \frac{5}{3}\lambda x^{2}q\) are not polynomial germs.
\end{exm}

The following theorem enables us to compute the standard basis for ideals in \(\E\) through computations in the fractional germs.
\begin{thm}\label{Thm2.17}
Consider \(I=\langle p_1,\ldots,p_n\rangle_{\E},\) \(\prec_{\alex}\) and \({\M_\RM}^{k}\subseteq\langle p_1, \ldots, p_n\rangle_{\RM}\subseteq \RM.\) Further,
suppose that \(\lbrace q_1,\ldots,q_m\rbrace\) is a standard basis for \(\langle p_1, \ldots, p_n\rangle_{\RM}.\) Then,
\(\lbrace q_1, \ldots, q_m\rbrace\) is a standard basis for \(I\).
\end{thm}
\bpr
Let \(g\in LT_\E(I)= \langle LT(f)\,|\, f\in I\rangle_\E.\) Since \({\M_\RM}^{k}\subseteq\langle p_1, \ldots, p_n\rangle_{\RM},\) part (1) in
Theorem \ref{Thm2.18} implies that \({\M_\E}^{k}\subseteq I\) and \({\M_\E}^{k}\subseteq LT_\E(I).\) Hence, without loss of generality we may assume that \(g=J^k(g).\) Thus,
\bes
g= \sum a_i LT(f_i)= J^k\sum a_i LT(f_i)= \sum J^k(a_i)LT(f_i)+ r,
\ees where \(a_i\in \E, f_i\in I\) and \(r\in {\M_\E}^{k+1}.\) Now we only need to prove that \(LT(f_i)\in \langle LT(q_j), j=1, \ldots m\rangle_\E\) when
\(\deg (LT(f_i))\leq k.\)

Since \(f_i\in I= \langle q_j, j=1, \ldots m\rangle_\E,\) \(f_i= \sum b_{ij} q_j\) for some \(b_{ij}\in \E.\) Thus,
\bes J^k(f_i)= J^k\sum J^k(b_{ij}) q_j= \sum J^k(b_{ij}) q_j+s
\ees where \(s\in {\M_{\RM}}^{k+1}.\) Therefore, \(J^{k}(f_i)\in \langle q_j, j=1, \ldots m\rangle_\RM.\) On the other hand,
\bas
LT(f_i) = LT(J^k(f_i))&\in& LT_\RM\langle q_1, \ldots, q_m\rangle_\RM
= \langle LT(q_1), \ldots, LT(q_m)\rangle_\RM
\\
&\subseteq& \langle LT(q_1), \ldots, LT(q_m)\rangle_\E.
\eas This completes the proof.
\epr
\section{Truncation degree and alternative ring computations}\label{TrunDegRing}

One of the main obstacles in working with the local rings is termination of the algorithms. There are methods in the literature for computing the {\it standard basis} (for fractional maps using Mora normal form) so that they solve the problem of algorithm terminations. Yet for the case of formal power series, no computer program can compute
and store infinite formal power series expansions and thus, their truncations up to certain degrees are mostly unavoidable. In order to circumvent this problem, we provide rigorous criteria so that the computations
can be performed modulo a sufficiently high degree; see Theorem \ref{Thm2.18a} (part 1). Further, we investigate
and compare the computations in the local (germ) rings of polynomials, the fractional germs and formal power series with those in the ring of smooth germs. We discuss the circumstances on which they can be alternatively used. This helps to efficiently use the different algorithms and yet ensure about correctness of the results.

The following theorems are two of our main contributions in this paper and provide important alternatives and criteria for our computations in different circumstances including ideals with infinite codimension.
\begin{thm}\label{Thm2.18}
\begin{itemize}
  \item[(1)] Suppose that \(G:=\{g_{1},\ldots,g_{n}\}\subset \RM.\) Then, the ideal \(\langle g_{1}, \ldots, g_{n}\rangle_{\E}\) has a finite
codimension iff \(\langle g_{1}, \ldots, g_{n}\rangle_{\RM}\) is a finite codimension ideal.
Assuming that \(\langle g_{1}, \ldots, g_{n}\rangle_{\RM}\) has a finite
codimension, for nonnegative integers \(i\) and \(j\), \({\M_\E}^{i}\langle \lambda^j\rangle_\E\subseteq \langle g_{1}, \ldots, g_{n}\rangle_{\E}\) iff
\({\M_{\RM}}^{i}\langle \lambda^j\rangle_\RM\subseteq \langle g_{1}, \ldots, g_{n}\rangle_{\RM}.\)
\item[(2)] Let \({\M_{K[x, \lambda]}}^{k}\subseteq\langle g_{1},\ldots, g_{n}\rangle_{K[x, \lambda]}\subseteq
K[x, \lambda].\) Then, \({\M_{K[x, \lambda]}}^{i}\langle \lambda^j\rangle_{K[x, \lambda]}\subseteq\langle g_{1},\ldots, g_{n}\rangle_{K[x, \lambda]}\) iff \({\M_\E}^{i}\langle \lambda^j\rangle_\E\subseteq\langle g_{1},\ldots, g_{n}\rangle_{\E}.\)

  \item[(3)] For a finite codimension ideal \(\langle g_{1}, \ldots, g_{n}\rangle_{\RM}\) in \(\RM\) and an ideal \(I\) in \(\E,\)
  \bes I= \langle g_{1}, \ldots, g_{n}\rangle_{\E} \quad \hbox{ iff }\quad  I\cap \RM= \langle g_{1}, \ldots, g_{n}\rangle_{\RM}.\ees
In particular, let \(I, J\) be two ideals in \(\E,\) \(I\) have a finite codimension and \(I\cap \RM= J \cap \RM.\) Then, \(I= J.\)
  \item[(4)] For an ideal \(I\) in \(\E,\) the following three conditions are equivalent.
  \bes
{\hbox{ (i)}}  \;\;   {\M_{\E}}^k\subseteq I, \qquad {\hbox{ (ii)}}\;\;   {\M_{\RM}}^k \subseteq I \cap \RM,  \qquad {\hbox{ (iii)}} \;\;
{\M_{K[x, \lambda]}}^k \subseteq I \cap K[x, \lambda].
  \ees
  \end{itemize}
  \end{thm}
  \bpr
  Part (1). The {\it if} part is trivial. Thus, we assume that \({\M_\E}^{k}\subseteq \langle g_{1}, \ldots, g_{n}\rangle_{\E}\) and prove that
\({\M_{\RM}}^{k}\subseteq
\langle g_{1}, \ldots, g_{n}\rangle_{\RM}\). Using Nakayama lemma \cite[Lemma 5.3, Page 71]{GolubitskySchaefferBook}, it is enough to verify that \({\M_\RM}^{k}\subset
\langle g_{1}, \ldots, g_{n}\rangle_{\RM}+{\M_\RM}^{k+1} \). Since \({\M_{\E}}^{k}\subset \langle g_{1}, \ldots, g_{n}\rangle_{\E}\), for some \(a_s\in \E\) we have
\begin{eqnarray*}
x^{l}\lambda^{k-l}&=&\sum_{s=1}^{n} a_{s}g_{s}=J^{k}\sum_{s=1}^{n}a_{s}g_{s}=\sum_{s=1}^{n}J^{k}\left(J^{k}(a_{s})g_{s}\right)= J^{k}\sum_{s=1}^{n}b_{s}g_{s}=\sum_{s=1}^{n}b_{s}g_{s}+r,
\end{eqnarray*}
where \(a_{s}\in\E, b_{s}:=J^k(a_s)\) and \(r\in{\M_{\RM}}^{k+1}.\) For the second claim, we do not need Nakayama lemma. We assume that \({\M_{\RM}}^{k}\subseteq
\langle g_{1}, \ldots, g_{n}\rangle_{\RM}.\) Similar to above, for \(i+j<k\) we have \(x^l\lambda^{i+j-l}= J^{k-1} (x^l\lambda^{i+j-l}) =\sum_{s=1}^{n}b_{s}g_{s}+r,\) where \(r\in {\M_{\RM}}^{k}\) and \(b_{s}\in \RM.\)

Part (2). Proof is similar to the proof of the second claim in part (1).

Part (3). Assume that \(I= \langle g_{1}, \ldots, g_{n}\rangle_{\E}.\) Trivially, \bes \langle g_{1}, \ldots, g_{n}\rangle_{\RM}\subseteq I \cap \langle g_{1}, \ldots, g_{n}\rangle_{\RM}.\ees Let \(f\in I \cap \RM.\) Thus, \(f= \sum a_i g_i\) for some \(a_i\in \E.\) So, \(f= J^k \sum J^{k} (a_i) g_i+ h\) for an \(h\in {\M_{K[x, \lambda]}}^{k+1}.\) Since the left hand side belongs to \(\langle g_{1}, \ldots, g_{n}\rangle_{\RM},\) the proof of the {\it if} part is complete.
Now assume that \(I\cap \RM= \langle g_{1}, \ldots, g_{n}\rangle_{\RM}\) and \(f\in I.\) By part (1) and \(\langle g_{1}, \ldots, g_{n}\rangle_{\E}\subseteq I,\) \(J^k(f)\in I\) for some \(k\in \mathbb{N}.\) Since \(J^k(f)\in I \cap \RM,\)
\bes
J^k(f)= \sum a_i g_i, \hbox{ for } a_i\in \RM.
\ees On the other hand \(f-J^k(f)\in \langle g_{1}, \ldots, g_{n}\rangle_{\E}.\) This completes the proof of part (3).
Part (4) is trivial.
  \epr
The hypothesis of part (1) in Theorem \ref{Thm2.18} allows us to use fractional germs for computations while the condition in part (2) permits the use of \Gr basis for our purposes despite the local nature of our results.

\begin{exm}\label{m15}
Part (1) from the previous theorem is not valid when \(\RM\) is replaced with \(K[x, \lambda]\). For instance, consider the example given
in \cite[Table 1, III.1 for \(k=5\)]{Keyfitz} and \cite[Page 77]{GolubitskySchaefferBook},
\bes I:=\langle g_1:=x^5+x^3\lambda+\lambda^2, g_2:=5x^5+3x^{3}\lambda, g_3:=5x^{4}\lambda
+3x^{2}\lambda^{2}\rangle_\E.\ees
The ideal \(I\) has a finite codimension since \({\M_{\E}}^{6}\subset I\). However, the ideal \(I_0:=\langle g_1, g_2, g_3 \rangle_{K[x, \lambda]}\) has an infinite codimension in \(K[x, \lambda]\). The reason for this is as follows. The reduced \Gr basis of \(I\) with respect to \(\lambda\prec_{\lex}x\) is given by
\bes G:=\big\lbrace 3125\lambda^{3}+108\lambda^{4}, 18\lambda^{3}+125\lambda^{2}x, 2x^{3}\lambda+5\lambda^{2}, 2x^{5}-3\lambda^{2}\big\rbrace.\ees
Further for any natural number \(n\geq 3,\) the remainder \(\Rem (x^n, G, \prec_{\lex})= c\lambda^3\) where
\(c\in K.\) Therefore, \(x^n\) does not
belong to \(I_0\) and \(I_0\) is an infinite codimensional ideal.
\end{exm}

\begin{thm}\label{Thm2.18a}
\begin{itemize}
  \item[(1)]
 Let \(g_{1}, \ldots, g_{n}\in\E\) and \(k\leqslant N,\) for \(k, N\in \mathbb{N}.\) Then,
\begin{itemize}
 \item[(a)] Either of the conditions
\bes {\M_\RM}^{k}\subseteq \langle J^{N}g_{1}, \ldots, J^{N}g_{n}\rangle_{\RM} \;\; \hbox{ and } \;\; {\M_{K[[x, \lambda]]}}^{k}\subseteq \langle J^{N}g_{1}, \ldots, J^{N}g_{n}\rangle_{K[[x, \lambda]]},
\ees is equivalent to
\({{\M_\E}}^{k}\subseteq \langle g_{1}, \ldots, g_{n}\rangle_{\E}.\)

\item[(b)] If \({\M_{K[x, \lambda]}}^{k}\subseteq\langle J^{N}g_{1},\ldots, J^{N}g_{n}\rangle_{K[x, \lambda]}\subseteq
K[x, \lambda].\) Then,
\bes
{\M_{K[x, \lambda]}}^{i}\langle \lambda^j\rangle_{K[x, \lambda]}\subseteq\langle J^{N}g_{1},\ldots, J^{N}g_{n}\rangle_{K[x, \lambda]} \;\; \hbox{ iff } \;\; {\M_\E}^{i}\langle \lambda^j\rangle_\E\subseteq\langle g_{1},\ldots, g_{n}\rangle_{\E}.
\ees
\end{itemize}
  \item[(2)] Consider the finite sequence of nonnegative integers \(k_i, l_i \in \mathbb{N}\) so that the sequence \(k_i+l_i\) is decreasing and \(l_i\) is increasing.
Let \(\langle f_1, f_2, \ldots, f_n\rangle_\E\) be an ideal that is not necessarily of finite codimension and
\be\label{InfEqu}
I:=\sum
{\M_\RM}^{k_i}\langle \lambda^{l_i}\rangle_{\RM}\subseteq \langle J^{N}f_{1}, \ldots, J^{N}f_{n}\rangle_{\RM}.
\ee Then, either of the following conditions
\begin{itemize}
  \item for each \(j\leq n,\) \(\Terms (f_j-J^{p_j}f_j)\subseteq \M_\RM I\) for some \(p_j\leq N.\)
  \item for each \(j\leq n,\) \(f_j\in \RM\) and \(f_j-J^{p_j}f_j\in \M_\RM I\) for some \(p_j\leq N.\)
\end{itemize} implies that
\be
\sum{\M_\E}^{k_i}\langle \lambda^{l_i}\rangle_{\E}\subseteq \langle f_1, f_2, \ldots, f_n\rangle_\E.
\ee
\end{itemize}
\end{thm}

\bpr
Part (1). Now the assumption \({\M_\RM}^{k}\subset \langle J^{N}g_{1}, \ldots, J^{N}g_{n}\rangle_{\RM}\) implies that for some \(a_{i}\in \RM,\) and any \(i=1, \ldots k,\) \begin{equation*}
x^{k-i}\lambda^{i}= J^k\sum_{i=1}^{n}a_{i}J^{N}(g_{i})= J^{k}\sum_{i=1}^{n}a_{i}g_{i}=
\sum_{i=1}^{n}a_{i}g_{i}+r, \quad \hbox{ for } \quad k\leqslant N,
\end{equation*}
where \(r\in{\M_\RM}^{k+1}\). This and Nakayama lemma conclude that \({{\M_\E}}^{k}\subseteq \langle g_{1}, \ldots, g_{n}\rangle_{\E}\). The converse and rest of the proof use similar arguments. Note that Nakayama lemma given in \cite[Lemma 5.3]{GolubitskySchaefferBook} is also true when \(\E\) is replaced with \(\RM\) and \(K[[x, \lambda]].\)

By the first part and \({\M_{K[x, \lambda]}}^{k}\subseteq\langle J^{N}g_{1},\ldots, J^{N}g_{n}\rangle_{K[x, \lambda]}\), we have \({\M_{\E}}^{k}\subseteq\langle g_{1},\ldots, g_{n}\rangle_{\E}\) and  \(\langle g_{1},\ldots, g_{n}\rangle_{\E}=\langle J^{N}g_{1},\ldots, J^{N}g_{n}\rangle_{\E}\). The rest of the proof is complete by Theorem \ref{Thm2.18}(2).

Since \({\M_{\E}}^{i}\langle \lambda^{j}\rangle_{\E}\subset \langle g_{1}, \ldots, g_{n}\rangle_{\E}\), for some \(a_s\in \E\) we have

\begin{eqnarray*}
x^{l}\lambda^{i+j-l}&=&\sum_{s=1}^{n}a_{s}g_{s}= J^{k-1}\sum_{s=1}^{n}a_{s}g_{s}=\sum_{s=1}^{n}J^{k-1}\left(J^{k-1}(a_{s})J^{N}(g_{s})\right)=J^{k-1}\sum_{s=1}^{n}b_{s}J^{N}(g_{s})\\&=&\sum_{s=1}^{n}b_{s}J^{N}(g_{s})+r,
\end{eqnarray*}
where \(b_{s}:=J^{k-1}(a_s)\) and \(r\in{\M_{K[x, \lambda]}}^{k}.\) Since \({\M_{K[x, \lambda]}}^{k}\subseteq\langle J^{N}g_{1},\ldots, J^{N}g_{n}\rangle_{K[x, \lambda]}\), from
the above equality one can conclude \({\M_{K[x, \lambda]}}^{i}\langle \lambda^j\rangle_{K[x, \lambda]}\subseteq\langle J^{N}g_{1},\ldots, J^{N}g_{n}\rangle_{K[x, \lambda]}\).

Part (2) is concluded by Nakayama lemma and the fact that
\bes
 x^{k_{i}-m}\lambda^{l_i+m}= \sum a_iJ^Nf_i= \sum a_if_i+r, \ees where \bes r\in\langle \Terms(f_j-J^{N}f_j)\rangle_\RM\subseteq \langle \Terms(f_j-J^{p_j}f_j)\rangle_\RM\subset \M_\RM I.\ees
\epr

\subsection{Computation}

The following items (I)-(XVI) provide our suggested algorithms that they are required for symbolic implementation of the results in a computer algebra system.
\begin{enumerate}[(I)]
\item {\bf Truncation degree and computational rings for ideals}

Depending on the type of input germs, our computations can be converted from \(\E\) to alternative smaller rings. This is feasible via Theorems \ref{Thm2.18} and \ref{Thm2.18a}. Now suppose that we are given an ideal \(\langle g_1, \ldots, g_n\rangle_{\E}\) with \(g_1,\ldots,g_n\in \E\). By multiplication matrix approach and standard basis computation, we can verify if the hypothesis of any part in
Theorems \ref{Thm2.18} and \ref{Thm2.18a} holds. In fact, from Theorem \ref{Thm2.18a}\((1)\)\ we find the permission for converting our computational ring to either of \(\RM\), \(K[[x, \lambda]]\), or \(K[x, \lambda]\).
For when \( g_1, \ldots, g_n\in \RM\), computations are converted to the computational rings \(\RM\) and/or \(K[x, \lambda]\) when part \((1)\) and/or \((2)\)\ in Theorem \ref{Thm2.18} holds, respectively. In the case that \(\langle g_1, \ldots, g_n\rangle_{\E}\) 
is not a finite codimensional ideal, we resort to Theorem \ref{Thm2.18a}\((2)\). These arguments give rise to the following procedure:
 \begin{proc}\label{TruncationAndRings}
 For the output computations regarding the ideal \(\langle g_1,\ldots,g_n \rangle_{\E}\) when
\begin{enumerate}
\item[(1)]
The inputs \(g_1,\ldots,g_n\in \E\), we have
\begin{enumerate}
\item Theorem \ref{Thm2.18a}(1) holds for \(k\leq 20\), return \(k\) as truncation degree and
\(K[x, \lambda]\), \(\RM\), \(K[[x, \lambda]]\), \(\E\) as permissible computational rings.
\item Theorem \ref{Thm2.18a}(1)(a) holds and \ref{Thm2.18a}(1)(b) fails for \(k\leq 20\),
we return \(k\) as truncation degree and \(\RM\), \(K[[x, \lambda]]\), \(\E\) as permissible computational rings.
\item Theorem \ref{Thm2.18a}(1) fails for \(k\leq 20\), the only permissible computational ring
is \(\E\) and no clear truncation degree is available. Yet Theorem \ref{Thm2.18a}(2) provides
a procedure for finding permissible truncation degree and the corresponding computational ring as \(\RM\) for computation of the intrinsic part of the ideal \(\langle g_1,\ldots,g_n \rangle_{\E}\).
\end{enumerate}
\item[(2)]
The inputs \(g_1,\ldots,g_n\in K[x, \lambda]\), we have
\begin{enumerate}
\item Theorem \ref{Thm2.18}(1)-(2) hold. Then, return \(K[x, \lambda]\), \(\RM\), \(\E\) as
 permissible computational rings and the maximal number \(k\) in part (2) is the permissible
 truncation degree.
\item Theorem \ref{Thm2.18}(1) holds, \ie
\be\label{eq}
\exists\,\, {\M_\RM}^{k}\subseteq \langle g_{1}, \ldots, g_{n}\rangle_{\RM}\quad \hbox{and}\quad k\leq 20.
\ee

and \ref{Thm2.18}(2) fails for \(k\). In this case, the only permissible computational rings are
\(\RM\) and \(\E\) while minimal number \(k\) in Equation \eqref{eq} represents the permissible
truncation degree.
\item Theorem \ref{Thm2.18}(1)-(2) fail. Since the inputs are polynomials, the truncation is
not needed. Hence, the number \(N\) in Equation \eqref{InfEqu} is taken as the {\tt Singularity}'s default upper bound \(N=20\) or the alternative choice given by the informative user.
\end{enumerate}
\end{enumerate}
\end{proc}
\suspend{enumerate}
\begin{rem}
When the codimension of a singular germ is too large or infinity, our algorithms need an upper bound for termination of computations. Therefore, the default upper bound for
computation of a singular germ is set to be 20. The truncation degree option in each procedure
provides our users to increase or decrease this default number; see \cite{GazorKazemiUser} for further details.
\end{rem}
\resume{enumerate}[{[(I)]}]
\item {\bf Truncation degree and computational rings for a smooth germ}

For computations regarding a singular germ when there is no specified input data from the user, there is a necessity for an algorithmic derivation of suitable truncation degree and computational rings. This procedure is given as follows.
\begin{proc}\label{TruncationAndRingGerm}
For an input germ \(g\in \E\) and the required output truncation degree \(k\) and computational ring \(\R\):
\begin{itemize}
\item Define the ideal \(I(g) := \langle xg, \lambda g, x^2 g_{x}, \lambda g_{x}\rangle_{\E}\); compare \(I(g)\) with the high order term ideal defined by equation \eqref{pg}. Next, follow the procedure \ref{TruncationAndRings} to find the truncation degree \(k\) and computational ring \(\R\) for the ideal \(I(g)\). Return \(k\) and \(\R.\)
\end{itemize}
\end{proc}

\item {\bf Remainder}

General smooth inputs (not necessarily polynomials) indeed are important
for many applications in the modeling and analysis of real life problems. Division algorithms in computers for general computational rings \(\E\) and \(K[[x, \lambda]]\) seem impractical in general. In this regard, our implementation of division algorithm for its
applications in bifurcation theory is feasible given Definition \ref{m11}, Theorems \ref{FcodRemainder}--\ref{Thm2.18a}), and in particular Procedures \ref{TruncationAndRings}--\ref{TruncationAndRingGerm}. This is one of our claimed novel contributions in this paper. Division algorithm computation of a smooth germ \(g\) divided by a set of smooth germs \(\lbrace f_1, \ldots, f_n\rbrace\) for the remainder is a fundamental tool in {\tt Singularity}. Our implementation adapts the algorithms based on three computational rings \(\RM\), \(K[[x, \lambda]]\) and \(\E\). The inputs \(g\) and \(\lbrace f_1, \ldots, f_n\rbrace\) can be arbitrarily chosen
from either of these rings. The division algorithm in \(\RM\) is known as Mora normal form
 \cite[Page 170]{CoxUsing}.
\item {\bf Standard Bases}

Standard bases are the basis of most computations in {\tt Singularity} in four computational rings \(K[x, \lambda]\), \(\RM\), \(K[[x, \lambda]]\), and \(\E\); see Definitions \ref{SBasis} and \ref{rsgb}, and Theorems \ref{m10}, \ref{m1}, \ref{m14}, \ref{FcodRemainder}(1), \ref{Thm2.17}, \ref{Thm2.18}, \ref{Thm2.18a}.
\section{Intrinsic ideal representation}\label{ItrRepresent}

An elegant approach for answering the ideal membership problem is to provide a good presentation for ideals.
In this section we define intrinsic ideals and use it for such representation. This plays a central role in developing \texttt{Singularity} and is a prerequisite for most computations of objects presented in Section \ref{Sec4}.

\begin{defn}
When only the trivial reparametrization \(\Lambda(\lambda)=\lambda\) is allowed in Equation \eqref{ContEqvl}, the associated relation \(\thicksim_s\)
is called {\em strongly equivalent relation}; see \cite[Page 51]{GolubitskySchaefferBook}. Then, an ideal \(I\) is called {\em intrinsic} \cite[Page 81]{GolubitskySchaefferBook} when \(I\)
includes all strongly equivalent classes of its members.
\end{defn}

By \cite[Proposition 7.1]{GolubitskySchaefferBook}, every finite codimension ideal $I$ in \(\E\) is intrinsic if and only if there exist nonnegative integers
\(s, m_i, n_i\) for \(i=0, \ldots, s\) so that
\be\label{Intrinsic}
I=\sum^s_{i=0} {\M_\E}^{m_i}\langle\lambda^{n_{i}}\rangle_\E,
\ee \(n_0=0,\) and the sequence \(n_i\) is strictly increasing while \(m_i+n_i\) is strictly decreasing. The conditions on \(m_i\) and \(n_i\) make the
representation \eqref{Intrinsic} unique. Equation \eqref{Intrinsic} for intrinsic ideals gives a convenient answer for the ideal membership problem.
The monomials \(x^{m_i}\lambda^{n_i}\) for \(i=0, \ldots, s\) are called {\it intrinsic generators} of the intrinsic ideal \(I.\)
For non-intrinsic ideals or more generally for a vector space \(I,\) we define their intrinsic part; \ie we denote \({\rm Itr}(I)\) for the largest
intrinsic ideal contained in \(I\).

\begin{lem}\label{LemIntDecom} For an intrinsic ideal \(I\), there always exists a reduced monomial standard basis for \(I\) that includes its
intrinsic generators.
Let \(I\) be a finite codimension vector subspace of \(\E\) or be a finitely generated ideal. Then, there exist nonnegative integers \(m\) and \(n,\)
a reduced monomial standard basis \(\{f_j\}^n_{j=1}\) for \(\Itr(I)\) and \(\{g_i\}^m_{i=1}\subset\E\) so that
\be\label{IntDecom} I= \Itr(I)+ \langle\{g_i\}^m_{i=1}\rangle_\R,\ee
where \(\R:= \E\) or \(\R:= K.\) Here, none of the terms in \(\Terms(g_i)\) is divisible by \(f_j,\) \ie \(\Itr(I)\cap (\cup^m_{i=0}\Terms(g_i))= \emptyset.\)
\end{lem}

\bpr
Let the intrinsic ideal \(I\) be given by \eqref{Intrinsic} and \(G\) be the set of its intrinsic generators. Then for \(i=0, \ldots, s,\)
we consecutively update \(G\) by adding the monomials of type \(x^\alpha\lambda^\beta\) to \(G\) where \(\alpha+\beta=m_i+n_i,\) \(\beta\geq n_i,\) and
\(x^\alpha\lambda^\beta\) is not divisible by the elements of \(G.\) This gives rise to a reduced standard monomial basis for \(I.\) For the second part,
we divide the generators of \(I\) by the reduced standard monomial basis of \(\Itr(I)\) and define the nonzero remainders as \(g_i.\)
\epr

A refinement of the decomposition \eqref{IntDecom} for finite codimension ideals is given as follows. We denote \(I^\perp\) for the set of all monomials
not in \(I,\) while \(\langle I^\perp\rangle_K\) stands for the vector space generated by \(I^\perp\). Now by \cite[Corollary 7.4]{GolubitskySchaefferBook}, we have
\be\label{IdealInt}
I=\Itr(I)\oplus \Le(I\cap \langle {\Itr(I)}^{\perp}\rangle_K\Ri).
\ee
\subsection*{Algorithms for intrinsic part computations}

Here we describe our suggested algorithm on how to compute the intrinsic ideal representation \eqref{Intrinsic}
for finite codimensional ideals.
\begin{itemize}
  \item
The first step is to find a lower and upper bound for the maximum natural number \(k\) so that
\({\M_\E}^{k}\subseteq I\lhd \E.\)
\end{itemize}

For any ideal \(J\) in a ring \(\R\) we define
\begin{equation}\label{R/J}
\varphi_{u, J}: \frac{\R}{J}\longrightarrow\frac{\R}{J}, \quad \varphi_{u, J}(f+J):=  uf+J,
\end{equation} for \(u\in\lbrace x, \lambda\rbrace.\) Obviously for
\bes J:=I \hbox{ and } \R:= \E,\ees \(\E/I\) is a finite dimensional vector space and \(\varphi_{u, I}\) is a linear nilpotent map.

Now we explain how to derive \(N_u\). Computing a (reduced) standard basis for \(I,\) say \(S=\lbrace g_i, i=1, \ldots n\rbrace,\)
we now introduce a vector space basis generator
\be\label{Iperp} B=\big\{w\,|\, w \hbox{ is a monomial, } w\notin\langle LT(g_i)\,|\,g_i\in S\rangle_\R\big\}. \ee
for \(\E/I\); also see \cite[Proposition 4, Pages 177--179]{CoxLittleIdeals} and \cite[Pages 128--129]{Armbruster}.

\cite[Lemma 3.3]{GazorKazemiArxiv} proves that for a finite codimension ideal \(I\) in \(\E\) and a standard basis \(S\) for \(I\), \(B+I\) is a \(K\)-vector space basis for \(\E/I\). This
readily provides a matrix representation
\be\label{E/ILem}
M_{u, I} := [\varphi_{u, I}]_{B, B}
\ee
for \(\varphi_{u, I}.\)
Since the nilpotent degree \(N_u\) for  \(\varphi_{u, I}\) is less than \(k\), we have
\bes
\max \{N_x, N_\lambda\}\leq k< N_{x}+N_{\lambda}.
\ees
Hence,
the natural number \(N,\)
satisfying \(\max \{N_x, N_\lambda\}\leq N\leq N_x+N_\lambda-1,\) is consecutively increased to obtain \(k=m_0.\) In fact the remainders of
\(x^i\lambda^{N-i}\) divided by the standard basis \(S\) conclude the result, thanks to Theorem \ref{m10} (part d).

\begin{itemize}
  \item Next we look for the maximum values of $m_i$ and \(n_i\) such that \({\M_\E}^{m_i}\langle \lambda^{n_i}\rangle_\E\subseteq I.\)
\end{itemize}
Our suggestion uses
the concept of colon ideals. Let \(J_1\) and \(J_2\) be two ideals in \(\E\). We define the {\it colon ideal} of \(J_1\) by \(J_2\)
(see \cite[Definition 5, Page 194]{CoxLittleIdeals}) as
\begin{equation*}
J_1:J_2=\langle f\in \E\mid fJ_2\subseteq J_1\rangle_\E.
\end{equation*}

By an inductive procedure and repeating the above for the colon ideal \(I: \langle \lambda^n\rangle_\E\) for each \(n\), we obtain \(m_i\) and \(n_i\) as desired.
Note that \(n_{i-1}< n_i\) and \(m_i+n_i< m_{i-1}+n_{i-1}<m_0\) for any \(i\). Therefore, the only remaining challenge is to compute the maximum
value \(m\) so that
\bes
{\M_\E}^m\subseteq I: \langle \lambda^n\rangle_\E\,\, \hbox{ or equivalently},\,\, {\M_\E}^m\langle \lambda^n\rangle_\E\subseteq I.
\ees
The procedure mentioned above by using \(J:= I: \langle \lambda^{n}\rangle_\E\) in Equation \eqref{R/J}
provides suitable lower and upper bound for \(m.\) The following lemma and its' follow up comments facilitate the computations of a finite generating set for colon ideals and next, their (reduced)
standard basis readily determines \(m.\)

\begin{lem}\label{ColonLem}
Suppose that \(I\) is a finite codimensional ideal in \(\E\), \(J:=I\cap K[x, \lambda]\) and \(J=\langle
 g_{1}, \ldots, g_{n}\rangle_{K[x, \lambda]}\). Then, \(I=\langle g_{1}, \ldots, g_{n}\rangle_{\E}.\)
Further, let \(J\cap \langle g\rangle_{K[x, \lambda]}=\langle h_{1}, \ldots, h_{n}\rangle_{K[x, \lambda]}\) for a monomial germ \(g\in K[x, \lambda]\). Then,
\be\label{Colon}
I:\langle g \rangle_{\E}=\left\langle \frac{h_1}{g}, \ldots, \frac{h_n}{g} \right\rangle_{\E}.
\ee
\end{lem}
\bpr
Let \(f\in I\). Since \(I\) has a finite codimension, \({\M_\E}^{k+1}\subseteq I\) for some \(k\).
Hence, \(J^k(f)\in I\cap K[x, \lambda]=J.\) Theorem \ref{Thm2.18} (part 3) completes the proof of the first part.

For any \(u\in I:\langle g\rangle_{\E}\), we have \(u\langle g\rangle_{\E}\subseteq I\). Thus,
\be\label{m9}
ug\in I\cap \langle g\rangle_{\E}.
\ee For any \(p\in I\cap \langle g\rangle_{\E}\cap K[x, \lambda],\) \(p= fg= J^{\deg p- \deg g}(f) g\) for some \(f\in \E.\) Thereby,
\bes
I\cap \langle g\rangle_{\E}\cap K[x, \lambda]= J\cap \langle g\rangle_{K[x, \lambda]}=\langle h_{1}, \ldots, h_{n}\rangle_{K[x, \lambda]}.
\ees
By the first part, we have
\bes
I\cap \langle g\rangle_{\E}=\langle h_{1}, \ldots, h_{n}\rangle_{\E}.
\ees
This and Equation \eqref{m9} imply that \(u\in\langle \frac{h_1}{g},\ldots, \frac{h_n}{g}\rangle_{\E}\) and \(I:\langle g\rangle_{\E}\subseteq
\langle \frac{h_1}{g},\ldots, \frac{h_n}{g}\rangle_{\E}\).

Now let \(u\in\langle \frac{h_1}{g},\ldots, \frac{h_n}{g}\rangle_{\E}\). Hence,
\(aug\in \langle h_1,\ldots, h_n\rangle_{\E}\) for any \(a\in \E\). On the other hand, \(\langle h_1,\ldots, h_n\rangle_{\E}\subseteq I\) shows that \(aug\in I\) and eventually \(u\in I:\langle g\rangle_{\E}\).
\epr
Lemma \ref{ColonLem} should be compared with its analogues in \cite[Theorem 11, Page 196]{CoxLittleIdeals} and \cite[Theorem 5.5, Page 185]{CoxUsing}. Indeed, our main contribution here is to justify the use of the polynomial germ ring instead of working with naive and alternative choice \(I\cap \langle g\rangle_\E.\) Remark that we need to compute \(J\cap \langle g\rangle_{K[x, \lambda]}\) for the monomials \(g:=\lambda^n.\) The benefit of polynomial germ ring is due to the existence of lexicographic ordering. In fact,
the computation of \(J\cap \langle g\rangle_{K[x, \lambda]}\) using the lexicographic ordering is efficient and classic (described below) while such an approach for \(I\cap \langle g\rangle_\E\) does not seem easy to work with standard basis with local orderings.

Finally, we recall the classical approach on the intersection computation of two ideals \(J_1\) and \(J_2\) in the polynomial germ ring \(K[x, \lambda]\);
see \cite[Page 187]{CoxLittleIdeals} for more details. Let \(J_1:= \la f_i \,|\, i=1, \ldots, n_1 \ra_{K[x, \lambda]}\) and
\(J_2:=\la g_i\,|\, i=1, \ldots, n_2 \ra_{K[x, \lambda]}.\) We define
\be
J_3:= \big\la t f_i, (1-t)g_j\,|\, \hbox{ for } i= 1, \ldots, n_1\hbox{ and } j=1, \ldots, n_2 \big\ra_{K[x, \lambda, t]}.
\ee Thus, \(J_3\cap K[x, \lambda]= J_1\cap J_2.\) We may compute a \Gr basis for \(J_3\) with respect to \(\prec_{\lex}\) (where \(x, \lambda\prec_{\lex} t\))
\bes G_3:= \big\{p_i(x, \lambda), q_j(x, \lambda, t)\,|\, \hbox{ for } i= 1, \ldots, n\hbox{ and } j=1, \ldots m\big\},\ees
here \(p_i\) represents those basis elements independent of \(t\) while \(q_j\) represents those explicitly depending on \(t.\) Since
\(\la LT(G_3)\ra_{K[x, \lambda, t]}= LT(J_3),\) we may conclude that
\be
G:= \big\{p_i(x, \lambda)\,|\, \hbox{ for } i= 1, \ldots, n\big\}
\ee is a \Gr basis for \(J_1\cap J_2.\) The reason is as follows. For any \(u\in LT(J_1\cap J_2)\subseteq LT(J_3)\), either \(LT(p_i)\)
or \(LT(q_j)\) for some \(i\) or \(j\) must divide \(u.\) Since \(t\) divides the leading term of \(q_j\) for any \(j,\) \(u\in \la LT(G)\ra_{K[x, \lambda]}.\)

\begin{rem}
Applications of the above procedure is when \(I\) is given by a finite set of generators \(\{g_i\}^n_{i=1}.\) Thus, it is essentially
useful (when it is possible) to instead use a truncated Taylor expansion. This is justified when the hypothesis of Theorem \ref{Thm2.18a}
(part 1) is satisfied. Parts (1) and (3) in Theorem \ref{Thm2.18}  indicate that for finite codimension ideals in \(\E,\) we may simply replace \(\R:= \E\) with either \(\RM\) or \(K[[x, \lambda]].\) Part (2) in Theorem \ref{Thm2.18a} provides an important computational criterion for ideals with infinite codimension in \(\E\).
\end{rem}
\item {\bf Multiplication matrices}

The aim of this tool is to find \(m_{i}\) for \(i=0,\ldots,s\) in Equation \eqref{Intrinsic} such that \(\M^{m_i}_{\E}\langle \lambda^{n_i}\rangle_{\E}\subset I\). Indeed, the degree of minimal polynomial associated with multiplication matrix of a given variable (either state variable or distinguished parameter) is the starting candidate in this detection. The algorithm is described through Equations \eqref{R/J}, \eqref{Iperp}, and \eqref{E/ILem}. This plays the key role in {\tt Verify} function in {\tt Singularity} to find permissible truncation degree and computational rings as well as the
estimations in the transformation functions (\(X(x, \lambda)\) and \(S(x, \lambda)\) in the Equation \eqref{ContEqvl}) transforming two equivalent singular germs
to each other. These transformations play the essential role in the analysis of real life problems. In fact, bifurcations arising
from a parametric normal form can be located in terms of the parameters of the original problem
through these transformations.
\item {\bf Colon ideals}

Colon ideals are one of the two main (together with Multiplication matrices) tools which enables us to construct intrinsic ideals (Equation \eqref{Intrinsic}), in particular, finding the maximal ideal terms of the form \(\M^{m_i}_{\E}\langle \lambda^{n_i}\rangle_{\E}\). Lemma  \ref{ColonLem} provides the criteria for the execution of the {\tt Intrinsic} command in {\tt Singularity}. In fact,  computations of colon ideals are converted from \(\E\) to \(K[x, \lambda]\) as long as those criteria hold. Hence, we take help from the efficient built-in {\sc Maple} command {\tt Intersect} and handle
the computation of colon ideals in \(\E\).

\item {\bf Intrinsic part of ideals and vector spaces}\label{intrinsicpartidealvector}

Intrinsic ideals are the main tools for a comprehensive representation of ideals. The command
{\tt Intrinsic} is designed for these purposes through several options for the computation of
intrinsic parts of finite and infinite codimension ideals and vector spaces. It is moreover used for computing intrinsic generators through the command {\tt IntrinsicGen} in our library. The command {\tt IntrinsicGen} is needed in the commands {\tt Normalform} and {\tt RecognitionProblem}.
\begin{itemize}
\item[(a)]
The {\tt Intrinsic} command computes the intrinsic part of a finite codimension ideal.
\item[(b)]\label{intrinsicb}
To obtain the intrinsic part of a vector space \(I+V\), where  \(I=\langle f_1,\ldots,f_n \rangle_{\E}\) and \(V=\mathbb{R}\lbrace g_1,\ldots, g_s\rbrace\), we first find \(Itr(I)\), that is the intrinsic part of \(I\). Following Lemma
\ref{LemIntDecom}, there exists a reduced monomial standard basis \(\lbrace p_{i}\rbrace\) for
\(Itr(I)\). Then, \(r_{i}:=\Rem(f_i,\lbrace p_i\rbrace,\prec)\) and \(\Rem(g_j,\lbrace \lbrace p_i\rbrace, \lbrace r_{i}\rbrace_{i=1}^{n}\rbrace,\prec)\) are replaced with \(f_{i}\) and \(g_{j}\), respectively. In the enlarging process of \(Itr(I)\) to \(Itr(I+V)\) the membership of elements \(x^{m}\lambda^{n}\) in \(I+V\) is checked by computing \(\Rem(x^{m}\lambda^{n},\lbrace f_i\rbrace,\prec)=r\) and finally verifying if \(r\in {\rm span}_{\mathbb{R}}\lbrace g_{j}\rbrace\).
\end{itemize}
\item {\bf Vector space complement computation for \(Itr(I)\) in \(I\) and for \(I\) in \(\E,\) when \(I\) is a finite codimension vector space or an ideal}\label{VectorComplement}
\begin{itemize}
\item[(a)] The classical basis complement space computation of a finite codimension ideal is well-known and in our context follows
Equation \eqref{Iperp}. This is implemented through the command {\tt Normalset} in {\tt Singularity}. This procedure is used in the computations of normal form, universal unfolding, and recognition problem.

\item[(b)] For a finite codimension vector subspace of \(\E\) or a finitely generated ideal, the algorithm follows Lemma \ref{LemIntDecom} while for the finite codimension ideals the computation follows Equation \eqref{IdealInt} that is \cite[Corollary 7.4]{GolubitskySchaefferBook}.

\item[(c)] The first step for complement space computation of a finite codimension vector space \(I\) in \(\E\) is to follow Equation \eqref{IntDecom}, that gives rise to the decomposition \(I=\langle f_1, \ldots, f_n\rangle_{\E}+\mathbb{R}\lbrace g_1, \ldots, g_m\rbrace\). Since \(\{f_i| i=1, \ldots n\}\) is a reduced standard basis, we apply the procedure described by Equation \eqref{Iperp} and find the monomial generators of \(\E/\langle f_1, \ldots, f_n\rangle_{\E}\), say \(\lbrace p_1, \ldots, p_k\rbrace\). Then, our algorithm initially assumes \(W:=\emptyset.\) Next we inductively from \(i=1\) to \(k\) look for the monomials \(p_{i}\not\in\mathbb{R}\lbrace g_1, \ldots, g_m, p_1, \ldots, p_{i-1}\rbrace\). All such monomials \(p_i\)-s are added to the updating complement space \(W\), \ie  \(W:=W\cup\lbrace p_i \rbrace\). The final set \(W\) is the desired complement space.
\end{itemize}
\section{Computations of objects in bifurcation theory}\label{Sec4}

In this section we recall the algebraic tools and present our suggested approaches that are needed for computation of {\it normal forms},
{\it universal unfolding,} and {\it persistent bifurcation diagram classification.}

\subsection{Normal form}

Given a singular germ \(g\in \E\), from \cite[Pages 88--89]{GolubitskySchaefferBook} we recall the intrinsic ideals \(\p(g)\) and \(\s(g)\) by
\be \label{pg}
\p(g):= \Itr\left(\langle xg,\lambda g,x^{2}g_{x},\lambda g_{x}\rangle_\E\right),
\ee and
\be \label{sg}
\s(g):=\Sigma_{(m_i, n_i)} {\M_\E}^{m_i}\langle\lambda^{n_i}\rangle_\E,
\ee where \(\frac{\partial^{m_i}}{\partial x^{m_i}}\frac{\partial^{n_i}}{\partial \lambda^{n_i}}(g)\, (0, 0)\neq0,\) and there would not exist
nonnegative integers \(p\) and \(q\) such that \(\frac{\partial^{p}}{\partial x^{p}}\frac{\partial^{q}}{\partial \lambda^{q}} (g)\,(0,0)\neq0,\)
\(q\leq n_i,\) \(p+q\leq m_i+n_i.\) The extra restrictions on \(m_i\) and \(n_i\) here make the presentation \eqref{sg} unique as of those
in Equation \eqref{Intrinsic}.

\cite[Proposition 8.6]{GolubitskySchaefferBook} indicates that for any germ \(p\in \p,\) \(g\pm p\) is contact-equivalent to \(g.\) Terms in
\(\p(g)\) are called {\it high order terms}; see \cite[Page 89]{GolubitskySchaefferBook}. Therefore, we may stay contact-equivalent to \(g\) by
removing all terms in \(\Terms(g)\cap \p(g)\) from Taylor expansion of \(g.\) Therefore, \(\p(g)\) has a finite codimension if and only if \(g\)
is {\it finitely determined}, \ie \(g\) is contact-equivalent to a polynomial germ.

\cite[Theorems 8.3 and 8.4]{GolubitskySchaefferBook} state that for any term \(x^m\lambda^n\) in \(\s(g)^\perp,\) we must have
\bes {\frac{\partial^{m}}{\partial x^{m}}}{\frac{\partial^{n}}{\partial \lambda^{n}}}(g)\, (0, 0)=0, \quad \hbox{ while }\quad
{\frac{\partial^{m}}{\partial x^{m}}}{\frac{\partial^{n}}{\partial \lambda^{n}}}(g)\, (0, 0)\neq 0\ees for any intrinsic generator
\(x^{m}\lambda^{n}\) in \(\s(g).\) Now we recall the {\it intermediate order terms} as terms belonging to
\bes
A:=\p(g)^\perp\setminus \Big(\s^\perp(g)\bigcup \big\{x^{m_i}\lambda^{n_i}\,|\, x^{m_i}\lambda^{n_i} \hbox{ is an intrinsic generator for } \s(g)\big\}\Big).
\ees

\subsubsection{Normal form computation}

Now we are ready to provide an algorithm for computing {\it normal form} of a given germ \(g.\)
Using the procedure given in Section \ref{ItrRepresent}, we may compute \(\p(g)\) and remove all terms in \(\Terms(g)\cap \p(g)\) from Taylor
expansion of \(g\) to obtain a more simplified contact-equivalent germ, say \(f\). Now it only remains to eliminate intermediate order terms \(A\)
from \(f\) as many as possible. Then, this gives rise to its {\it normal form}.

When \(A\) is empty, \(f\) will be called {\it normal form} of \(g.\) Otherwise, we may use suitable polynomial change of variable \(X(x, \lambda)\)
and positive polynomial germ \(S(x, \lambda)\) to eliminate some intermediate order terms from \(f\). For example, we may replace \(x\) by \(ax+b\lambda\)
in \(f\) where \(a>0\) and \(b\) are arbitrary constant coefficients. This gives rise to a system of linear equations and a maximal solvable subsystem
leads to the further elimination of negligible terms in \(f\) and thus, the normal form computation of \(g.\) This approach needs to be systematically adapted along with standard basis computations for high order terms in multi-state variable cases. This is because complete algebraic characterization for high order terms is not yet known for many multi-state variable cases. This will be addressed 
in our upcoming result for multi-state variable cases.

\item {\bf High order terms ideal}\label{highorder}

This ideal is given by Equation \eqref{pg} and its computation follows \ref{intrinsicpartidealvector}.
This command is required
for deriving normal forms, transformations transforming equivalent germs, and optimal truncation degree.
\item {\bf Smallest intrinsic ideal containing a germ}\label{S}

This ideal is given by Equation \eqref{sg} and is performed in {\tt Singularity} by the command {\tt S}. This command is needed for the computation of  normal forms and recognition problem (described below). The computation starts with the derivation of intrinsic generators and then a generating set for the smallest intrinsic ideal. Finally, the procedure described by Lemma \ref{LemIntDecom} gives rise to a reduced monomial standard basis. The latter prevents redundant terms in intrinsic representation.

\subsection{Universal unfolding}

In this section we recall the algebraic formulation needed for computation of the universal unfolding of a singular germ \(g.\)

The restricted tangent space ideal of \(g\) is defined by \(RT(g):=\langle g, xg_{x}, \lambda g_{x}\rangle_{\E}\) while the tangent space \(T(g)\) is defined by
\be\label{TangSpace}
T(g):= RT(g) \oplus K\{g_{x}, g_{\lambda},\lambda g_{\lambda}, \ldots, \lambda^{\ell}g_{\lambda}\},
\ee for sufficiently large \(\ell.\) For a finite codimension \(g,\) there exists a natural number \(\ell\) so that
\(\lambda^{\ell}g_{\lambda}\notin \langle g, xg_{x}, \lambda g_{x}\rangle_{\E}\) and
\(\lambda^{l}g_{\lambda}\in \langle g, xg_{x}, \lambda g_{x}\rangle_{\E}\) for \(l>\ell\); see \cite[Page 127]{GolubitskySchaefferBook}.

Next, a universal unfolding of \(g\) is defined by
\be\label{UniUnfold} G(x, \lambda, \mu):= g(x, \lambda)+\sum_{i=1}^{k}\alpha_{i} p_{i}(x, \lambda),\ee
where \(\alpha:=(\alpha_1, \ldots, \alpha_k),\) $p_{i}$-s form a basis for a complement space of \(T(g).\) Thus, we may choose \(p_i\in T(g)^\perp.\)
The number \(k\) is called {\it codimension} of \(T(g)\) or equivalently {\it codimension} of \(g\).

The following describes the algorithms needed for computation of universal unfolding.
\item {\bf Tangent space}\label{tangentspace}

Procedure \ref{intrinsicpartidealvector} together with Procedure \ref{VectorComplement}(b) yield a complete description for restricted tanget space computation.

The tangent space \(T(g)\) is a vector space defined by Equation \eqref{TangSpace}. The algorithm initially uses the procedure described in \ref{intrinsicpartidealvector}(b) in order
to derive the maximal intrinsic ideal subset of \(T(g)\) for a given singular germ \(g\). Once \(\Itr(T(g))\) is computed, Procedure \ref{VectorComplement}(b) derives the complement
subspace \(W\) of \(\Itr(T(g))\) in \(T(g)\) and thus, \(T(g)\) is given by
\be\label{decompose}
\Itr(T(g))\oplus W.
\ee
This task is performed in {\tt Singularity} with the command {\tt T}.

\item {\bf Universal unfolding}

The universal unfolding of a singular germ \(g\) is given by Equation \eqref{UniUnfold}
where $p_{i}$-s constitute a vector monomial basis for the complement space of
\(T(g)\). Thus, Procedure \ref{tangentspace} is first used to obtain \(\Itr(T(g))\) and \(W\) in Equation \eqref{decompose}, that is, a representation of the form described in Equation \eqref{IntDecom}. Therefore, the hypothesis of Procedure \ref{VectorComplement}(c) is satisfied. Hence, Procedure \ref{VectorComplement}(c) derives the monomial basis \(p_i\) and finally the universal unfolding \eqref{UniUnfold} is computed. This is available through the command {\tt UniversalUnfolding} in our library. The {\tt list} option can be used with the {\tt UniversalUnfolding} command to return a list of alternative universal unfoldings for an input singular germ \(g\); see \cite [Page 7]{GazorKazemiUser}. This list is generated by changing the ordering of \(x\) and \(\lambda\) in the standard basis computations.

\subsection{Persistent bifurcation diagram classification} \label{DiagramClass}

Given our description in Section \ref{Sec1Int}, persistent bifurcation diagram classification is complete by simplifying the defining equations of transition
set \(\Sigma\) and then, choosing one parameter from each connected component of \(\Sigma^c\); see \cite[Page 140]{GolubitskySchaefferBook}.
The latter is enabled in {\tt Singularity} by using the Maple package {\tt RegularChains}. Let \(g\) be a singular germ of finite codimension and \(G(x, \lambda, \alpha), \alpha\in \mathbb{R}^k,\) be a universal unfolding of the germ \(g.\)
We recall that the transition set
\(\Sigma:= \mathscr{B}\cup \mathscr{H}\cup \mathscr{D},\) where
\ba\nonumber
\mathscr{B}&:=&\{\alpha\in\mathbb{R}^{k}\mid G=G_{x}=G_{\lambda}=0 \hbox{ at }  (x,\lambda,\alpha)
\hbox { for some }(x, \lambda)\in \mathbb{R}\times \mathbb{R}\},\\\label{Trans}
\mathscr{H}&:=&\{\alpha\in\mathbb{R}^{k}\mid G=G_{x}=G_{xx}=0 \hbox{ at } (x, \lambda, \alpha) \hbox{ for some }(x, \lambda)\in \mathbb{R}\times \mathbb{R}\},
\\\nonumber
\mathscr{D}&:=& \{\alpha\in\mathbb{R}^{k}\mid G =G_{x}=0 \hbox{ at } (x_{i}, \lambda, \alpha) \hbox{ for } i=1, 2 \hbox{ and } x_{1}\neq x_{2}\}.
\ea
Possible reduction of variables in the defining equations (in particular removing \(x\) and \(\lambda\) from the equations) in \eqref{Trans}
is desirable.
In order to achieve this goal, we let \(I\subset K[x, \lambda, \alpha]\) and \(J\subset \E_{x, \lambda, \alpha}\) be the ideals
generated by polynomial defining Equations \eqref{Trans} for either of \(\mathscr{B}\) and \(\mathscr{H}.\)
Here \(\E_{x, \lambda, \alpha}\) stands for all germs of smooth functions of \(x, \lambda,\) and \(\alpha.\)
For the case of \(\mathscr{D}\), a new variable
$\zeta$ is introduced and the quadratic germ \(1-\zeta(x_{1}-x_{2})\) is also added to the generators of the ideal \(I\) and \(J\).
This is due to the fact that \(x_1\neq x_2.\) Next, we compute the \Gr basis \(G\) for \(I\) in \(K[x, \lambda, \alpha], \alpha\in \mathbb{R}^k,\) with respect to \(\alpha_j\prec_\lex \lambda\prec_\lex x\) for \(j= 1, \ldots, k.\) Thus by \cite[Theorem 2]{CoxLittleIdeals}, \(G\cap K[\alpha]\) is a \Gr basis for
\(I\cap K[\alpha].\) Hence, \(G\cap K[\alpha]\subseteq J\cap \E_\alpha.\) Additional restrictions on the transition sets are also obtained by the other
elements of \Gr basis \(G.\) This justifies how we are able to convert the computation of elimination ideal from \(\E\) to its elimination ideal analogue in \(K[x, \lambda]\).

\item {\bf Elimination ideals}

Parametric smooth germs are converted to equivalent parametric polynomial germs
using the command {\tt UniversalUnfolding} (with the option \verb|normalform|; see \cite[Page 7]{GazorKazemiUser}). Then, using the procedure described above and the
{\tt EliminationIdeal} function in {\sc Maple}, we efficiently handle the possible reduction of variables for each of the Equations in \eqref{Trans}. The significance of this computational tool in {\tt Singularity} is not only for systematic derivation of transition sets but also for the cases that hand and numeric calculations fail.
\item {\bf Choosing points from connected components}

In order to classify persistent bifurcation diagrams, we need to pick one
point from each connected component in the complement of the transition set in the parameter space. To bring this in algorithmic fashion we use the function {\tt CylindricalAlgebraicDecompose} available in the {\tt RegularChains} library in {\tt Maple}. The drawback of this function for our purpose is that it generates more than one point in each connected component. We have written a program to reduce unnecessary and extra points, although further refinements is yet necessary.

\item {\bf Recognition problem}
\begin{itemize}

\item[(a)] Recognition problem for a normal form of a singular germ \(g\) determines the generic and degenerate conditions that a singular germ \(f\) is contact-equivalent to \(g\). The generic and degenerate conditions refer to a list of zero and non-zero conditions on derivatives of \(f\); see \cite[Page 88]{GolubitskySchaefferBook}. The command {\tt RecognitionProblem} in our library finds this list through the computation of \(\s\) and \(\s^{\perp}\). Next, intrinsic
generators of \(\s\) and monomials in \(\s^{\perp}\) derive non-zero and zero conditions, respectively.


\item[(b)] To find a solution for the recognition problem of a universal unfolding, we take a non-parametric germ, say \(g,\) as an input
and assume \(G\) as a parametric unfolding for \(g.\) Next normal form of \(g\) is computed and replaced with \(g\). Thus, the tangent space of \(g\) is computed via the command {\tt T}. Afterwards, the procedure described in \ref{VectorComplement}(a) is applied on the intrinsic part of the tangent space. This gives rise to a list of
monomials. Following \cite[Page 139]{GolubitskySchaefferBook} the projection map \(J\) from the
Taylor series of \(g\) to those monomials is derived. Then the list of non-zero elements \(J(x^{r}\lambda^{s}g)\), \(J(x^{r}\lambda^{s}g_{x})\), and \(J(\lambda^{s}g_{\lambda})\) are computed. Now we have the required information to apply \cite[Page 139, part (iii)]{GolubitskySchaefferBook} to obtain the desired matrix determinant associated with
\(G\) and \(g\). These are performed in {\tt Singularity} using the command {\tt RecognitionProblem} with the option of {\tt universalunfolding}. We remark that the
default function of {\tt RecognitionProblem} is for normal form that is described in part \((a)\).
A similar procedure enables us to verify if a parametric germ is a universal unfolding
for its folded singular germ.

\end{itemize}

\end{enumerate}

\section{Main features of \texttt{Singularity}}\label{SecFeatures}

The readers are referred to our user-guide \cite{GazorKazemiUser} for a list of all functions, their capabilities, options and comprehensive information on how to work with \texttt{Singularity}. In this section we merely describe the main features of \texttt{Singularity}. The main functions (not all) are given in Tables \ref{MainFuncSing} and \ref{MainFuncGeom}. \texttt{Singularity} has been tested by all scalar examples (nonsymmetric and without modal parameters) and classifications given in \cite{GovaertsBook,Melbourne87,Keyfitz,GolubitskySchaefferBook} and a few (differences, error or not already reported data due to computational burden) are verified in our favor. To illustrate how these functions work, let
\bes h_1(x, \lambda):= x^{5}+\lambda x+\lambda^{2}.\ees Then, \(\verb|AlgObjects|(h_1, [x, \lambda])\) returns
\bes
 \p:=\M^{6}+\M^{2}\langle\lambda\rangle+\langle\lambda^{2}\rangle, RT:=\M^{5}+\M\langle\lambda\rangle, T:=\M^5+\M\langle\lambda\rangle+K\lbrace x+2\lambda, x^4+\frac{1}{5}\lambda\rbrace,
 \ees
 \bes
\E/T:=K\{1, \lambda, x^2, x^3\}, \s:=\M^{5}+\M\langle \lambda\rangle, \s^{\perp}:=K\lbrace 1, \lambda, x, x^2, x^3, x^4\rbrace,
\ees
and intrinsic generators of \(\s\) are \(x^5\) and \(x\lambda.\) For an example of computations in infinite codimensional ideals, consider restricted tangent space of \(h_2(x, \lambda)= \lambda^3\sin(x).\) \({\tt RT(h_2, [x,\lambda], InfCodim)}\) generates the restricted tangent space of \(h_2\) as \(\M\langle \lambda^3\rangle.\)

\begin{table}[h]
\caption{\texttt{Singularity}'s main functions in singularity theory.}\label{MainFuncSing}
\begin{center}\footnotesize
\renewcommand{\arraystretch}{1.3}
\begin{tabular}{|ll|ll|}\hline
\multicolumn{2}{|c|}{{\bf Function}} & \multicolumn{2}{c|}{{\bf Description}}\\ \hline
\verb"Verify"&& information on truncation degree and computational ring. &
\\ \hline
\verb|AlgObjects|& & \(\p\), \(RT\), \(T\), \(\E/T\), \(\s\),\(\s^{\perp},\) intrinsic generators of \(\s\).&\\ \hline
\verb|Normalform|& & normal forms of a given germ. &\\ \hline
\verb|UniversalUnfolding|& & universal unfoldings of a given germ. &\\ \hline
\verb|RecognitionProblem|& & for normal forms and universal unfoldings.
 &\\ \hline
\verb|TransitionSet|& & transition set are computed and plotted or animated. &\\ \hline
\verb|PersistentDiagram|& & plots or animates all persistent bifurcation diagrams. &\\ \hline
\verb|Transformation|& & estimates \(S\) and \(X\) relating two contact-equivalent germs. &\\ \hline
\verb|Intrinsic|& & intrinsic part of a given ideal or vector subspace of \(\E.\) &\\ \hline
\end{tabular}
\end{center}
\end{table}

\verb"Verify"(\(g\), \verb"Vars") checks a germ \(g\) for its bifurcation analysis while \verb"Verify"(\(G\), \verb"Vars", \verb"Ideal") checks germ generators \(G\) of an ideal \(\langle G\rangle_\E\) for divisions or standard basis computations. In either case, it suggests suitable computational rings and a truncation degree and it verifies that their use does not lead to error according to Theorems \ref{Thm2.18} and \ref{Thm2.18a}. Here \verb"Vars", say \([x, \lambda],\) stand for the state variable \(x\) and the distinguished parameter \(\lambda\). For instance, \(\verb|Verify|(\sin(x) \tanh(x^3-\lambda), [x, \lambda])\) returns computational rings as smooth germs, formal power series, fractional maps and truncation degree 5. Now consider the ideal \(I\) given in Equation \eqref{specialexm}. Then,
\(\verb|Verify|(I, [x, \lambda], \verb|Ideal|)\) returns computational rings as smooth germs, formal power series, fractional maps and truncation degree 6. In other words, the polynomial germ ring is not allowed for this  example. The command \(\verb|Intrinsic|(I, [x, \lambda])\) gives rise to \(\M^6+\M^2 \langle \lambda^2\rangle\).

\(\verb|Normalform|(h_1, [x, \lambda])\) computes the normal form \(x^5+\lambda x\) while its universal unfolding
\be\label{G1sec} H_1:=x^5+\lambda x+\alpha_{1}+\alpha_2 \lambda+\alpha_3 x^2+\alpha_4 x^3\ee is derived by \({\tt UniversalUnfolding}(h_1, [x, \lambda], \verb"normalform").\)

The command \verb"RecognitionProblem"(\(x^3+\lambda^2 \cos(x), [x, \lambda]\), 6, \verb"Formal") answers the recognition problem by
\begin{center}
 "nonzero condition=", \([\frac{\partial^2}{\partial\lambda^2}f\neq 0, \frac{\partial^3}{\partial x^3}f\neq 0]\)\\
"zero condition=", \([f=0, \frac{\partial}{\partial \lambda}f=0 , \frac{\partial}{\partial x}f=0,\frac{\partial^2}{\partial x^2}f =0, \frac{\partial^2}{\partial x \partial \lambda}f =0]\).
\end{center} while with the universal unfolding option \verb"RecognitionProblem"(\(x^3+\lambda, [x, \lambda]\), \verb"UniversalUnfolding", \verb"Fractional") gives
\bes
\det\left( {\begin{array}{cc}
g_{\lambda}(0) & g_{\lambda, x}(0)\\
G_{\alpha_1}(0) & G_{\alpha_1, x}(0)
\end{array}}
\right)\neq 0.
\ees\\

The command \(\verb|TransitionSet|(H_1, [\alpha_1, \alpha_2, \alpha_3, \alpha_4], [x, \lambda])\) derives the transition set associated with \(H_1\) as \(\B:=\lbrace (\alpha_1, \alpha_2, \alpha_3, \alpha_4)\,|\, \alpha_{1}=\alpha_{2}^{5}+\alpha_{2}^{3}\alpha_{4}-\alpha_{2}^{2}\alpha_{3}\rbrace,\) \(\mathscr{H}\) and  \(\mathscr{D}\) given in Appendix \ref{HD}.


\begin{figure}[h]
\begin{center}
\subfigure[]{\includegraphics[width=.13\columnwidth,height=.13\columnwidth]{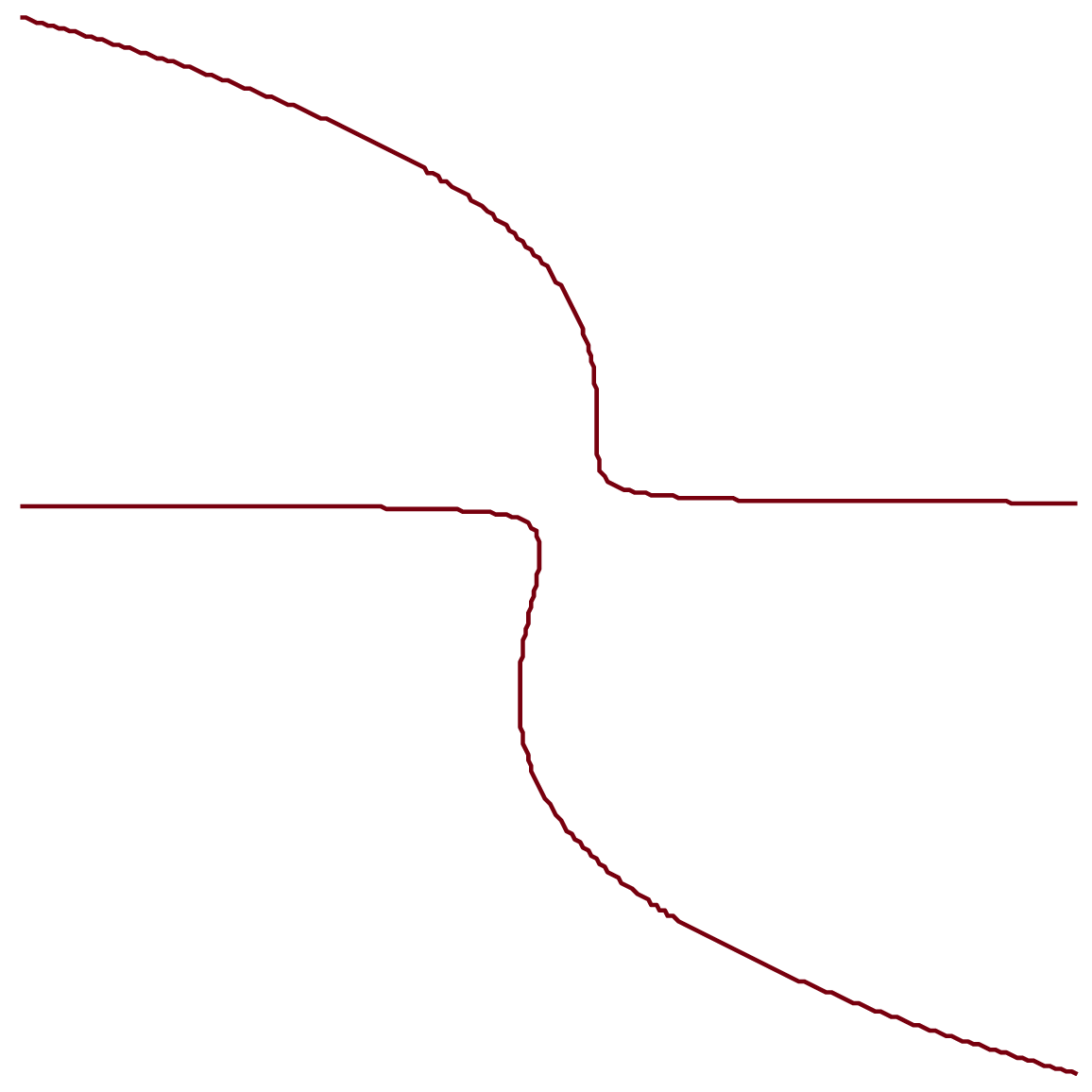}}
\subfigure[]{\includegraphics[width=.13\columnwidth,height=.13\columnwidth]{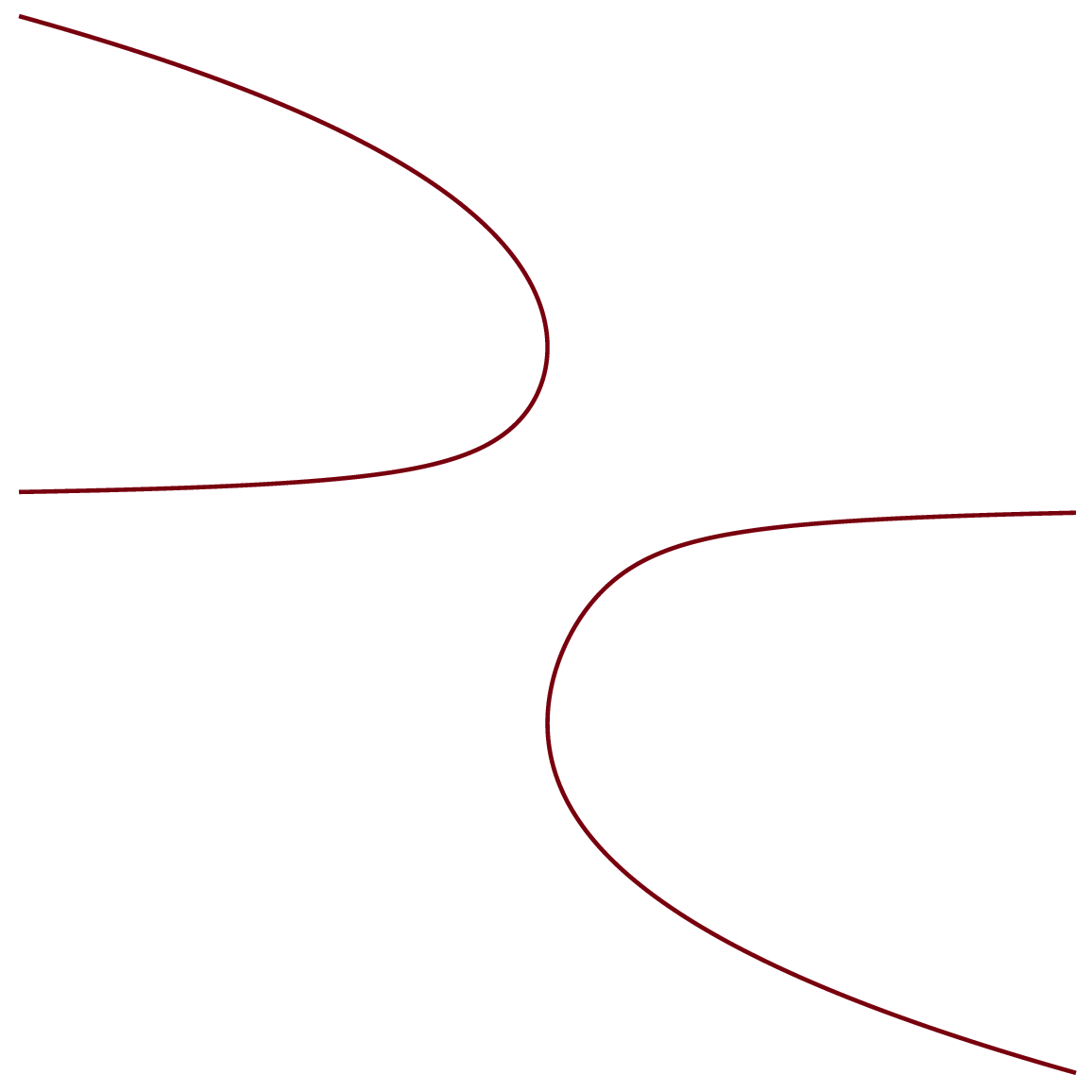}}
\subfigure[]{\includegraphics[width=.13\columnwidth,height=.13\columnwidth]{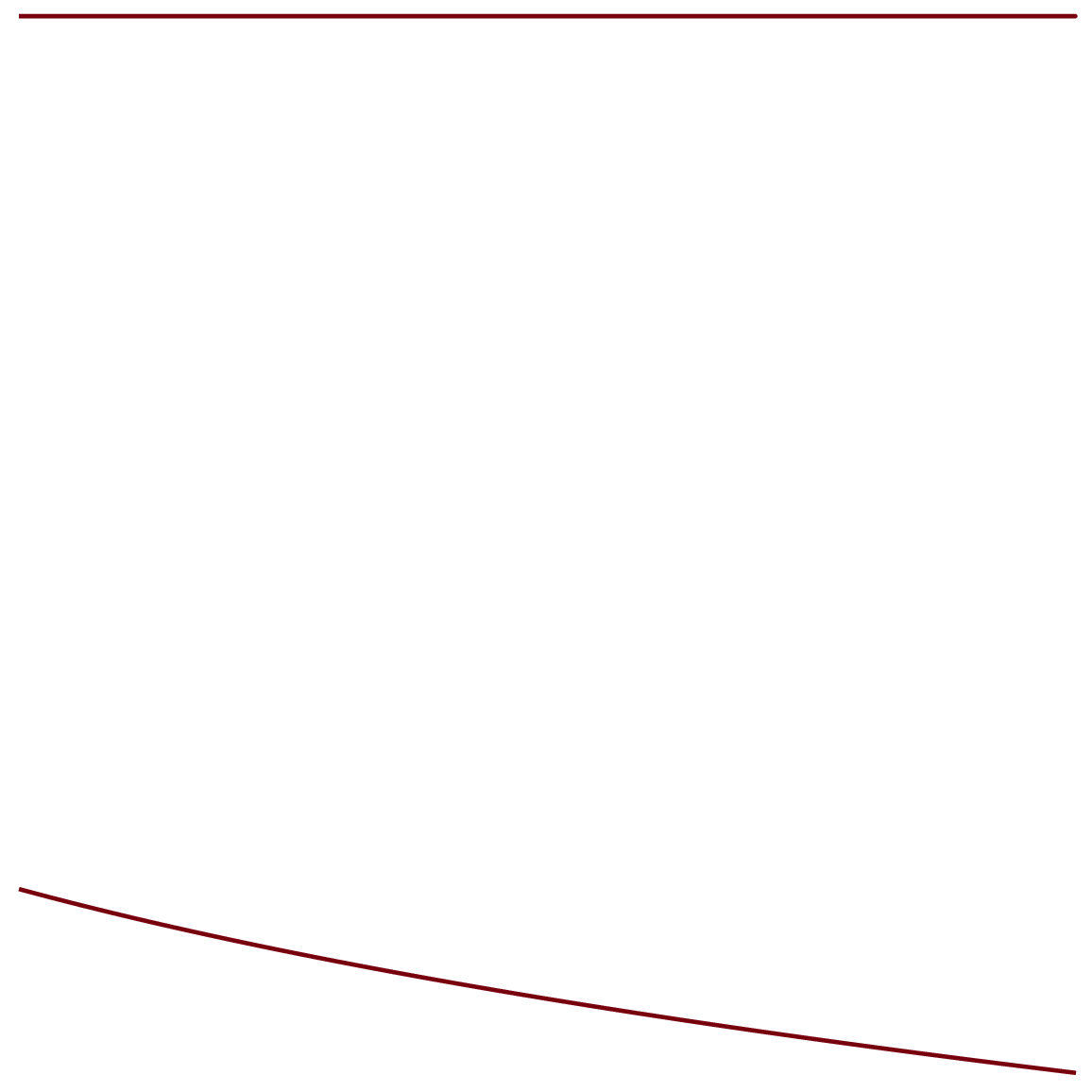}}
\subfigure[]{\includegraphics[width=.13\columnwidth,height=.13\columnwidth]{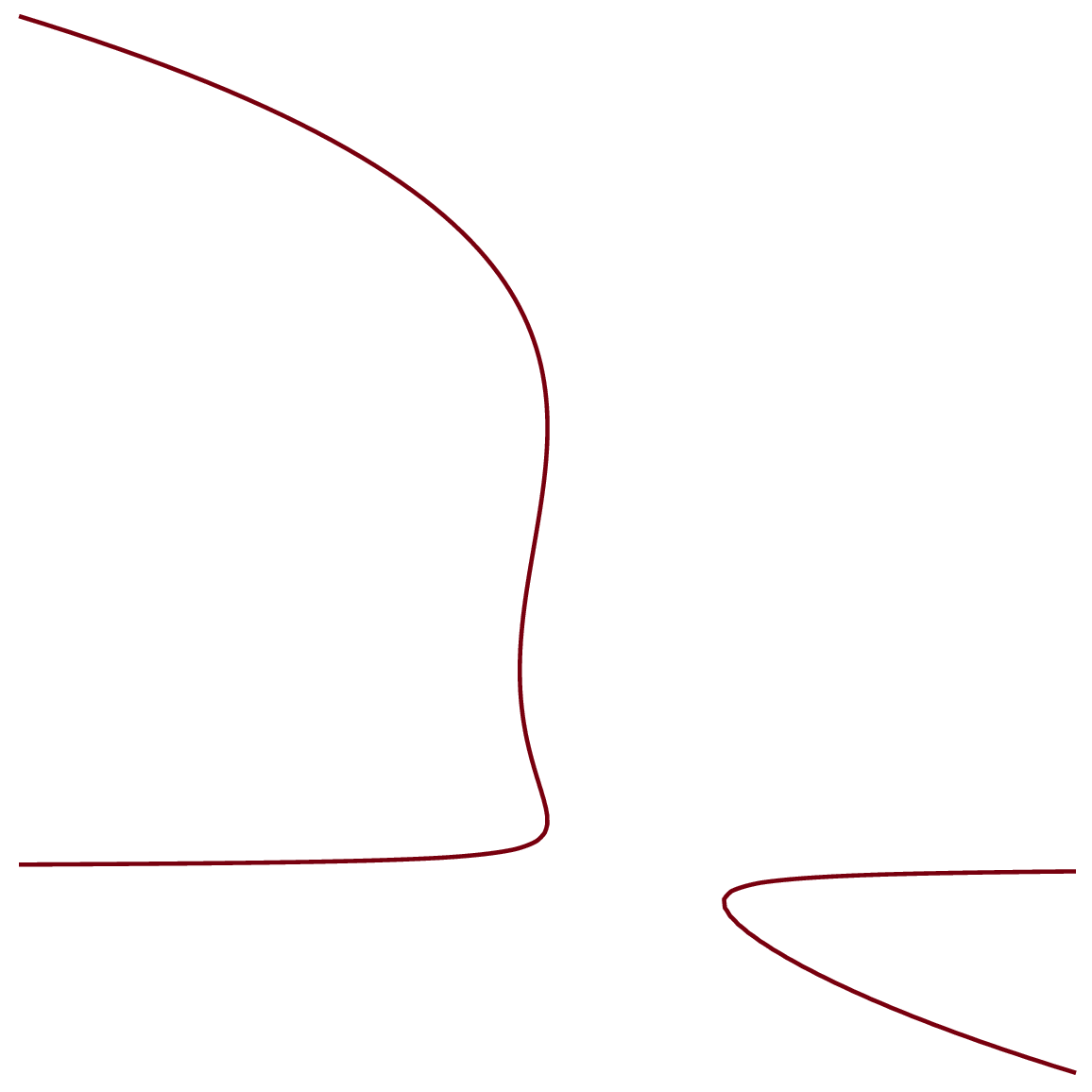}}
\subfigure[]{\includegraphics[width=.13\columnwidth,height=.13\columnwidth]{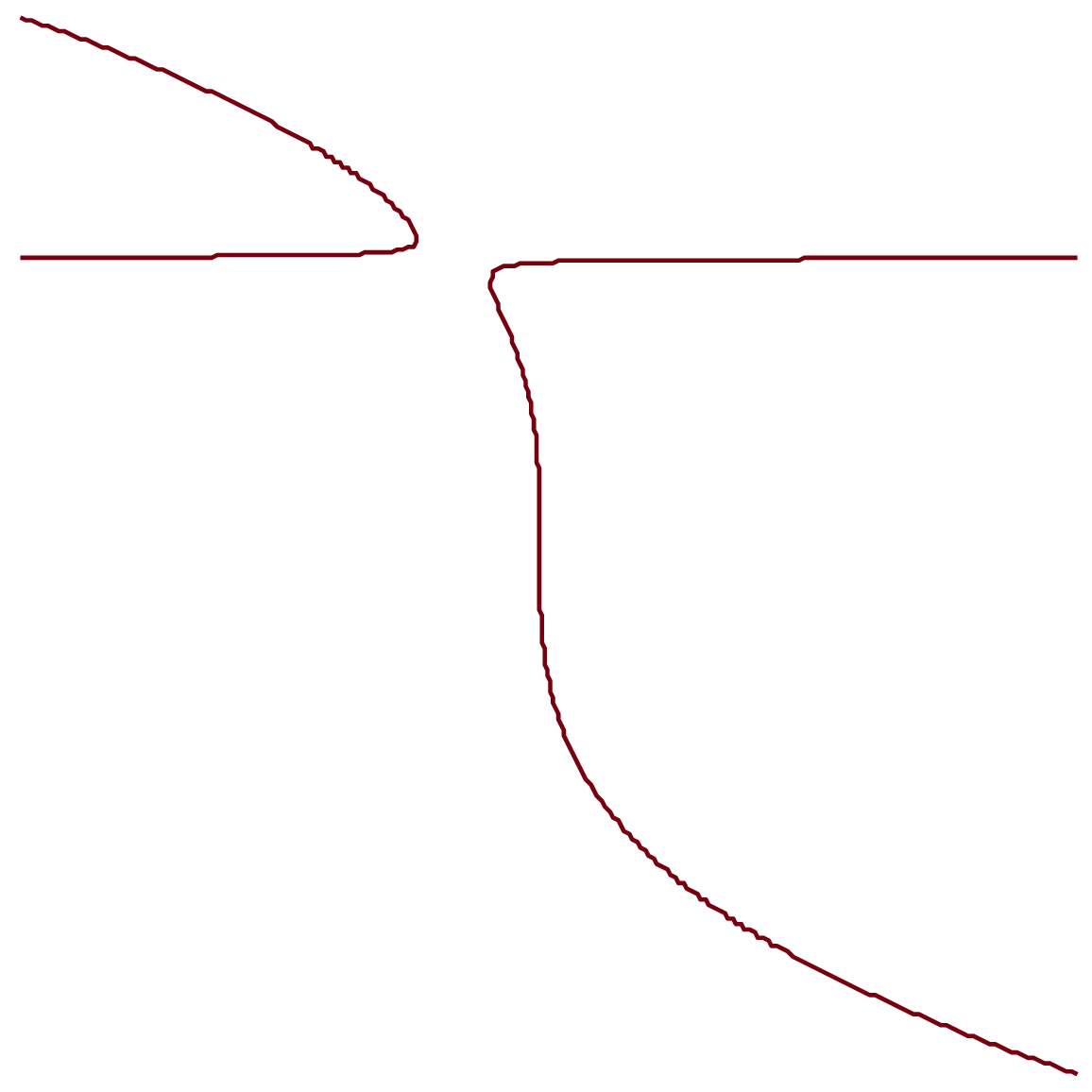}}
\subfigure[]{\includegraphics[width=.13\columnwidth,height=.13\columnwidth]{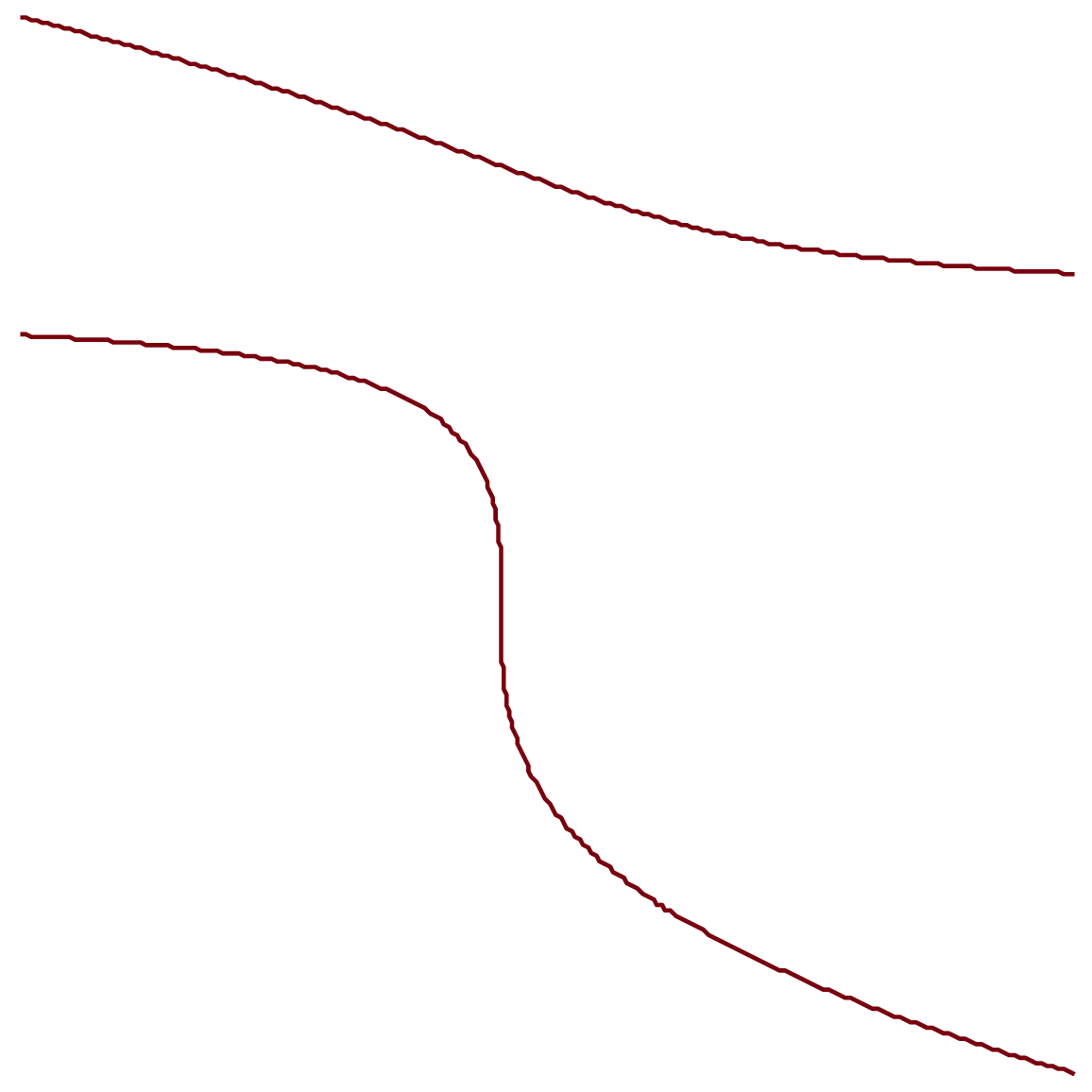}}
\subfigure[]{\includegraphics[width=.13\columnwidth,height=.13\columnwidth]{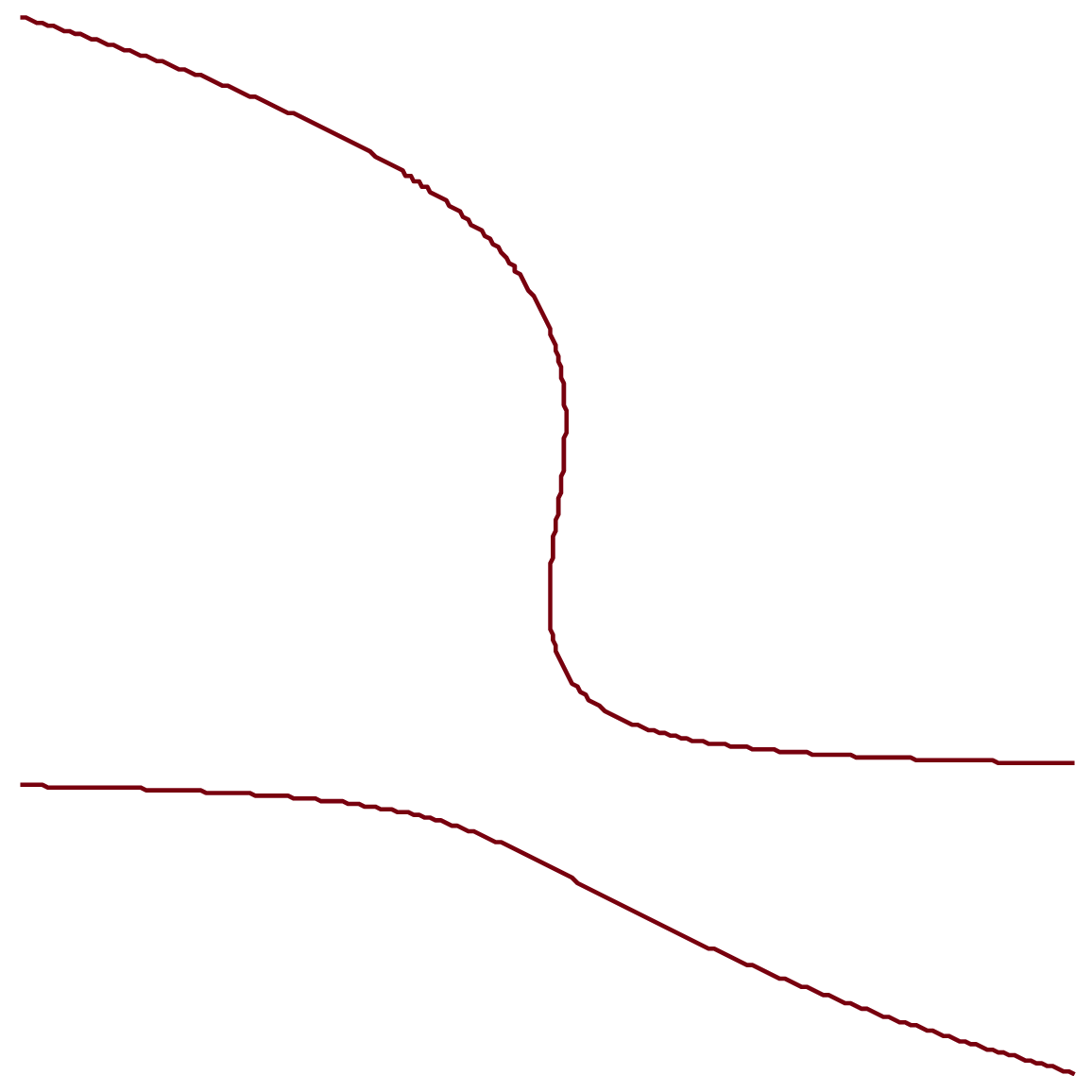}}
\subfigure[]{\includegraphics[width=.13\columnwidth,height=.13\columnwidth]{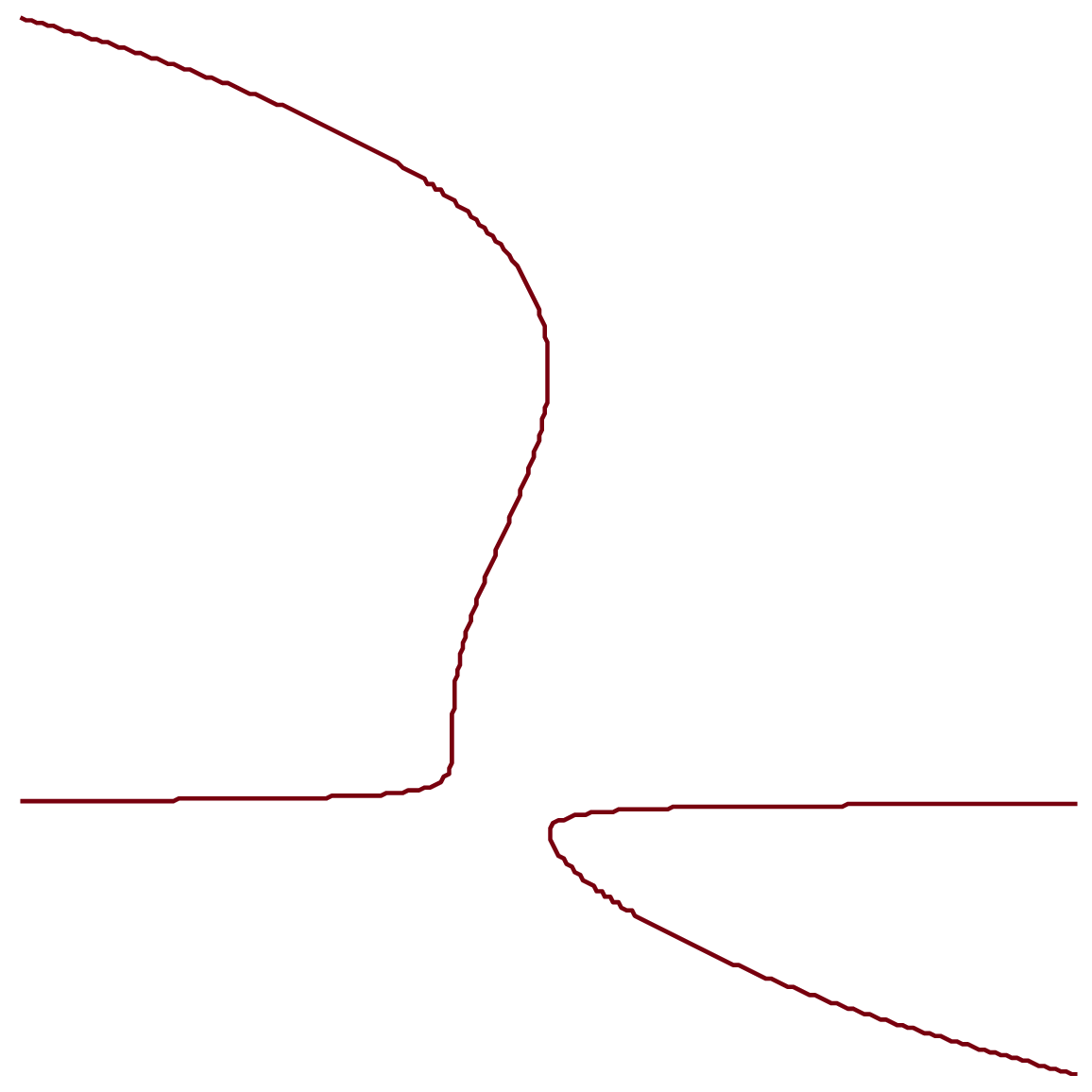}}
\subfigure[]{\includegraphics[width=.13\columnwidth,height=.13\columnwidth]{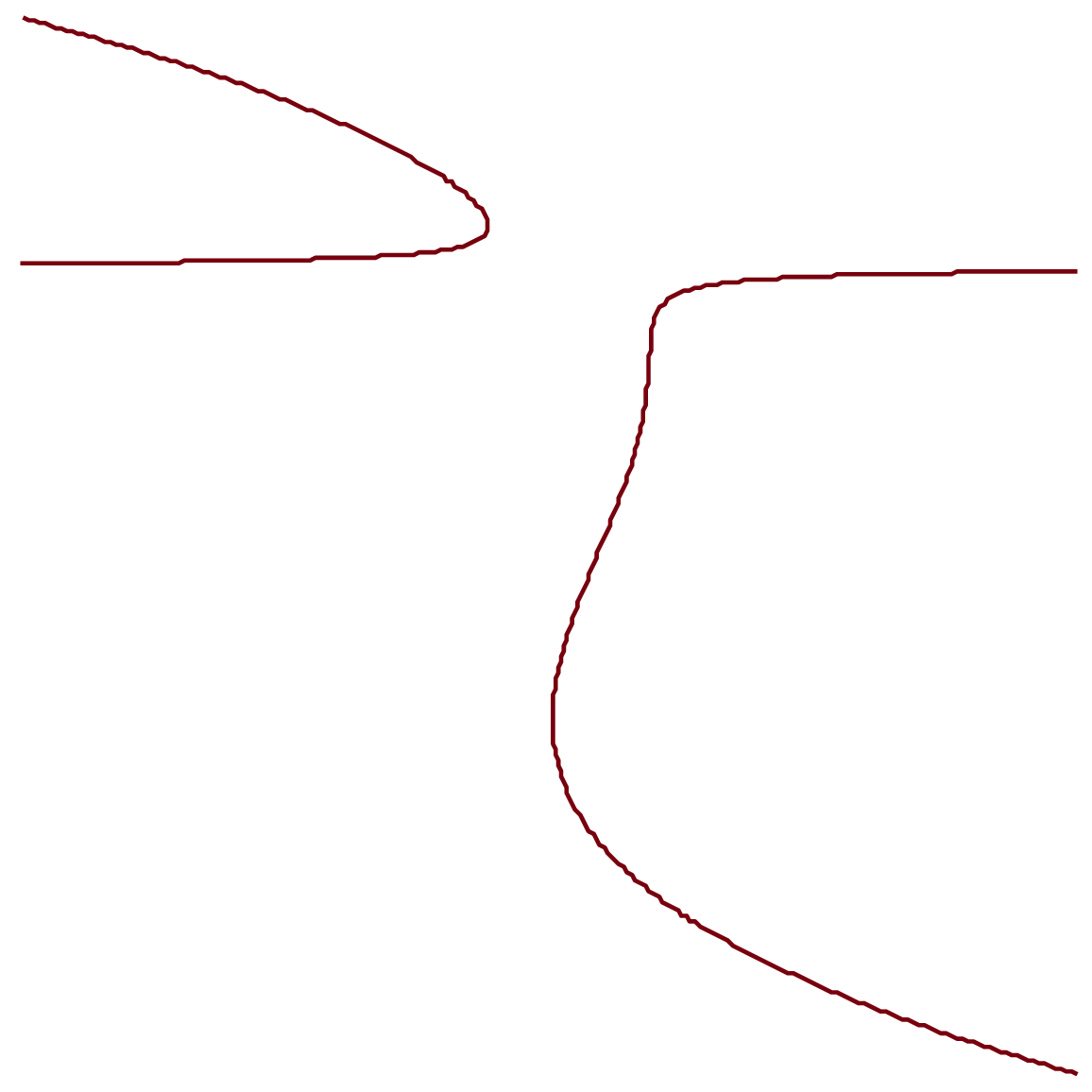}}
\subfigure[]{\includegraphics[width=.13\columnwidth,height=.13\columnwidth]{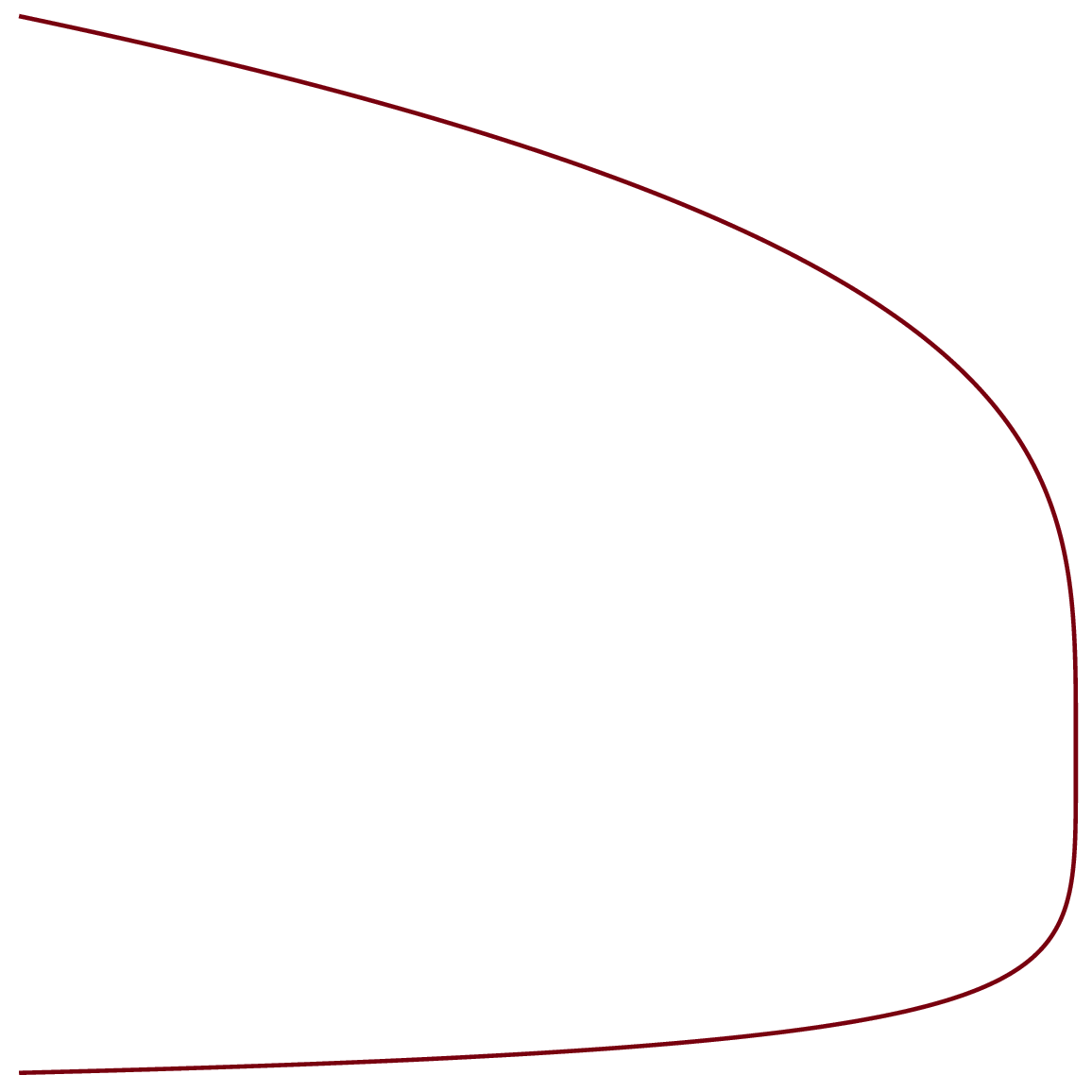}}
\subfigure[]{\includegraphics[width=.13\columnwidth,height=.13\columnwidth]{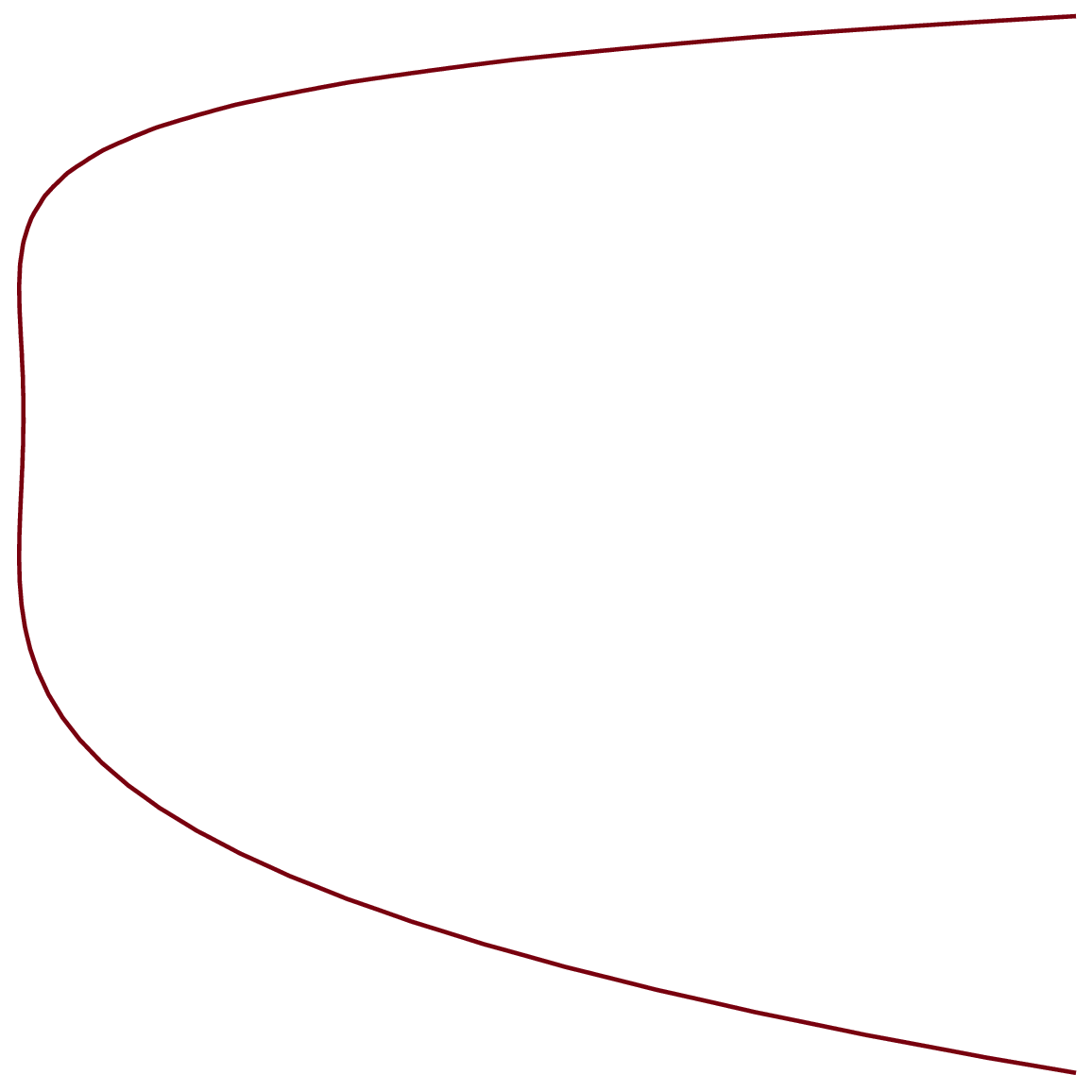}}
\subfigure[]{\includegraphics[width=.13\columnwidth,height=.13\columnwidth]{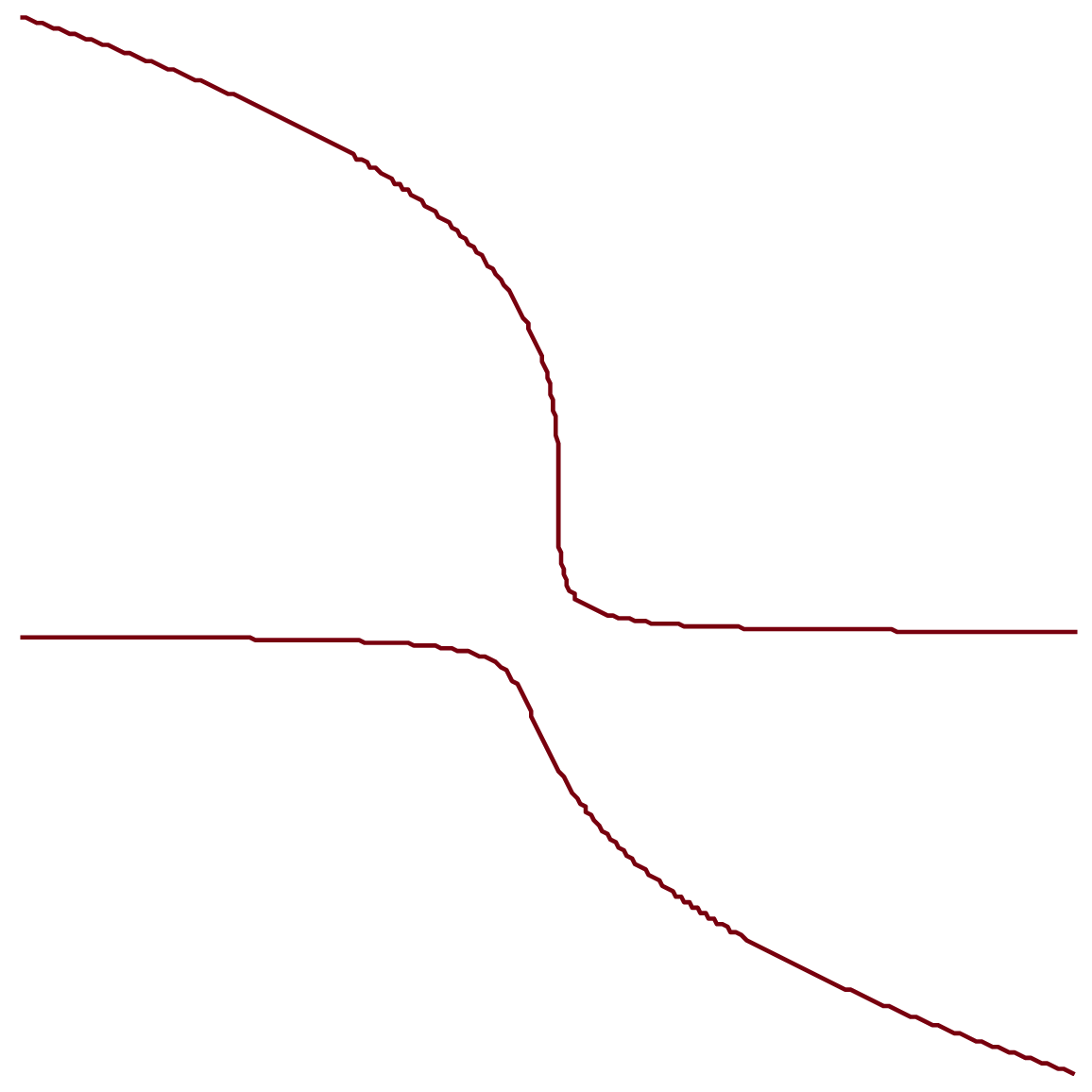}}
\subfigure[]{\includegraphics[width=.13\columnwidth,height=.13\columnwidth]{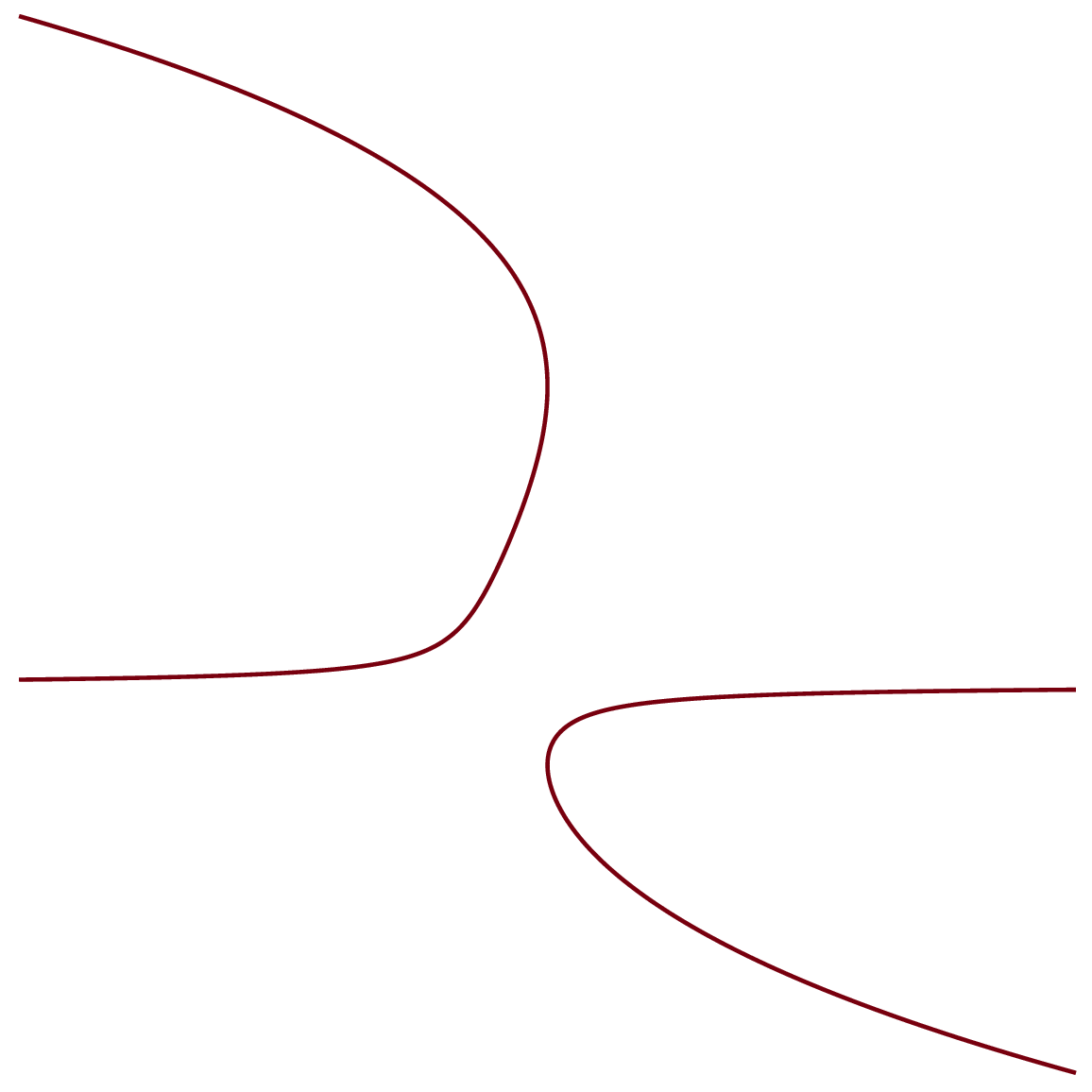}}
\subfigure[]{\includegraphics[width=.13\columnwidth,height=.13\columnwidth]{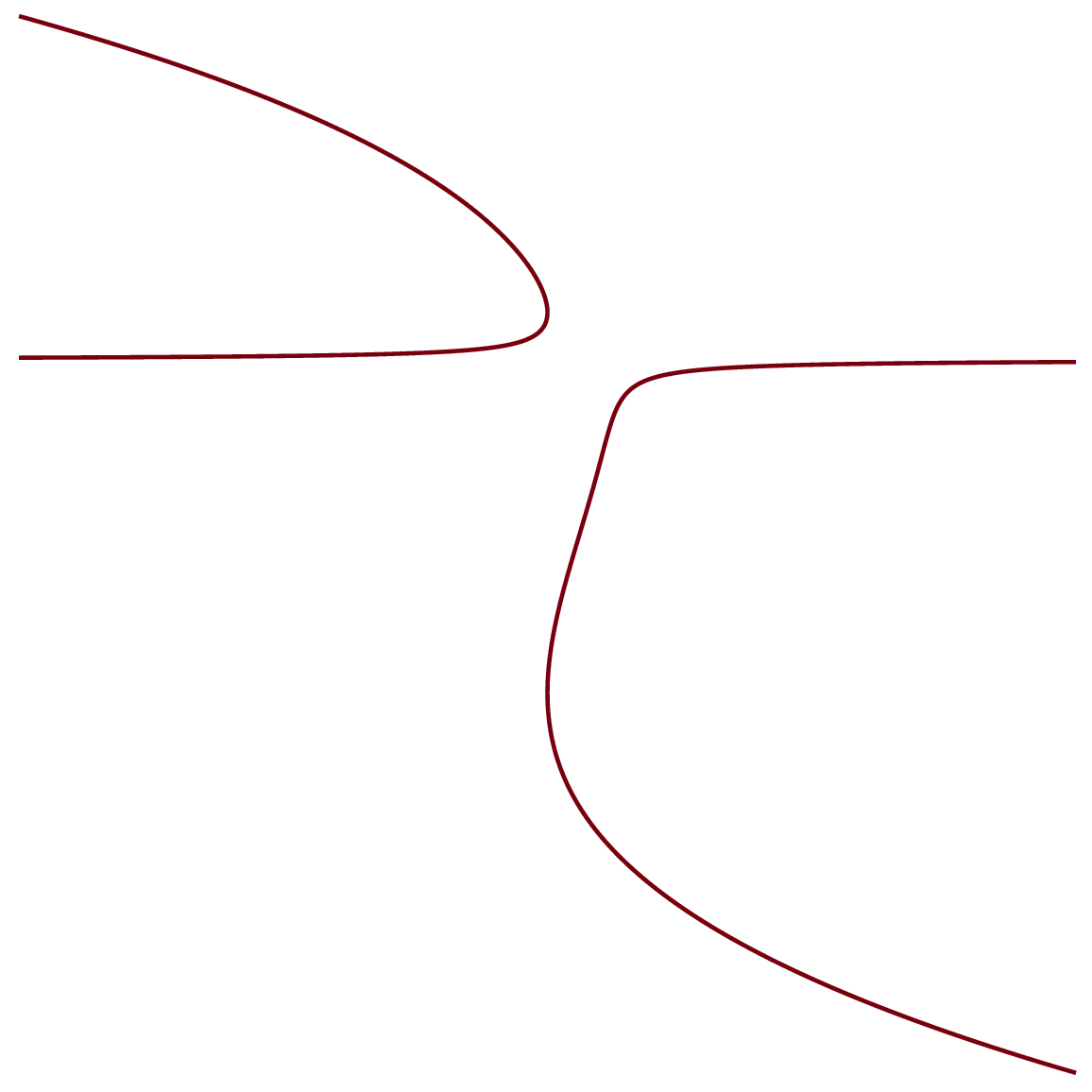}}
\caption{\small Persistent bifurcation diagrams for parametric singularity \(H_2\) in \eqref{G2}.}\label{Fig3}
\end{center}
\end{figure}

Now we consider the parametric germ
\be\label{G2} H_{2}=x^4+\lambda x+\alpha_{1}+\alpha_{2}\lambda+\alpha_{3}x^{2},\ee that is a universal unfolding for germ \(x^4+\lambda x\). Next, \(\verb|TransitionSet|(H_2, [\alpha_1, \alpha_2, \alpha_3], [x, \lambda])\) gives rise to
 \(\B:=\lbrace (\alpha_1, \alpha_2, \alpha_3)\,|\,\alpha_{2}^{4}+\alpha_{2}^{2}\alpha_{3}+\alpha_{1}=0 \rbrace\), \(\mathscr{H}:=\lbrace
(\alpha_{1}, \alpha_{2}, \alpha_{3})\,|\,128\alpha_{2}^{2}\alpha_{3}^{3}+3\alpha_{3}^{4}
+72\alpha_{1}\alpha_{3}^{2}+432\alpha_{1}^{2}=0\rbrace\), and \(\mathscr{D}:=\lbrace
(\alpha_{1}, \alpha_{2}, \alpha_{3})\,|\,\alpha_{3}^{2}-4\alpha_{1}=0, \alpha_{3}\leq 0\rbrace.\) A list of persistent bifurcation diagrams is generated by the command \(\verb|PersistentDiagram|(H_2, [x, \lambda], \verb|plot|, \verb|IntermediateList|)\); Figure \ref{Fig3} is some inequivalent diagrams chosen from this list.

\begin{table}[h]
\caption{Main functions of \texttt{Singularity} associated with the local rings \(K[[x, \lambda]], \RM\) and \(\E\).} \label{MainFuncGeom}
\begin{center}\footnotesize
\renewcommand{\arraystretch}{1.3}
\begin{tabular}{|ll|ll|}\hline
\multicolumn{2}{|c|}{{\bf Function}} & \multicolumn{2}{c|}{{\bf Description}}\\ \hline
\verb|StandardBasis|& & standard basis in either of the rings.&\\ \hline
\verb|ColonIdeal|& & computes the colon ideal given in Equation \eqref{Colon}. &\\ \hline
\verb|Division|& & remainder of a germ \(g\) divided by a set \(G\).
 &\\ \hline
\verb|MultMatrix|& & computes the matrix \(\varphi_{u, J}\) in Equations \eqref{R/J}. &\\ \hline
\verb|Normalset|& & finds a basis for complement space of an ideal \(I.\) &\\ \hline
\end{tabular}
\end{center}
\end{table}

In order to develop \texttt{Singularity}, we have implemented some local tools (suitable in our context) from computational algebraic geometry including {\it division remainders,} {\it standard bases,} {\it elimination ideals,} {\it ideal membership problem} and {\it colon ideals} for the local rings of fractional (germ) maps, formal power series and smooth maps. These are accordingly implemented in \texttt{Singularity}. We hereby thank Amir Hashemi for his frequent fruitful discussions. We also acknowledge Benyamin M.-Alizadeh's helps and ideas in the early stages of this project. They were essentially helpful in a fast learning of the concepts from algebraic geometry and their programming. Let
\bes G:=\{g_1\!:=\sin(\lambda^{7}+x)+\exp(x^{4})-x-1-\lambda^{9}, g_2\!:=x^{5}-\lambda^{2}, g_3\!:=\cos(x^{6})-\lambda-1\}\ees
and \(k\) be a sufficiently large truncation degree. Then using the commands
\bes \verb|StandardBasis|(G, [x, \lambda], k, \verb|Fractional|), \verb|StandardBasis|(G, [x, \lambda], k, \verb|Formal|),\ees and \(\verb|StandardBasis|(G, [x, \lambda], k, \verb|SmoothGerms|),\) we obtain the Standard basis of \(G\) in the rings of fractional maps, formal power series and smooth germs \(\E.\) In this example, the standard basis in either of these rings is given by \(\lbrace \lambda, x^{3}\rbrace.\) Note that the options \(\verb|FormalSeries|\) and \(\verb|Fractional|\) in \(\verb|StandardBasis|\) use a truncated Taylor series expansion of germs (for non-polynomial and non-fractional germs) and are adapted according to Theorem \ref{Thm2.18a}. We have already developed some {\sc Maple} programs for (parametric and orbital) normal form computation, bifurcation analysis and control of {\it singular differential equations} \cite{GazorSadri2017,GazorSadri,GazorMoazeni,GazorMokhtari2015,GazorMokhtariInt,GazorMokhtariSanders,GazorYuSpec} and their integration with \texttt{Singularity} shall lead to a toolbox for local bifurcation control and analysis of singularities.

\section{An illustrating buckling example}\label{Sec8}

\begin{figure}[h]
\begin{center}
\subfigure[\label{Model1}]{\includegraphics[width=.25\columnwidth,height=.2\columnwidth]{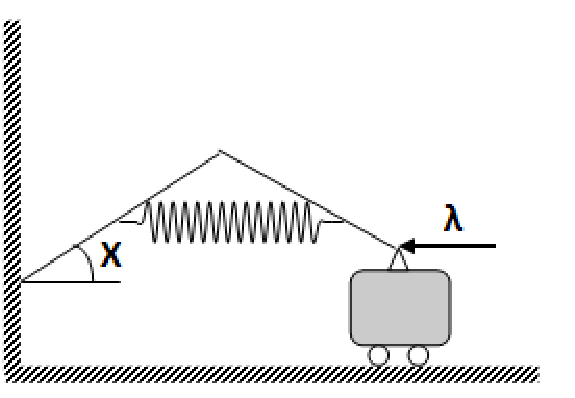}}
\subfigure[\label{Model23}]{\includegraphics[width=.25\columnwidth,height=.2\columnwidth]{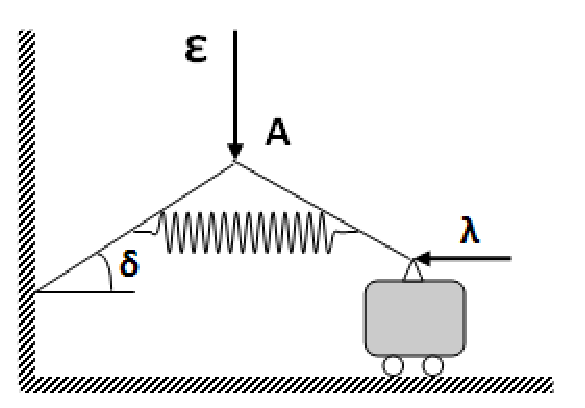}}
\subfigure[\label{pich_fork}]{\includegraphics[width=.2\columnwidth,height=.2\columnwidth]{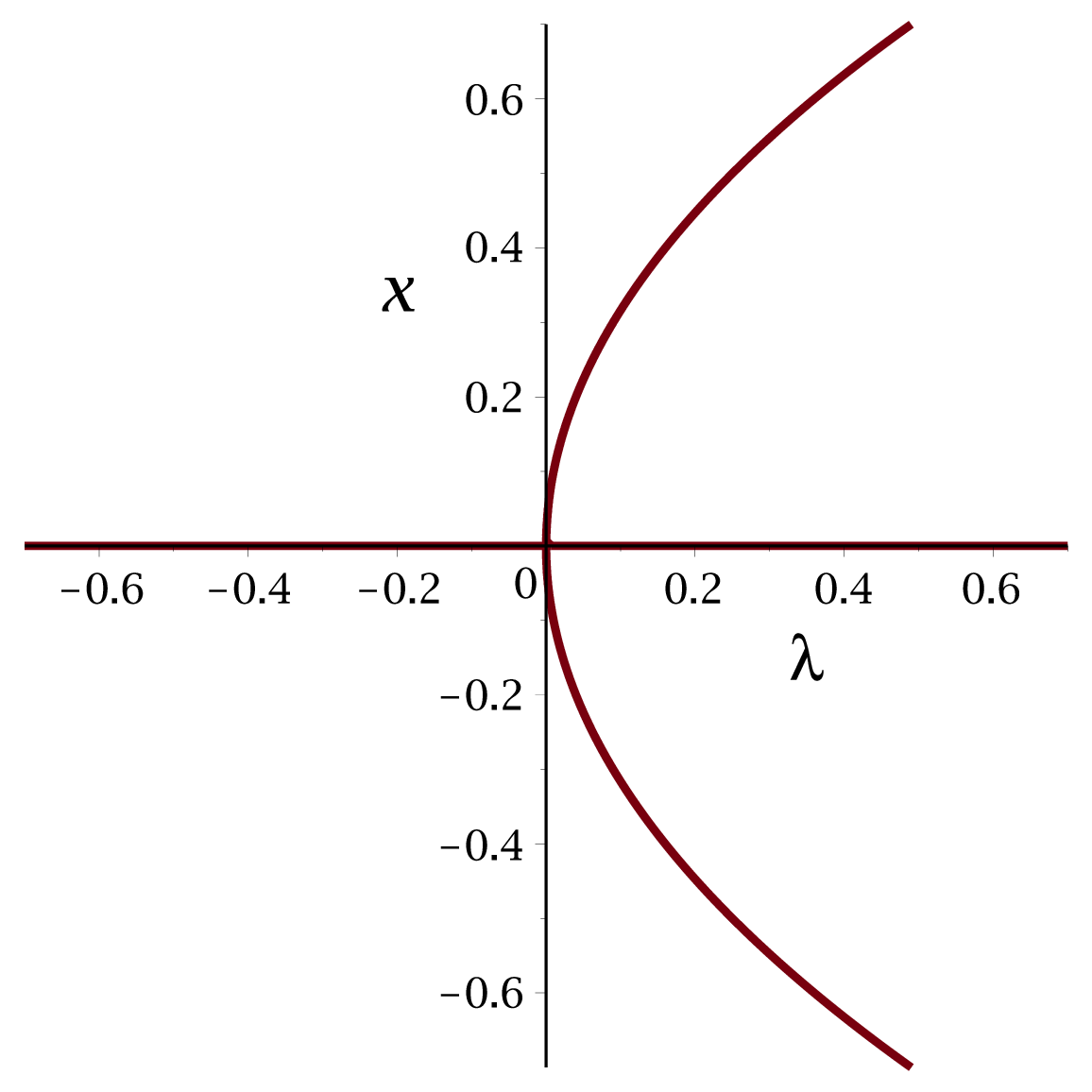}}
\subfigure[\label{transition_set_1}]{\includegraphics[width=.2\columnwidth,height=.2\columnwidth]{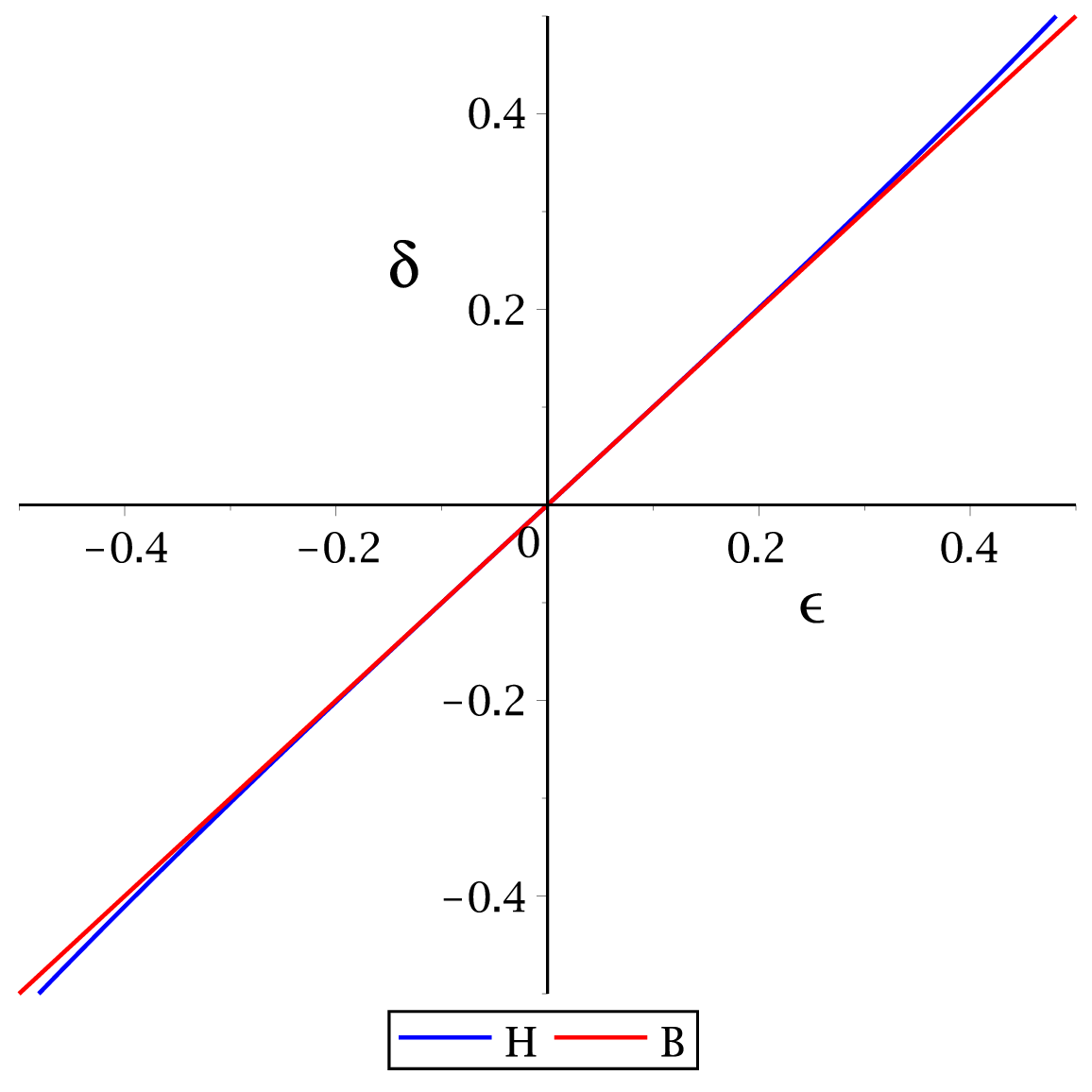}}
\subfigure[\label{first1}]{\includegraphics[width=.2\columnwidth,height=.2\columnwidth]{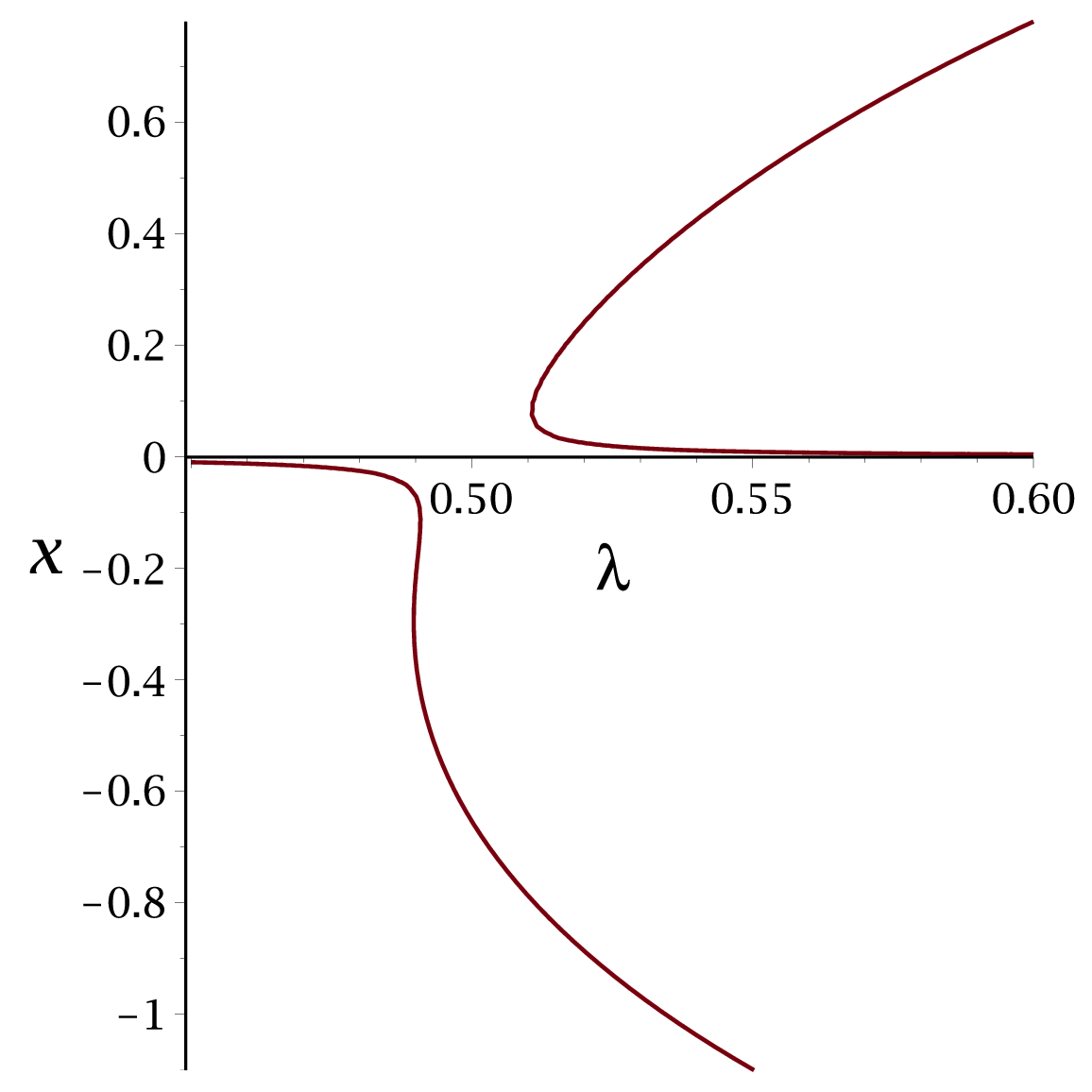}}
\subfigure[\label{first2}]{\includegraphics[width=.2\columnwidth,height=.2\columnwidth]{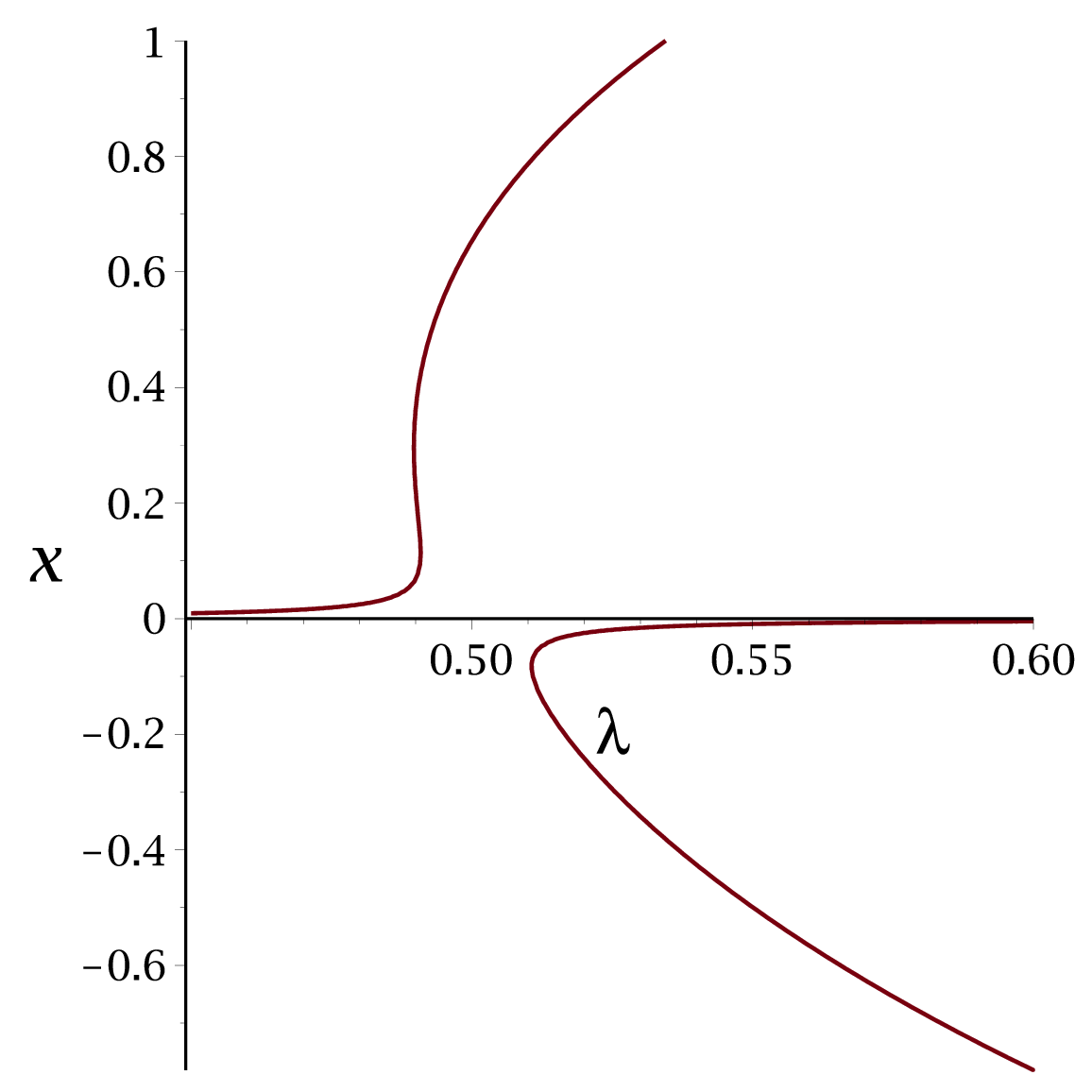}}
\subfigure[\label{first3}]{\includegraphics[width=.2\columnwidth,height=.2\columnwidth]{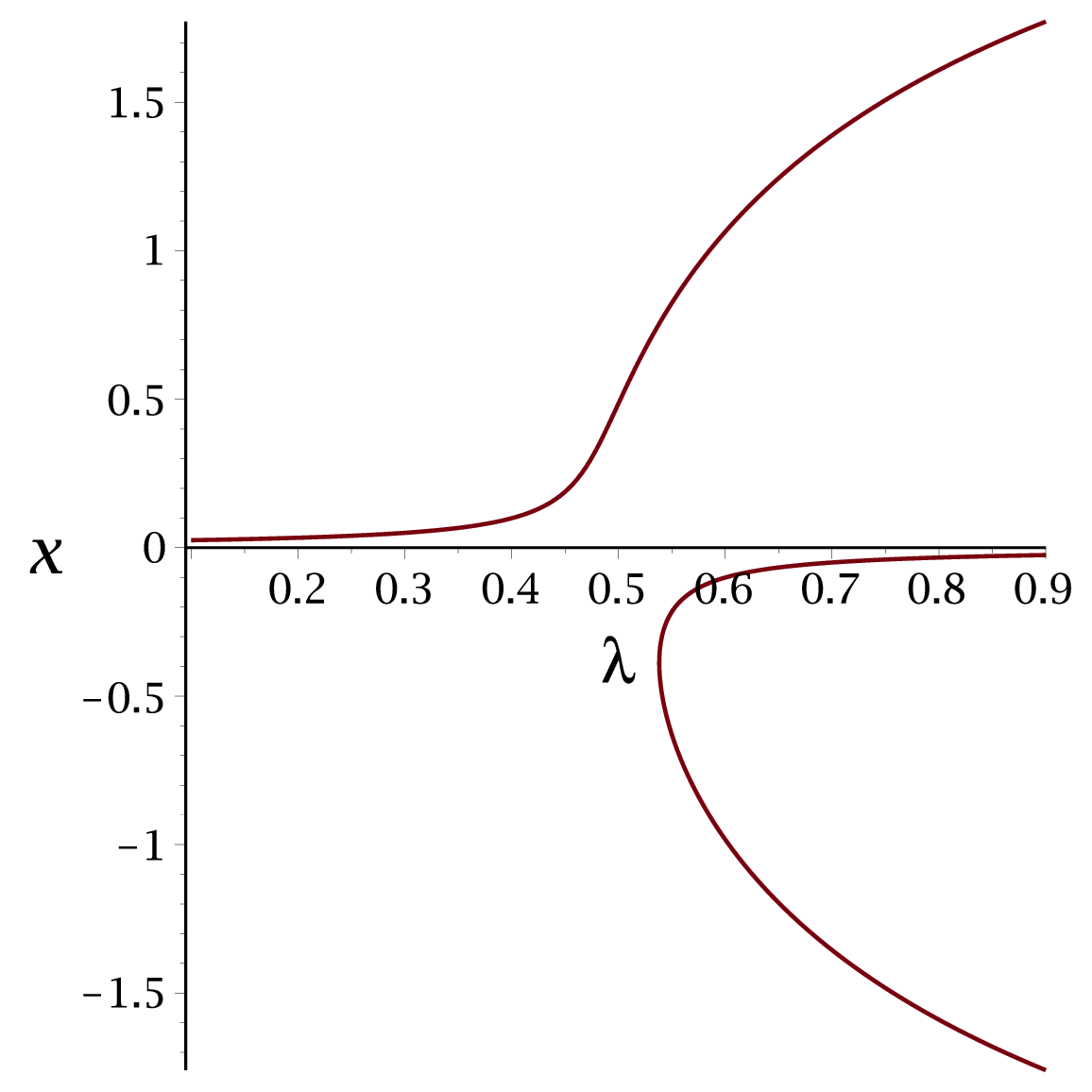}}
\subfigure[\label{first4}]{\includegraphics[width=.2\columnwidth,height=.2\columnwidth]{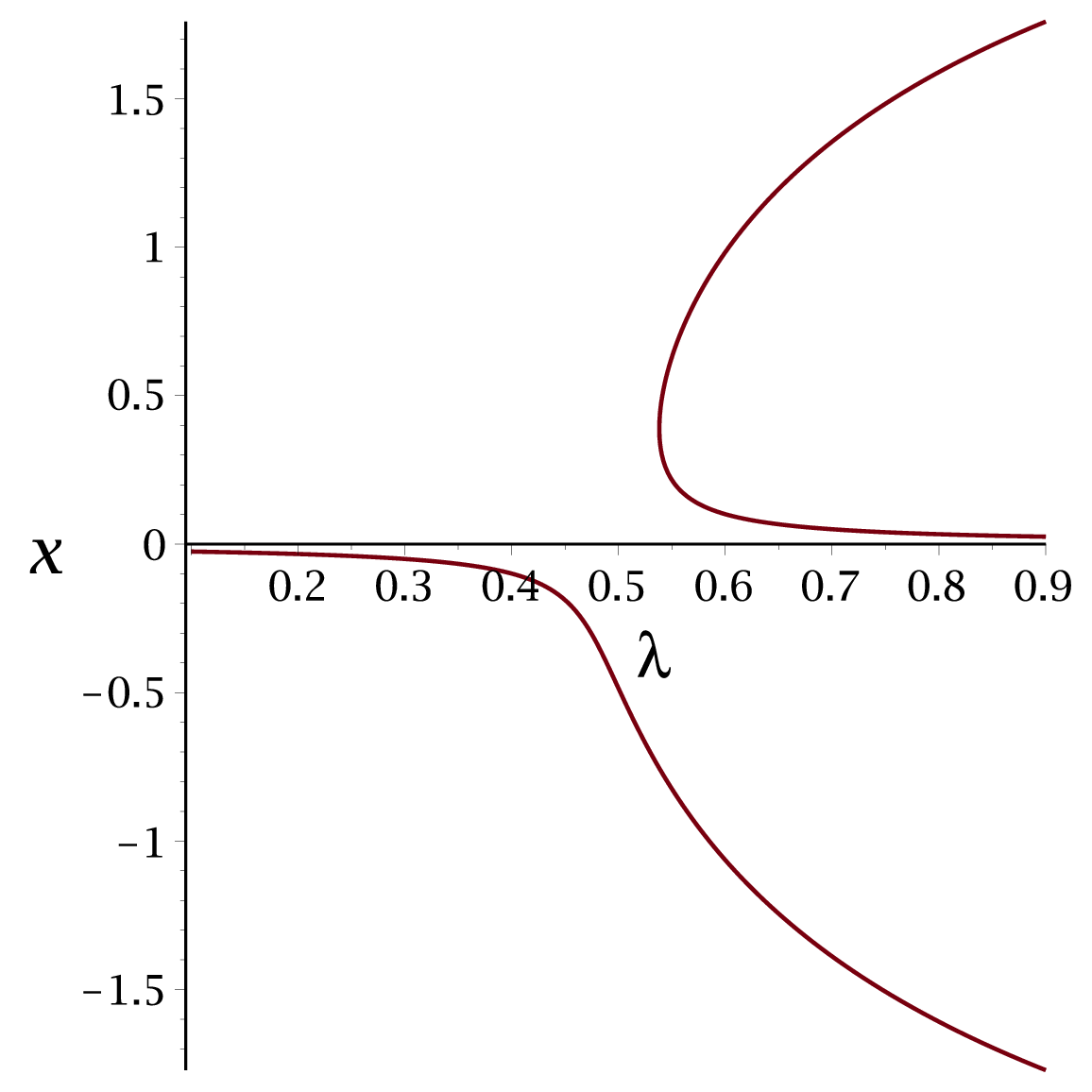}}
\caption{The buckling models, pitchfork bifurcation diagram, transition set and numerical persistent bifurcation diagrams}\label{transition_first_model}
\end{center}
\end{figure}

In order to illustrate how {\tt Singularity} can be used in an application, we borrow a {\it finite element analogue of Euler buckling} example from \cite[pages 3-10, figures 1.2 and 1.4]{GolubitskySchaefferBook}, \ie figures \ref{Model1} and \ref{Model23}. The associated
steady-state solutions (figure \ref{Model1}) follow
\be
g(x, \lambda)= \frac{\partial V}{\partial x}(x, \lambda)=x-2\lambda \sin x.
\ee The command {\tt CheckSingularity}\((g, [x, \lambda])\) provides the singular point \((x, \lambda)=(0, \frac{1}{2})\), where \(g\) is singular.
\verb"Verify" derives the permissible truncation degree and computational rings, \eg

\verb"Verify"(\(g, [x, \lambda]\), {\tt 'SingularPoint'}='\([0, \frac{1}{2}]\)') leads to
\begin{itemize}
  \item The following rings are allowed as the means of computations: Ring of smooth germs, Ring of formal power series, Ring of fractional germs, Ring of polynomial germs.
\item The truncation degree must be: 4.
\item Recommended rings are: Fractional ring, Polynomial ring.
\end{itemize}

\noindent Now the normal form type is derived by {\tt Normalform}\((g, [x, \lambda], 4,\) '{\tt SingularPoint}'='\([0, \frac{1}{2}]\)') as
\(x^3-\lambda x.\) The associated bifurcation diagram can be made available using {\tt PersistentDiagram}\((x^3-\lambda x, [x, \lambda]),\) see figure \ref{pich_fork}.

Since \(g\) is not a universal unfolding, any possible small perturbations of \(g\) (such as modeling imperfections) has the potential to dominate the associated dynamics of the problem. Thus, the actual dynamics of the buckling problem does not follow the pitchfork bifurcation diagram \ref{pich_fork}. To treat this caveat, an approach is to determine the universal unfolding of \(g\) by the command \verb"UniversalUnfolding" and accordingly design a universal unfolding model for the buckling experiment via imposing additional parameters into the problem; see \cite{GolubitskySchaefferBook,GolubitskySchaefferImperfection}. One way to feed the additional unfolding parameters into the problem is to consider them as the controller inputs for a controlled experimental problem. Hence, we may refer to these parameters as the controller inputs. The command \verb"UniversalUnfoldingPars" can be used in order to find all alternative lists of parameters, where parameters from each list play the role of the universal unfolding parameters; also see \cite[section 6]{GazorSadri} and \cite{GazorKazemiUser}. This capability is very helpful in a possible alternative modeling refinement of a real life problem. Excluding the extra parameters, the universal unfolding model is obtained. The transition set splits the parameter space into a finite number of connected components. Thereby, the controller inputs should be chosen from the interior part of each connected component (distanced from the component's boundary). These types of choices for the controller inputs prevent that the qualitative dynamics of the controlled system would be influenced by small perturbations, \ie the controlled problem is robust against small modeling imperfections.

Using \verb"RecognitionProblem"(\(x^3-\lambda x, [x, \lambda]\), \verb"UniversalUnfolding", \verb"Fractional", \verb"subs"), we obtain
\bes
\det\left( {\begin{array}{cccc}
0 & -1 & 0 & 6\\
0 & 0 & -1 & 0\\
G_{\alpha_1}(0) & G_{\lambda, \alpha_1}(0)& G_{x, \alpha_1}(0)&G_{x, x, \alpha_1}(0)\\
G_{\alpha_2}(0) & G_{\lambda, \alpha_2}(0)&G_{x, \alpha_2}(0) &G_{x, x, \alpha_2}(0)
\end{array}}
\right)\neq 0.
\ees We follow \cite[Page 7]{GolubitskySchaefferBook} to consider a small vertical force \(\epsilon\) applied on the point \(A\) and assuming a zero torque at \(x=\delta\), \ie the spring is at rest when \(x:=\delta.\) Hence, the potential function is
\bes
V(x, \lambda, \epsilon, \delta)=\frac{(x-\delta)^2}{2}+2\lambda(\cos x-1)+\epsilon \sin x,
\ees
with equilibrium equation  \(G(x, \lambda, \epsilon, \delta):= x-\delta-2\lambda \sin x+\epsilon \cos x=0\).

The command
{\tt CheckUniversal}\((G, [\delta, \epsilon], [x, \lambda], 3\), {\tt 'SingularPoint'}='\([0, \frac{1}{2}\)]') returns "yes", namely \(G\) plays the role of universal unfolding for the singular germ \(g\).

{\tt TransitionSet}\((G, [\delta, \epsilon], [x, \lambda]\), {\tt 'SingularPoint'}='\([0, \frac{1}{2}]\)') returns
\bes
\mathscr{B}:= \{(\epsilon, \delta)\,|\,\delta-\epsilon=0 \}, \quad \mathscr{H}:=\{(\epsilon, \delta)\,|\, \epsilon^3-6\delta+6\epsilon=0\}, \quad \mathscr{D}:=\emptyset,\ees
and \(\Sigma:= \mathscr{B}\cup \mathscr{H},\) while {\tt TransitionSet}\((G, [\delta, \epsilon], [x, \lambda]\), {\tt 'SingularPoint'}='\([0, \frac{1}{2}]\)', {\tt plot}) gives figure \ref{transition_set_1}. Finally, the persistent bifurcation diagrams \ref{first1}-\ref{first4} are derived via

{\tt PersistentDiagram}(\(G\), \([\delta, \epsilon]\), \([x, \lambda]\), '{\tt values}'='\([\delta=0.2346, \epsilon=0.2345]\)', '{\tt IntervalPlot}'='\([x = -2 .. 2, \lambda = 0.4 .. 0.7]\)') and the same command using the values \((\delta, \epsilon)= (-0.2346, -0.2345),\) \((0.01, -0.01),\) and \((-0.01, 0.01),\) respectively.

\section{APPENDIX}\label{Append}

\subsection{Parameter spaces \(\mathscr{H}\) and \(\mathscr{D}\) for \(G_1\) given in Equation \eqref{G1sec} }\label{HD}

\(\mathscr{H}:=\lbrace (\alpha_1, \alpha_2, \alpha_3, \alpha_4)\,|\, 0=3375\alpha_2^{3}\alpha_{3}^{2}\alpha_{4}^{3}+10125\alpha_{1}\alpha_{2}^{2}\alpha_{4}^{4}+
16875\alpha_{2}^{3}\alpha_{3}^{4}-675\alpha_{2}^{2}\alpha_{3}^{3}\alpha_{4}^{2}+288\alpha_{2}
\alpha_{3}^{2}\alpha_{4}^{4}+67500\alpha_{1}\alpha_{2}^{2}\alpha_{3}^{2}\alpha_{4}-900\alpha_{1}
\alpha_{2}\alpha_{3}\alpha_{4}^{3}+864\alpha_{1}\alpha_{4}^{5}+1080\alpha_{2}\alpha_{3}^{4}
\alpha_{4}-32\alpha_{3}^{3}\alpha_{4}^{3}+45000\alpha_{1}^{2}\alpha_{2}\alpha_{4}^{2}+
13500\alpha_{1}\alpha_{2}\alpha_{3}^{3}+3300\alpha_{1}\alpha_{3}^{2}\alpha_{4}^{2}
-108\alpha_{3}^{5}+30000\alpha_{1}^{2}\alpha_{3}\alpha_{4}+50000\alpha_{1}^{3}\rbrace\)
and
\(\D:=
\lbrace (\alpha_1, \alpha_2, \alpha_3, \alpha_4)\,|\,2000\alpha_{2}^{3}\alpha_{4}^{6}+111000\alpha_{2}^{3}\alpha_{3}^{2}\alpha_{4}^{3}
-6400\alpha_{2}^{2}\alpha_{3}\alpha_{4}^{5}+64\alpha_{2}\alpha_{4}^{7}+28000\alpha_{1}
\alpha_{2}^{2}\alpha_{4}^{4}-16875\alpha_{2}^{3}\alpha_{3}^{4}-43200\alpha_{2}^{2}
\alpha_{3}^{3}\alpha_{4}^{2}+5472\alpha_{2}\alpha_{3}^{2}\alpha_{4}^{4}-128\alpha_{3}\alpha
_{4}^{6}+45000\alpha_{1}\alpha_{2}^{2}\alpha_{3}^{2}\alpha_{4}-69600\alpha_{1}
\alpha_{2}\alpha_{3}\alpha_{4}^{3}+256\alpha_{1}\alpha_{4}^{5}+34020\alpha_{2}\alpha_{3}^{4}
\alpha_{4}-1728\alpha_{3}^{3}\alpha_{4}^{3}+130000\alpha_{1}^{2}\alpha_{2}\alpha_{4}^{2}-81000
\alpha_{1}\alpha_{2}\alpha_{3}^{3}+43200\alpha_{1}\alpha_{3}^{2}\alpha_{4}^{2}-5832\alpha_{3}^{5}
-180000\alpha_{1}^{2}\alpha_{3}\alpha_{4}+200000\alpha_{1}^{3}=0\rbrace.\)
\end{document}